\documentclass[a4paper,french,twoside,11pt]{smfartVS}

\usepackage{smfenum,smfthm,amssymb,mathrsfs,euscript}
\usepackage[latin1]{inputenc}
\usepackage{comment}
\usepackage{url}
\usepackage[usenames]{color}
\numberwithin{equation}{section}
\theoremstyle{plain}
\newtheorem*{thintro}{Théorème}

\def\FF{\mathbf{F}}
\def\NN{\mathbf{N}}

\def\RR{\mathrm{R}}
\def\ZZ{\mathbf{Z}} 

\def\A{{\rm A}}
\def\B{{\rm B}}
\def\C{{\rm C}}

\def\E{{\rm E}}
\def\F{{\rm F}}
\def\G{{\rm G}}
\def\H{{\rm H}}
\def\I{{\rm I}}

\def\K{{\rm K}}
\def\L{{\rm L}}
\def\M{{\rm M}}
\def\N{{\rm N}}
\def\P{{\rm P}}
\def\Q{{\rm Q}}
\def\R{{\rm R}}
\def\SS{{\rm S}}
\def\T{{\rm T}}
\def\U{{\rm U}}
\def\V{{\rm V}}
\def\W{{\rm W}}
\def\X{{\rm X}}
\def\Y{{\rm Y}}
\def\Z{{\rm Z}}

\def\Aa{\EuScript{A}}
\def\Bb{\mathscr{B}}
\def\Cc{\EuScript{C}}
\def\Dd{\mathscr{D}}
\def\Ee{\mathscr{E}}
\def\Ff{\mathscr{F}}

\def\Hh{\EuScript{H}}

\def\Ll{\mathscr{L}}
\def\Mm{\mathscr{M}}
\def\Nn{\mathscr{N}}
\def\Oo{\EuScript{O}}
\def\Pp{\mathscr{P}}
\def\Rr{\mathscr{R}}

\def\Tt{\mathscr{T}}

\def\Ww{\mathscr{W}}

\def\Ga{\Gamma}

\def\a{\alpha}

\def\e{\varepsilon}

\def\g{\gamma}
\def\h{\varphi}
\def\k{k}
\def\l{\lambda}

\def\p{\mathfrak{p}}
\def\s{\sigma}
\def\t{\theta}
\def\v{\upsilon}
\def\w{\varpi}

\def\ie{c'est-à-dire }

\def\>{\geqslant}
\def\<{\leqslant}

\def\Hom{{\rm Hom}}
\def\End{{\rm End}}

\def\GL{{\rm GL}}

\def\Ker{{\rm Ker}}
\def\Im{{\rm Im}}

\def\id{{\rm id}}

\def\ind{{\rm ind}}

\def\mult#1{{#1}^{\times}}
\def\red#1{\textcolor{red}{#1}}
\def\blue#1{\textcolor{blue}{#1}}
\def\mag#1{\textcolor{magenta}{#1}}

\def\1{{\bf 1}}
\def\AA{\EuScript{A}}
\def\BB{\EuScript{B}}
\def\FF{\EuScript{F}}

\def\TT{\EuScript{T}}
\def\Mod{\Mm}

\def\kb{\overline{k}}

\def\ii{{\bf I}}
\def\jj{{\bf J}}
\def\rr{{\bf R}}

\def\ip{\boldsymbol{i}}
\def\jp{\boldsymbol{j}}

\def\mm{\mathfrak{m}}
\def\nn{\mathfrak{n}}

\def\HP{\H'}
\def\CP{\C'}
\def\ModP{\Mod'}
\def\eP{\e_1}

\def\1{{\bf 1}}

\def\UB{{\U^{-}}}

\def\Dd{\EuScript{D}}

\def\Root{{\rm\Phi}}

\def\coroot{\alpha}

\def\root{\check\coroot}

\def\Cf{\C}

\def\CcI{\Cc'}
\def\HhI{\Hh'}

\def\t{\boldsymbol{\tau}}

\def\Bp{\B(\F)}
\def\Ip{\I_{1}}
\def\Gp{\tilde{\G}}
\def\Iw{\I}

\def\Xb{\X}
\def\Yb{\Y}
\def\Zb{\Z}
\def\Tb{\Omega}
\def\Sb{\SS}
\def\VV{\V}
\def\GB{{\rm G}}
\def\TB{{\rm T}}
\def\BB{{\rm B}}
\def\UB{{\rm U}}
\def\fpb{\overline{\mathbf{F}}_p}

\def\fp{\mathbf{F}_p}

\def\To{\T}
\def\cs{{\chi^{s}}}
\def\e{\boldsymbol{\varepsilon}}

\def\kb{\fpb}
\def\GB{\GL_{n}}

\def\Gp{{\tilde\G}}
\def\Gf{{\G}}
\def\Up{{\tilde\U}}
\def\Uf{{\U}}
\def\Tp{{\tilde\T}}
\def\Tf{{\T}}
\def\Bp{{\tilde\B}}
\def\Wp{{\tilde\W}}

\def\L{{\rm\Lambda}}

\def\Zp{Z}
\def\ori{\sigma_{0}}

\def\RR{\fpb}

\def\Up{\tilde{\U}}
\def\Tp{\tilde{\T}}

\author{Rachel Ollivier}
\address{Columbia University, 2990 Broadway, New York, NY 10027 \textbf{\&}}

\address{Laboratoire de Mathématiques de Versailles, 
Université de Versailles Saint-Quentin,
45 avenue des Etats-Unis,
78035 Versailles Cedex}
\email{ollivier@math.columbia.edu}

\author{Vincent S\'echerre}
\address{Laboratoire de Mathématiques de Versailles, 
Université de Versailles Saint-Quentin,
45 avenue des Etats-Unis,
78035 Versailles Cedex}
\email{vincent.secherre@math.uvsq.fr} 

\title{Modules universels de $\GL_{3}$ sur un corps $p$-adique en 
  caractéristique $p$} 

\begin{altabstract} 
Let $\F$ be a $p$-adic field with residue class field $k$.
We investigate the structure of certain mod $p$ universal modules for 
$\GL_{3}(\F)$ over the corresponding Hecke algebras.

To this  end, we first study the  structure of some mod  $p$ universal modules
for the finite group $\GL_{n}(k)$ as modules over the corresponding Hecke 
algebras.  We then  relate  this finite  case  to the  $p$-adic  one by  using
homological coefficient systems on the the affine Bruhat-Tits building of $\GL_{3}$.

Suppose that $k$ has cardinality $p$. We prove 
that the mod $p$  universal module of $\GL_{3}(\F)$ relative to the Iwahori subroup $\I$  
is flat and projective over the Iwahori-Hecke algebra: this universal module is the
 space of compactly supported functions on $\GL_3(\F)$ with values in 
$\fpb$ and invariant by translations by $\I$.  
When replacing the Iwahori subgroup of $\GL_3(\F)$ by its pro-$p$-radical, 
we prove that the corresponding module is flat over the pro-$p$ Iwahori-Hecke 
algebra if and only if $p=2$. 
\end{altabstract}

\thanks{Le second auteur est financé en partie par l'Agence Nationale de la Recherche 
(ANR-08-BLAN-0259-01).} 

\begin{document}

\maketitle
\tableofcontents
\setcounter{tocdepth}{2} 

\section*{Introduction}

Désignons par $\F$ un corps localement compact 
non archimédien de caractéristique 
résiduelle $p$ et de corps résiduel $k$. 
Soit $n\geq 1$ un entier.
Les représentations considérées sont à c\oe fficients dans une clôture 
algébrique $\fpb$ de $k$. 
Si ce sont des représentations de ${\rm GL}_n(\F)$, elles sont de plus 
lisses.  
La question de la platitude du module universel relatif à un sous-groupe 
compact maximal de 
${\rm GL}_n(\F)$ sur l'algèbre de Hecke sphérique 
a été fructueusement abordée (\cite{BO}, \cite{GK}).  
Nous considérons dans cet article 
le module universel relatif au pro-$p$-sous-groupe d'Iwahori standard $\I_1$ 
(respectivement au sous-groupe d'Iwahori standard $\I$) de ${\rm GL}_n(\F)$ en 
caractéristique $p$. Le principal résultat obtenu concerne le cas $n=3$. 

Associer à une représentation de ${\rm GL}_n(\F)$ son espace des invariants 
sous l'action de $\I_1$ définit un foncteur naturel de la catégorie des 
représentations vers celle des modules (à droite) sur l'algèbre de Hecke $\Hh$ 
du pro-$p$-sous-groupe d'Iwahori de ${\rm GL}_n(\F)$.
Notons $\Cc$ l'induite compacte du caractère trivial du pro-$p$-sous-groupe 
d'Iwahori de ${\rm GL}_n(\F)$~: c'est une représentation de ${\rm GL}_n(\F)$
et un $\Hh$-module à gauche que l'on appellera \emph{module universel}. 
\emph{Via} le produit tensoriel d'un $\Hh$-module à droite 
par $\Cc$ au dessus de
$\Hh$, on définit un adjoint à gauche $\Tt$ du précédent foncteur des 
pro-$p$-invariants.  
L'objet de cet article est d'étudier l'exactitude de cet adjoint, autrement 
dit la platitude du module universel $\Cc$ sur l'algèbre de Hecke $\Hh$.  

Considérons d'abord le cas $n=2$.
Il est prouvé  dans \cite{Platitude} que le foncteur  $\Tt$ est exact (lorsque
l'on se restreint  à la catégories des représentations  et modules avec action
scalaire du centre de ${\rm GL}_2(\F)$) si et seulement si le corps résiduel 
de $\F$ est de cardinal égal à $p$.  
Quand ce foncteur n'est pas exact, il est donc vain
d'espérer que le foncteur des pro-$p$-invariants fournisse une équivalence 
entre les $\Hh$-modules et les représentations engendrées par leurs 
pro-$p$-invariants (avec action scalaire du centre de ${\rm GL}_2(\F)$). 
Et parmi les candidats restants, $\mathbb Q_p$ est l'unique corps 
$\F$ pour lequel cette
équivalence de  catégories est réalisée  (voir \cite{Inv}, et un  argument non
publié de Paskunas 
pour le cas des extensions totalement ramifiées de $\mathbb Q_p$). 
En 2004, Vignéras a établi la classification des modules simples de l'algèbre
de Hecke du pro-$p$-Iwahori de ${\rm GL}_2(\F)$ (\cite{VigGL2}), 
classification qui ne fait pas intervenir la nature du corps $\F$, 
\ie qui ne dépend pas de la caractéristique ($0$ ou $p$) de $\F$. 
Associée à
l'équivalence de catégorie suscitée, cette classification permet de retrouver
celle des représentations irréductibles  de ${\rm GL}_2(\mathbb Q_p)$ obtenue
par Barthel-Livné et Breuil (\cite{BL}, \cite{B}). 
 
 \medskip
 
Passons maintenant au cas général de ${\rm GL}_n(\F)$, $n\>2$.
La  combinatoire des  modules sur  l'algèbre de  Hecke du  pro-$p$-Iwahori de
${\rm GL}_n(\F)$ est désormais bien comprise 
(\cite{VigProp}, \cite{OGL3}, \cite{OParab}). 
Comme dans le cas $n=2$, 
la classification des modules simples met en évidence une 
coïncidence numérique  dont les conséquences du côté  des représentations de
${\rm GL}_n(\F)$, $n\>3$ sont encore mystérieuses~: les classes
d'isomorphismes    de    modules    simples    de   dimension    $n$    dits
\emph{supersinguliers}  sont en bijection  avec les  classes d'isomorphismes
des représentations irréductibles de dimension $n$ du groupe de Galois
absolu de $\F$ (\cite{OParab}).  

\medskip

Outre le module universel $\Cc$,
définissons  l'induite compacte  $\Cc'$  du caractère  trivial du  sous-groupe
d'Iwahori de ${\rm GL}_n(\F)$ et  l'algèbre de Hecke $\Hh'$ des endomorphismes
de $\Cc'$ comme représentation de ${\rm GL}_n(\F)$.   
Pour étudier  la platitude du  $\Hh$-module $\Cc$ et du  $\Hh'$-module $\Cc'$,
nous proposons une approche qui souligne que les critères obtenus proviennent 
 des propriétés des objets  analogues dans le cas \emph{fini} que nous définissons ci-après.
Suivant cette méthode, le résultat principal obtenu est le suivant~:
 \begin{thintro} 
Supposons que le corps résiduel de $\F$ est de cardinal $p$.
\begin{enumerate}
\item 
L'induite compacte $\Cc'$ du caractère trivial du sous-groupe d'Iwahori de 
${\rm GL}_3(\F)$ est un module plat et même projectif sur l'algèbre de 
Hecke-Iwahori $\Hh'$. 
\item 
L'induite compacte $\Cc$ du caractère trivial du pro-$p$-sous-groupe d'Iwahori
de ${\rm GL}_3(\F)$ est un module plat sur la pro-$p$-algèbre de Hecke $\Hh$
si et seulement si $p=2$ auquel cas il est même projectif.  
\end{enumerate}
\end{thintro}
 
Remarquons que lorsque l'on travaille avec des représentations complexes, il 
est connu (\cite{Borel}) que $\Cc'$ est un module plat sur $\Hh'$, et ce quels 
que soit $n$ et le cardinal de $k$, 
et le foncteur des $\I$-invariants fournit une équivalence 
entre la catégorie des représentations de ${\rm GL}_n(\F)$ engendrées par leurs 
vecteurs $\I$-invariants et les $\Hh$-modules. 
 
 \medskip
 
Décrivons maintenant la méthode suivie pour démontrer le théorème.
Désignons par $\Oo$ l'anneau des entiers de $\F$, par $\varpi$ une 
uniformisante, et par $q$ le cardinal du corps résiduel $k$ où $q$ est une 
puissance de $p$. 
Soit $\K$ le sous-groupe compact maximal ${\rm GL}_n(\Oo)$ de
 ${\rm GL}_n(\F)$. Par réduction modulo $\varpi$, on dispose d'un morphisme 
surjectif de $\K$ dans le groupe réductif fini 
${\rm GL}_n(k)$. Soit $\B$ le sous-groupe de Borel de
 ${\rm GL}_n(k)$ des matrices triangulaires supérieures et
 $\U$ son radical unipotent.  
 Les images inverses de $\B$ et $\U$ dans $\K$ par la réduction modulo $\varpi$
 sont respectivement le sous groupe d'Iwahori $\I$ de ${\rm GL}_n(\F)$ 
et son pro-$p$ sous-groupe de Sylow $\I_1$.
On note respectivement $\C$ et $\CP$ l'induite à ${\rm GL}_n(k)$ du caractère 
trivial 
de $\U$ et de celui de $\B$, et $\H$ et $\HP$ les algèbres de Hecke de leurs 
${\rm GL}_n(k)$-endomorphismes. 
Les parties \ref{Sec1} à \ref{CasGL3} 
de cet article étudient certaines catégories de 
représentations modulo $p$ de ${\rm GL}_n(\F)$ et explorent la dialectique 
entre représentations et modules sur l'algèbre de Hecke dans ce cas dit \emph{fini}. 

Une classification des représentations irréductibles modulo $p$ de 
${\rm GL}_n(k)$ est donnée par \cite{CL}~: toute représentation irréductible y
est explicitée comme quotient de la représentation universelle $\C$~; de plus, 
 le foncteur  qui à une représentation associe  son sous-espace $\U$-invariant
 induit une bijection entre les représentations irréductibles de 
${\rm GL}_n(k)$ et les $\H$-modules simples. 
 Restreignons ce foncteur des $\U$-invariants à la sous-catégorie pleine $\Ee$
 des représentations (de dimension finie) engendrées par leurs vecteurs
 $\U$-invariants et  notons-le alors $\FF$.  Comme dans le cas  $p$-adique évoqué plus haut, le
 produit tensoriel par  $\C$ au dessus de $\H$ permet de  définir un adjoint à
 gauche de  $\FF$ que l'on note  $\TT$. En s'appuyant  fondamentalement sur un
 résultat  de  \cite{CE}  qui  exploite  la  propriété  d'auto-injectivité  de
 l'algèbre de Hecke  finie $\H$, on démontre que $\FF$  est une équivalence de
 catégories si et seulement si $\TT$ est un foncteur exact (proposition
 \ref{DS2019}).  
 Et il s'avère que cette condition n'est vérifiée que dans les seuls cas où 
 $n=1$, ou bien $n=2$ et $q=p$ (propositions \ref{PlatGL22} et 
 \ref{OldStuff2}).  
 Les cas de $n=2$ et $n=3$ sont traités directement en travaillant sur les 
 suites de décompositions des séries principales. 
Le cas de $n\geq 4$ s'en déduit par un processus d'induction parabolique, en 
étudiant le module universel et l'algèbre de Hecke relatifs à un 
sous-groupe de Lévi standard $\M$ de ${\rm GL}_n(k)$ (corollaire 
\ref{ThomasMarvel}).  

La partie \ref{RepPadic} établit le premier pas du passage du cas fini au cas
$p$-adique.  
On note $\K_1$ le sous-groupe de congruence de $\K$ des matrices égales à 
l'identité modulo $\varpi$.  
On démontre 
que le $\Hh$-module à gauche $\Cc^{\K_1}$ des $\K_1$ invariants de $\Cc$
est plat (et même projectif) si et seulement si $\C$ est un $\H$-module plat
(proposition \ref{LockedOut}).  
Plus généralement, on déduit de l'étude de
$\H$-modules de la forme $\C^\N$ (où $\N$ est un sous-groupe de $\U$) 
 des critères de platitude 
pour les $\Hh$-modules à gauche de la forme $\Cc^{\Nn}$ (où $\Nn$ est le 
sous-groupe de $\K_1$ image réciproque de $\N$ par réduction modulo 
$\varpi$). 

Les parties \ref{sec6} et \ref{IAS}
sont respectivement consacrées à l'étude de la
platitude du $\Hh$-module universel $\Cc$ pour ${\rm GL}_2(\F)$ et 
${\rm GL}_3(\F)$. Elles reposent de façon essentielle sur un résultat de
Schneider  et Stuhler  (\cite{SSReso})  qui nous  assure  que le  $\Hh$-module
universel $\Cc$ est l'homologie de niveau $0$ du système de c\oe fficients 
$\Hh$-équivariant sur l'immeuble qui lui est naturellement associé.
Ainsi, on dispose d'une résolution de  $\Cc$ par des sommes directes de copies
de $\Hh$-modules à gauche de la forme $\Cc^{\Nn}$ dont nous savons déterminer
la platitude par la partie \ref{RepPadic}.
 
On établit dans la partie \ref{sec6} que lorsque $n=2$ le $\Hh$-module $\Cc$ est
plat si et seulement si $q=p$.  
Notons que ce résultat étend celui de 
\cite{Platitude} précédemment cité, puisqu'ici on travaille avec le module 
universel pour ${\rm GL}_2(\F)$ et non celui de ${\rm GL}_2(\F)/\varpi^\mathbb{Z}$.  
Ajoutons que lorsque $\F=\mathbb Q_p$, le résultat de Schneider et Stuhler 
associé à l'équivalence de catégories entre représentations et modules 
permet de donner une autre preuve de l'existence d'une 
\emph{présentation standard} pour les représentations admissibles de 
${\rm GL}_2(\mathbb Q_p)$, résultat prouvé pour la première fois par Colmez 
(\cite{C}) (également prouvé dans \cite{BP}, \cite{VigCrit}, \cite{Hu}, \cite{GK2}).

Dans la partie \ref{IAS}, on exploite les résultats de la partie
\ref{RepPadic} et une propriété fine 
de la combinatoire des sommets de l'immeuble affine de 
${\rm GL}_3(\F)$ exhibée dans \cite{BO} pour prouver le théorème annoncé. 

\section*{Remerciements}

Nous remercions Marc Cabanes, Florian Herzig et Shaun Stevens pour de
fructueuses discussions à propos de ce travail. 

Le second auteur remercie l'University of East Anglia pour son accueil et son 
soutien financier (EPSRC grant GR/T21714/01). 
Il remercie aussi l'Uni\-ver\-sité de la 
Méditerranée Aix-Marseille~2 et l'Institut de Mathématiques de Luminy où la 
majeure partie de ce travail a été réalisée. 

Le premier auteur remercie le Laboratoire de Mathématiques de Versailles qui 
lui a fourni un cadre de travail 
chaleureux et dynamisant durant la gestation de cet article. Un séjour au 
\emph{Radcliffe Institute for advanced Study} à l'université de Harvard au 
printemps 2010 lui a également offert une alliance de calme et de stimulation 
qui fut sensiblement propice à la finalisation de plusieurs arguments de cet 
article.  

\section{Représentations de $\GL_{n}$ sur un corps fini}
\label{Sec1}

\subsection{Catégories de représentations}
\label{CatRep}

Soit $p$ un nombre premier, soit $\k$ un corps fini de caractéristique 
$p$ et soit $\fpb$ une clôture algébrique de $\k$. 
Étant donné un entier $n\>1$, on pose $\G=\GL_{n}(k)$, et on note $\U$ 
le sous-groupe de $\G$ cons\-ti\-tué des matrices unipotentes supérieures,
qui est un $p$-sous-groupe de Sylow de $\G$. 

\medskip

Dans cette section, et ce jusqu'à la fin de la section \ref{CasGL3}, 
par \textit{re\-pré\-sen\-ta\-tion} on entend 
re\-pré\-sen\-ta\-tion de dimension finie à coef\-fi\-cients dans $\RR$.

\begin{defi}
On note $\Rr_{\RR}(\G)$, ou plus simplement $\Rr$,
la catégorie des représentations de $\G$ 
(qui sont de dimension finie et à coef\-fi\-cients dans $\RR$),
et on note $\Ee$ la sous-catégorie pleine de $\Rr$ constituée 
des représentations engendrées par leurs vecteurs $\U$-in\-va\-riants. 
\end{defi}

\'Etant données des représentations $\V_1,\V_2$ dans $\Rr$, on note 
$\Hom_{\G}(\V_1,\V_2)$ l'espace des homomorphismes de représentations 
de $\V_1$ dans $\V_2$, \ie l'espace des applications $\R$-linéaires 
$\G$-équivariantes de $\V_1$ dans $\V_2$. 
On note $\End_{\G}(\V_1)$ la $\R$-algèbre des endomorphismes de $\V_1$.

\medskip

Le résultat classique suivant est essentiel dans l'étude de la catégorie 
$\Rr$.
Si $\V$ est une repré\-sentation de $\G$, on note $\V^\U$ l'espace de 
ses vecteurs $\U$-in\-va\-riants. 
(Pour une preuve, on renvoie à \cite{Serre}, paragraphe 8.3, proposition 26.)

\begin{lemm}
\label{QuintEssence}
Pour toute représentation $\V$ dans $\Rr$, on a $\V^\U=0$ si et seulement si 
$\V=0$. 
\end{lemm}

On note $\C$
la représentation de $\G$ induite à partir du caractère trivial de $\U$. 
Il s'agit de l'espace des fonctions de $\G$ dans $\R$ qui sont 
invariantes par translations à gauche par $\U$, muni de l'action de $\G$
par translations à droite.
C'est une représentation appartenant à $\Ee$, puisqu'engendrée par 
la fonction ca\-rac\-té\-ris\-ti\-que de $\U$, notée $\1_\U$, qui est 
$\U$-invariante. 
(Plus généralement, étant donnée une partie $\X$ de $\G$, on note $\1_{\X}$ 
sa fonction carac\-téristique.)

\medskip

Si $\V$ est une représentation de $\G$ et si $f\in\Hom_\G(\C,\V)$, 
le vecteur $f(\1_\U)$ est invariant par $\U$. 
En outre, l'application $f\mapsto f(\1_\U)$ 
est un isomorphisme de $\R$-espaces vectoriels de 
$\Hom_{\G}(\C,\V)$ dans $\V^\U$ (réciprocité de Frobenius).
On a le résultat suivant. 

\begin{lemm}
\label{CritEng}
Une représentation de $\G$ appartient à $\Ee$ si et seulement s'il 
existe un entier $b\>1$ tel qu'elle soit isomorphe à un quotient de $\C^b$. 
\end{lemm}

\begin{proof}
Soit $\V$ une représentation de $\G$, et soit $b\>1$ un entier.
Par réciprocité de Frobenius, on a un isomorphisme de $\R$-espaces vectoriels~:
\begin{equation*}
\Hom_{\G}(\C^b,\V)\to(\V^\U)^b
\end{equation*}
associant à tout homomorphisme de $\C^b$ dans $\V$ une famille de $b$
vecteurs $\U$-in\-variants de $\V$, et un tel homomorphisme est surjectif 
si et seulement si la famille qui lui correspond engendre $\V$ com\-me 
représentation de $\G$. 
\end{proof}

\begin{rema}
La sous-catégorie $\Ee$ est stable par quotients dans $\Rr$.  
On verra plus loin qu'elle n'est pas toujours stable 
par sous-objets dans $\Rr$, \ie qu'une sous-re\-pré\-sentation dans $\Rr$ 
d'un objet de $\Ee$ n'appartient pas toujours à $\Ee$.
\end{rema}

Si $\V$ est dans $\Rr$, on désigne par $\V^{\dag}$ la sous-représentation 
de $\V$ engendrée par $\V^{\U}$.
Elle est dans $\Ee$, et l'espace de ses vecteurs $\U$-invariants est 
égal à $\V^{\U}$.
Le lemme sui\-vant compare les notions de noyau et d'image dans $\Rr$ et 
dans $\Ee$.

\begin{lemm}
\label{KerIm}
Soient $\V_1,\V_2$ des représentations dans $\Ee$.  
Tout homormophisme de re\-pré\-sen\-ta\-tions 
$f\in\Hom_{\G}(\V_1,\V_2)$ a un noyau et une image dans
$\Ee$, et on a~:
\begin{equation*}
\Ker_{\Ee}(f)=\Ker_{\Rr}(f)^{\dag}
\quad\text{et}\quad
\Im_{\Ee}(f)=\Im_{\Rr}(f).
\end{equation*}
\end{lemm}

\begin{proof}
Puisque $\Im_{\Rr}(f)$ est un quotient de $\V_1$ et que $\V_1$ est dans 
$\Ee$, le lemme \ref{CritEng} im\-pli\-que que $\Im_{\Rr}(f)$ est dans 
$\Ee$, ce qui donne la seconde égalité. 

Pour la première égalité, il suffit de vérifier que si $\V$ est une 
sous-représentation de $\Ker_{\Rr}(f)$ en\-gen\-drée par $\V^\U$, 
alors $\V$ est en fait une sous-représentation de 
$\Ker_{\Rr}(f)^{\dag}$.
\end{proof}

\begin{coro}
\label{Cordoba}
Soient $\V_1$ dans $\Ee$ et $\V_2$ dans $\Rr$.  
On a un isomorphisme canonique~:
\begin{equation*}
\Hom_{\Rr}(\V_1,\V_2)\simeq\Hom_{\Ee}(\V_1,(\V_2)^\dag)
\end{equation*}
de $\RR$-espaces vectoriels.
\end{coro}

\begin{proof}
Soit $f\in\Hom_{\Rr}(\V_1,\V_2)$.
On a l'inclusion $\Im_{\Rr}(f)^{\dag}\subseteq(\V_2)^{\dag}$,
ce qui entraîne l'inclusion $\Im_{\Rr}(f)\subseteq(\V_2)^{\dag}$ 
d'après le lemme \ref{KerIm}.
\end{proof}

En d'autres termes, le foncteur $\V\mapsto\V^\dag$ est adjoint à droite 
au foncteur d'inclusion de $\Ee$ dans $\Rr$.

\medskip

On note $\V\mapsto\V^{\vee}=\Hom_{\RR}(\V,\RR)$ le foncteur exact de
la catégorie $\Rr$ dans elle-même associant à toute représentation 
de $\G$ sa représentation con\-tra\-gré\-dien\-te.

\begin{lemm}
\label{Isher}
La représentation $\C$ est isomorphe à sa contragrédiente.
\end{lemm}

\begin{proof}
Si l'on identifie $\C$ à l'espace des fonctions
de $\U\backslash\G$ dans $\RR$, l'espace dual est en\-gen\-dré par 
les fonctions d'évaluation ${\rm ev}(x):f\mapsto f(x)$, 
pour $x\in\U\backslash\G$.
On vérifie que~:
\begin{equation*}
g\cdot{\rm ev}(x)={\rm ev}(xg^{-1}), 
\end{equation*}
pour tout $g\in\G$, de sorte que l'homomorphisme $\RR$-linéaire 
$\1_{x}\mapsto{\rm ev}(x)$, as\-so\-ciant à la fonction caractéristique 
$\1_{x}$
de la classe $x\in\U\backslash\G$ sa fonction d'éva\-lu\-a\-tion, 
est un isomorphisme de re\-pré\-sen\-ta\-tions de $\G$.
\end{proof}

On vérifie au moyen des lemmes \ref{CritEng} et \ref{Isher} 
qu'une re\-pré\-sen\-ta\-tion 
de $\Rr$ appartient à $\Ee$ si et seu\-le\-ment s'il existe un entier $b\>1$
tel que sa représentation contragrédiente soit iso\-morphe à une 
sous-représentation de $\C^b$.

\begin{defi}
\label{CatBiz}
On note $\Bb$ la plus grande sous-catégorie pleine de $\Ee$ qui soit
stable par le foncteur $\V\mapsto\V^{\vee}$.
\end{defi}

D'après ce qui précède, $\Bb$ est la sous-catégorie pleine de $\Rr$
constituée des représen\-ta\-tions qui sont à la fois quotients et
sous-re\-pré\-sen\-ta\-tions de $\C^b$ pour un entier $b\>1$ assez 
grand, ou encore, si l'on préfère, des représentations images d'un 
élément de $\End_\G(\C^b)$ pour un entier $b\>1$.

\medskip

Si $\V$ est dans $\Ee$, on désigne par $\V^{\ddag}$ le plus grand 
quotient de $\V$ appartenant à $\Bb$, qui s'identifie à 
$\V^{\vee\dag\vee}$. 
Le foncteur $\V\mapsto\V^\ddag$ est adjoint 
à gauche au foncteur d'in\-clu\-sion de $\Bb$ dans $\Ee$.

\subsection{L'algèbre de Hecke\label{HeckeFinie}}

On note~: 
\begin{equation*}
\H=\End_\G(\C)
\end{equation*}
la $\RR$-algèbre des $\G$-endomorphismes de $\C$. 
Par réciprocité de Frobenius, elle s'identifie canoni\-que\-ment à l'espace 
$\C^\U$ des fonctions de $\G$ dans $\RR$ invariantes par $\U$ par 
translations à droite et à gauche, muni du produit de convolution d'unité 
$\1_\U$.  

\medskip

Soient $\T$ le tore déployé de $\G$ composé des matrices dia\-go\-nales 
et $\B$ le sous-groupe de Borel standard composé des matrices triangulaires 
supérieures de $\G$. 
On note $\Root=\Root^+\cup\Root^-$ l'ensemble des racines (décomposé en
racines po\-si\-ti\-ves et négatives) et $\Pi$ l'ensemble des racines simples.
On note $\W_0$ le sous-groupe de $\G$ constitué des matrices de 
permutation. 
C'est un groupe de Coxeter de système générateur 
$\SS_{0}=\{s_{1},\dots,s_{n-1}\}$, où $s_{i}$ désigne, 
pour tout entier $i\in\{1,\dots,n-1\}$, 
la transposition entre $i$ et $i+1$. 
On note~:
\begin{equation*}
\ell_0:\W_0\to\ZZ_{\>0}
\end{equation*}
l'application longueur qui lui correspond. 
Le groupe $\W_{0}$ agit naturellement sur 
l'ensemble des racines 
$\Root$, et l'on désignera par $w\cdot\a$ la racine conjuguée de 
$\alpha\in\Root$ par  $w\in \W_{0}$. 
L'entier $\ell_{0}(w)$ est le nombre de racines 
positives rendues négatives par l'action de $w$,
\ie le cardinal de $(w\cdot\Phi^+)\cap\Phi^-$.
On note~:
\begin{equation*}
\W=\W_0\rtimes\T
\end{equation*}
le produit semi-direct de $\W_0$ par $\T$ (où $\W_0$ agit sur $\T$ 
par conjugaison). 
On définit une action de $\W$ sur $\Root$, 
ainsi qu'une application lon\-gueur 
$\ell$ sur $\W$, par inflation à partir de celles sur $\W_{0}$.
Les éléments de longueur nulle de $\W$ sont donc exactement les éléments 
de $\T$.
Le groupe $\W$ constitue un système de re\-pré\-sen\-tants des doubles classes 
de $\G$ modulo $\U$. 
En particulier, les fonc\-tions ca\-rac\-té\-ris\-ti\-ques~:
\begin{equation}
\label{Diplo}
\t_{w}=\1_{\U w\U},
\quad 
w\in\W,
\end{equation}
forment une base de $\H$ comme $\RR$-espace vectoriel, et on a~: 
\begin{equation}
\label{pawf}
\U w\U w'\U=\U ww'\U
\quad\Leftrightarrow\quad
\t_{w}\t_{w'}=\t_{ww'}
\quad\Leftrightarrow\quad
\ell(w)+\ell(w')=\ell(ww')
\end{equation}
pour tous $w,w'\in\W$. 
L'ensemble $\{\t_{s},\t_{t}\ |\ s\in\SS_{0},\ t\in\T\}$ constitue 
donc un système générateur de la $\RR$-algèbre $\H$.

\medskip

On note $\Mod(\H)$, ou plus simplement $\Mod$, 
la catégorie des modules à droite de type fini sur $\H$. 
Si $\mm_1,\mm_2$ sont deux modules dans $\Mod$,
on note $\Hom_{\H}(\mm_1,\mm_2)$ l'espace des 
homomorphismes de $\H$-modules à droite de $\mm_1$ 
dans $\mm_2$.
On a un foncteur~:
\begin{equation}
\label{FoncF}
\V\mapsto\Hom_{\G}(\C,\V)
\end{equation}
de $\Rr$ dans $\Mod$, admettant un adjoint à gauche défini par~:
\begin{equation}
\label{FoncT}
\mm\mapsto\mm\otimes_{\H}\C,
\end{equation}
où $\C$ est considéré à la fois comme un $\H$-module à gauche 
(de type fini) et une représentation de $\G$.
L'espace $\Hom_\G(\C,\V)$ s'identifie à l'espace
$\V^\U$ des vecteurs $\U$-in\-va\-riants de $\V$ et, si $\h$ est un
homomorphisme entre deux représentations $\V_1,\V_2$ de $\Rr$, 
il lui correspond l'ap\-pli\-cation $\H$-linéaire de $(\V_1)^\U$ dans 
$(\V_2)^{\U}$ induite par restriction.

\medskip

Dorénavant, on se cantonne à l'étude de la catégorie $\Ee$, et on 
note $\FF$ la restriction de (\ref{FoncF}) à $\Ee$, \ie le composé
de (\ref{FoncF}) avec le foncteur d'inclusion de $\Ee$ dans $\Rr$.

\begin{lemm}
\label{CVS}
Le foncteur $\FF:\Ee\to\Mod$ est fidèle.
\end{lemm}

\begin{proof}
Soient $\V_1,\V_2$ dans $\Ee$ et $\h$ un homomorphisme de $\V_1$ dans 
$\V_2$ tel que $\FF(\h)=0$, \ie tel que $\h$ soit trivial sur $(\V_1)^\U$.
Il l'est donc également sur la sous-re\-pré\-sen\-ta\-tion de $\V_1$ 
engendrée par $(\V_1)^\U$, qui est égale à $\V_1$ par hypo\-thèse. 
En conclusion, on a $\h=0$, ce dont on déduit que $\FF$ est fidèle.
\end{proof}


\begin{rema}
Soit $\V$ une représentation de $\G$. 
Par adjonction, il correspond à l'identité de $\Hom_{\H}(\V^\U,\V^\U)$ 
un homomorphisme $\V^\U\otimes_{\H}\C\to\V$ de représentations de $\G$, 
dont l'image est égale à la représentation $\V^\dag$ définie au paragraphe 
\ref{CatRep}.
\end{rema}

On note $\TT$ l'adjoint à gauche de $\FF$, de $\Mod$ dans $\Ee$, 
défini par (\ref{FoncT}).

\medskip

Concentrons-nous momentanément sur la catégorie $\Bb$.
On note $\FF_\Bb$ la res\-triction de $\FF$ à $\Bb$, 
\ie le composé de $\FF$ et du foncteur d'inclusion de 
$\Bb$ dans $\Ee$. 
Il admet un adjoint à gauche $\TT_\Bb$, qui est le composé 
de $\TT$ et du foncteur $\V\mapsto\V^\ddag$,
et qui n'est pas toujours égal à $\TT$.

\begin{rema}
Puisque la dualité $\V\mapsto\V^\vee$ préserve l'irréductibilité et 
que toute re\-pré\-sen\-ta\-tion ir\-ré\-duc\-ti\-ble de 
$\G$ est engendrée par ses vecteurs $\U$-invariants (voir le 
lemme \ref{QuintEssence}), toute 
re\-pré\-sen\-ta\-tion ir\-ré\-duc\-ti\-ble de $\G$ est 
dans $\Bb$.
\end{rema}

Le théorème suivant est dû à Cabanes et 
Enguehard (voir \cite[Theorem 1.25]{CE}).

\begin{theo}[Cabanes-Enguehard]
\label{EqDeCat}
Le foncteur $\FF_\Bb$ est une équi\-va\-len\-ce de catégories entre
$\Bb$ et $\Mod$. 
\end{theo}

D'après \cite[Proposition 6.11]{CE}, $\H$ est une algèbre de 
Frobenius.  
Elle est donc auto-injective et on a la propriété suivante 
(voir \cite[Lemma 1.26]{CE}).

\begin{enonce}{Fait}
\label{SelfInj} 
Pour tout $\H$-module à droite de type fini $\mm$, 
il existe un homo\-mor\-phis\-me in\-jectif de $\mm$ dans $\H^d$, 
pour un entier $d\>1$ convenable. 
\end{enonce}

Étant donné un $\H$-module à droite de type fini 
$\mm$, on choisit un entier $d\>1$ 
et un homomorphisme in\-jectif $\iota$ de $\mm$ dans $\H^d$. 
Ceci permet d'identifier $\mm$ à un sous-espace de $\Hom_\G(\C,\C^d)$, 
et de former l'image de l'application naturelle~:
\begin{equation*}
\mm\otimes_{\H}\C\rightarrow \H^d\otimes_{\H}\C\simeq\C^d,
\end{equation*}
qu'on note $\V=\V(\mm,\iota)$. 
On a le résultat suivant 
(voir \cite[Lemma 1.27]{CE}).

\begin{enonce}{Fait}
\label{QuasiInverse}
La représentation $\V$ est dans $\Bb$, sa classe d'isomorphisme 
ne dépend pas de $\iota$ et le $\H$-module à droite $\FF_{\Bb}(\V)$ 
est isomorphe à $\mm$.
\end{enonce}

Voici une autre conséquence de l'auto-injectivité de $\H$ qui sera utile au 
paragraphe \ref{Fairway}.

\begin{lemm}
\label{Ohio}
Soient $\nn$ et $\mm$ des $\H$-modules à gauche.
On suppose que $\nn$ est un sous-module de $\mm$ 
et un $\H$-module projectif de type fini. 
Alors $\nn$ est un facteur direct de $\mm$. 
\end{lemm}

\begin{proof}
D'après \cite[Proposition 1.6.2 (ii)]{Benson}, 
le $\H$-module $\nn$ est injectif. 
\end{proof}

Signalons que, si $\G'$ est un sous-groupe de $\G$, alors $\C^{\G'}$
est un sous-$\H$-module à gauche de $\C$.

\begin{coro}
\label{Wyoming} 
Soit $\U'$ un sous-groupe de $\U$.
Alors le $\H$-module à gauche $\C^{\U}\simeq\H$ 
est un facteur direct de $\C^{\U'}$. 
En particulier, c'est un facteur direct de $\C$.
\end{coro}

Enfin,  puisque $\H$ est 
noethérien en tant que $\H$-module à droite, on a le lemme suivant
(voir par exemple \cite[Theorem 4.38]{LMR}).

\begin{lemm}
\label{lam}
Tout $\H$-module plat de type fini est pro\-jectif.
\end{lemm}

\subsection{Critères de platitude}

Le foncteur $\FF$ est fidèle et 
essentiellement surjectif (voir le lemme \ref{CVS}, 
la seconde assertion découlant par exemple du théorème 
\ref{EqDeCat}). 
Nous établissons des con\-ditions pour que, de surcroît, il soit  plein.

\begin{prop}
\label{MarquiseDeMerteuil}
\label{DS2019}
Les assertions suivantes sont équivalentes~:
\begin{enumerate}
\item
le foncteur $\FF:\Ee\to\Mod$ est plein~;
\item
le foncteur $\TT:\Mod\to\Ee$ est exact~;
\item
les catégories $\Ee$ et $\Bb$ coïncident~;
\item
le $\H$-module à gauche $\C$ est plat~;
\item
le $\H$-module à gauche $\C$ est projectif.
\end{enumerate}
\end{prop}

\begin{proof}
$(3\Rightarrow1,2)$
On suppose que les catégories $\Ee$ et $\Bb$ coïncident.
Alors le foncteur $\FF=\FF_\Bb$ est une équivalence de 
catégories, de sorte que son adjoint $\TT$ 
est exact et $\FF$ est plein.

$(1\Rightarrow3)$
On suppose que $\FF$ est plein.
Si $\V$ est dans $\Ee$, alors il existe un homo\-mor\-phis\-me 
injectif de $\H$-modules de $\FF(\V)$ dans $\H^d$ pour un entier 
$d\>1$ (voir le fait \ref{SelfInj}).
Cet homomorphisme est de la forme $\FF(\h)$ 
pour $\h\in\Hom_\G(\V,\C^d)$, et il reste à voir que $\h$ est 
injectif.
Par hypothèse, la restriction $\FF(\h)$ de $\h$ à $\V^\U$ est
injective, de sorte que $\Ker_{\Rr}(\h)^\U$ est trivial.
On en déduit que $\Ker_{\Rr}(\h)$ est trivial
(voir le lemme \ref{QuintEssence}), 
\ie que $\h$ est injective, donc que $\V$ est dans $\Bb$.

$(2\Rightarrow3)$
On suppose que $\TT$ est exact.
D'après le fait \ref{SelfInj}, pour tout $\H$-module à droite de type 
fini $\mm$,
on a un homomorphisme injectif de $\mm$ dans $\H^d$ 
pour un entier $d\>1$.
En ap\-pli\-quant $\TT$, on obtient un homomorphisme 
injectif de représentations de $\G$ de $\TT(\mm)$ dans $\C^d$, 
ce qui prouve que $\TT$ est à valeurs dans $\Bb$.
Il s'agit donc de l'adjoint de $\FF_\Bb$, qui est une équivalence de 
catégories d'après le théorème \ref{EqDeCat}.
\'Etant donné $\V$ dans $\Ee$, on a une suite exacte dans $\Rr$~:
\begin{equation}
\label{Gyoza}
0\to\Q\to\TT\, \FF(\V)=\V^\U\otimes_\H\C\to\V\to0,
\end{equation}
où $\Q$ désigne le noyau de la projection canonique de 
$\TT\, \FF(\V)$ sur $\V$.
Puisque $\TT$ est à valeurs dans $\Bb$, le foncteur 
$\FF\, \TT$ est isomorphe au foncteur identité de $\Mod$, 
de sorte que, en appliquant $\FF$ à (\ref{Gyoza}), on obtient 
la suite exacte de $\H$-modules~:
\begin{equation*}
0\to\FF(\Q)\to\FF(\V)\to\FF(\V),
\end{equation*}
où l'homomorphisme de droite est l'identité.  
On en déduit que $\FF(\Q)$, donc $\Q$, est trivial.
Ainsi $\TT\, \FF(\V)$ est ca\-no\-ni\-que\-ment isomorphe à $\V$, 
donc $\V$ est dans $\Bb$.

$(2\Leftrightarrow4)$
Si $\C$ est plat comme $\H$-module, alors le foncteur de $\Mod$ dans 
$\Rr$ défini par (\ref{FoncT}) est exact, donc $\TT$ est exact. 
Inversement, si $\nn$ est un $\H$-module à droite et $\mm$ un 
sous-module de $\nn$, on a un homomorphisme~:
\begin{equation}
\label{MadeleineR}
\mm\otimes_\H\C\to\mathfrak{n}\otimes_\H\C
\end{equation}
dans $\Rr$. 
Si $\TT$ est exact, alors le noyau de (\ref{MadeleineR}) 
dans $\Ee$ est trivial, ce qui équivaut à dire que son noyau dans $\Rr$
est trivial (voir le lemme \ref{KerIm}). 
Le foncteur (\ref{FoncT}) est donc exact, 
\ie que $\C$ est plat comme $\H$-module à gauche.

$(4\Leftrightarrow5)$
C'est une conséquence du lemme \ref{lam}, selon lequel 
tout $\H$-module (à gauche ou à droite) est 
plat si et seulement s'il est pro\-jectif.
\end{proof}

\begin{rema}
\label{Aristophon}
On  peut  ajouter à  la  proposition  
\ref{MarquiseDeMerteuil} les assertions équivalentes
sui\-van\-tes~: 
\begin{enumerate}
\item[(6)]
la représentation $\C$ est quasi-projective de type fini \cite{Selecta}~;
\item[(7)]
la sous-catégorie $\Ee$ est stable par sous-objets dans $\Rr$.
\end{enumerate}
Rappelons que $\C$ est quasi-projective comme représentation de $\G$ si, 
pour tout homomorphisme surjectif $\h\in\Hom_{\G}(\C,\V)$,
l'homomorphisme $\FF(\h)\in\Hom_{\H}(\H,\V^\U)$ qui s'en déduit est surjectif. 
D'après le lemme \ref{QuintEssence}, la représentation $\C$ est sans 
$\C$-torsion dans la terminologie de \cite{Selecta}, \ie que, 
pour toute sous-représentation non nulle $\V$ de $\C$, 
le $\H$-module $\FF(\V)$ est non nul. 
On déduit du théorème 9 de l'appendice de \cite{Selecta} 
que les conditions (1) et (6) sont équivalentes. 

D'autre part, si les catégories $\Ee$ et $\Bb$ coïncident, alors $\Ee$
est stable par sous-objets dans $\Rr$ puisque c'est le cas
de $\Bb$.
Inversement, si $\Ee$ est stable par sous-objets dans $\Rr$, 
la conjonction de la proposition 2 et du 
théorème 4 de l'appendice de \cite{Selecta} implique la condition (1).
\end{rema}

Dès que le $\H$-module $\C$ n'est pas plat, il y a donc des représentations 
dans $\Ee$ qui ne sont pas dans $\Bb$.
On en donne un exemple à la proposition \ref{eXaMpLe}.

\subsection{Décomposition de $\C$ et de $\H$}
\label{Borel}

Rappelons que $\B$ désigne
le sous-groupe des matri\-ces triangulaires 
supérieures de $\G$ et notons $\hat\To$ le groupe 
des $\RR$-caractères du tore $\T$. 
Par transitivité de l'induction, la représentation $\C$ de 
$\G$ se décompose sous la forme~:
\begin{equation}
\label{DecompC}
\C=\bigoplus\limits_{\chi\in\hat\To}\C_{\chi},
\end{equation}
où $\C_{\chi}$
désigne la représentation de $\G$ 
induite à partir du caractère de $\B$ obtenu en composant $\chi$
avec la surjection canonique de $\B$ sur $\T$.
(Il s'agit de la représentation de $\G$ par translations à droite sur 
l'espace des fonctions $f$ de $\G$ dans $\R$ qui vérifient 
$f(tug)=\chi(t)f(g)$ pour tous $t\in\T$, $u\in\U$, $g\in\G$.)
On note $\e_{\chi}$ la projection de $\C$
sur $\C_{\chi}$ dé\-fi\-nie par (\ref{DecompC}).
La famille des $\e_{\chi}$, pour $\chi\in\hat\To$, est une 
famille d'idempotents orthogonaux de $\H$ qui décomposent 
l'unité $\1_{\U}\in\H$.

\medskip

Le résultat suivant est dû à Carter-Lusztig \cite{CL}.

\begin{prop}[\cite{CL}]
\label{CL1}
Soit $\V$ une représentation irréductible de $\G$.
\begin{enumerate}
\item 
Il existe un unique caractère $\chi$ de $\T$
tel que $\V$ soit isomorphe à un quotient de $\C_\chi$.
\item 
L'espace des vecteurs $\U$-invariants de $\V$ est de 
dimension $1$ et la représenta\-tion de $\T$ sur $\V^{\U}$
est égale à $\chi$.
\end{enumerate}
\end{prop}

Pour $w\in \W_0$ et $\chi\in\hat\To$, on note $\chi^w$ le caractère 
de $\T$ défini par $t\mapsto\chi(wtw^{-1})$.
Ceci définit une action de $\W_0$ sur $\hat\To$.
On note $\Ga$
l'ensemble des orbites de $\hat \To$ sous l'action de $\W_0$. 
Si $\g$ est une telle orbite, on note $\e_{\g}$ la somme des 
$\e_{\chi}$ pour $\chi\in\g$.  

\begin{rema} 
\label{Fourier}
\begin{enumerate}
\item 
On vérifie que $\e_{\chi}$ se décompose dans la base (\ref{Diplo})
 sous la forme~:
\begin{equation}
\label{Fourier1}
\e_{\chi}=(-1)^n\sum_{t\in\T}\chi(t)\t_{t}
\end{equation}
et qu'on a l'égalité $\e_{\chi}\t_{t}=\chi(t)^{-1}\e_{\chi}$ pour tout $t\in\T$.
\item 
Pour tous $\chi\in\hat\To$ et $s\in\SS_{0}$, 
on a $\t_{s}\e_{\cs}=\e_{\chi}\t_{s}$ et, d'après \cite[Theorem 4.4]{CL}, on a~:
\begin{equation}
\label{Fourier2}
\t_{s}^2\e_{\chi}^{}=\left\lbrace 
\begin{array}{cl}
-\t_{s}\e_{\chi}&\textrm{ si } \cs=\chi,\cr
0&\textrm{ sinon. }\cr
\end{array}\right.
\end{equation} 
\item 
L'ensemble $\{\t_{s},\e_{\chi}\ |\ s\in\SS_{0},\ \chi\in\hat\To\}$ est 
un système générateur de l'algèbre $\H$. 
\end{enumerate}
\end{rema}

Par un argument classique, étant donnés 
$\chi,\chi'\in\hat\T$, l'espace
$\Hom_\G(\C_\chi,\C_{\chi'})$ est 
non nul si et seulement si $\chi$ et $\chi'$ appartiennent à la 
même orbite sous $\W_0$.
On en déduit que les $\e_{\g}$, $\g\in\Ga$, 
forment une famille d'idempotents centraux orthogonaux de 
$\H$ qui décomposent l'unité $\1_\U\in\H$.
On pose $\H_{\g}=\H\e_{\g}$, et on note $\C_\g$ la somme directe des 
$\C_\chi$ pour $\chi\in\g$, qui est un $\H_{\g}$-module.  
L'algèbre de Hecke $\H$ se décompose en la somme directe de 
$\RR$-algèbres~:
\begin{equation*}
\label{DecompH}
\H=\bigoplus\limits_{\g\in\Ga}\H_{\g}.
\end{equation*}

\begin{prop}
\label{CetSesCopains}
Le $\H$-module $\C$ est plat si et seulement si 
le $\H_{\g}$-module $\C_{\g}$ est plat 
pour toute orbite $\g\in\Ga$.
\end{prop}

\subsection{Représentations ayant des vecteurs invariants 
par le sous-groupe de Borel} 
\label{Borel2}

On note $1$ le caractère trivial de $\T$. 
Son orbite sous $\W_0$ est réduite à un singleton. 
Si $\chi=1$, on note $\CP$ et $\HP$ plutôt que $\C_1$ et $\H_1$
pour éviter d'éventuelles confusions avec certaines notations 
des sections \ref{RepPadic}, \ref{sec6} et \ref{IAS}.
L'algèbre $\HP$ est isomorphe à l'algèbre de Hecke de $\G$ relative 
au sous-groupe de Borel $\B$, \ie à l'espace des fonctions de $\G$ 
dans $\RR$ invariantes par $\B$ par translations à droite et à gauche, 
muni du produit de convolution d'unité égale à $\1_\B$.

\medskip
 
On note $\ModP$ la sous-catégorie pleine de $\Mod$ formée 
des $\H$-modules à droite de type fini $\mm$ tels que 
$\mm\cdot\eP=\mm$.
Elle s'identifie na\-tu\-rellement à la catégorie $\Mod(\HP)$ 
des $\HP$-modules à droite de type fini.
D'après la formule (\ref{Fourier1}), si $\V$ est une représentation de $\G$, 
alors $\V^\U$ est dans $\ModP$ si et seulement si $\V^\U=\V^\B$.  

\begin{prop}
\label{GTDMS}
Soit $\V$ dans $\Bb$.
Les conditions suivantes sont équivalentes~: 
\begin{enumerate}
\item 
on a $\V^\U=\V^\B$~; 
\item
il existe $b\>1$ tel que $\V$ soit isomorphe à l'image d'un endomorphisme 
de $\CP^b$.
\end{enumerate}
\end{prop}

\begin{proof}
Si $\V$ satisfait à la condition (2), elle est isomorphe à une 
sous-représentation de $\CP^b$ pour un entier $b\>1$.  
La propriété $\V^\U=\V^\B$ suit alors de ce que $\CP^\U=\CP^\B$.  

Si $\V$ satisfait à la condition (1), alors le module $\mm=\V^\U$ est dans 
$\ModP$.  
D'après le fait \ref{QuasiInverse}, si $\iota$ est un homomorphisme 
injectif de 
$\H$-modules de $\mm$ dans $\H^d$ pour un entier $d\>1$ convenable, 
l'image de l'application naturelle de $\mm\otimes_{\H}\C$ dans $\C^d$ est 
isomorphe à $\V$.  
Puisqu'on a $\mm\cdot\eP=\mm$, la représentation $\mm\otimes_{\H}\C$ 
s'identifie à $\mm\otimes_{\HP}\CP$, et son image dans $\C^d$ est 
incluse dans $\CP^d$, ce qui prouve que $\V$ satisfait à la condition (2).  
\end{proof}

On note $\Ee'$ la sous-catégorie pleine de $\Rr$ formée des représentations 
engendrées par leurs vec\-teurs $\B$-invariants et $\Bb'$ la sous-catégorie 
pleine de $\Bb$ formée des représentations $\V$ telles que $\V^\U=\V^\B$.

\begin{prop}
\label{NonNominatus}
On suppose que $\CP$ est un $\HP$-module plat. 
Alors~: 
\begin{enumerate}
\item 
les catégories $\Bb'$ et  $\Ee'$ coïncident~;
\item 
le foncteur $\FF$ induit une équivalence entre $\Ee'$ et $\ModP$. 
\end{enumerate}
\end{prop}

\begin{proof} 
D'après le théorème \ref{EqDeCat} et la proposition \ref{GTDMS}, 
le foncteur $\FF_\Bb$ induit une équi\-valen\-ce de catégories 
entre $\Bb'$ et $\ModP$, et la secon\-de assertion découle de la 
première. 

Soit $\V$ une représentation engendrée par 
l'espace $\mm$ de ses vecteurs $\B$-invariants.
Soit $\iota$ un homomorphisme injectif 
de $\H$-modules de $\mm$ 
dans $\H^d$ pour $d\>1$ (voir le fait \ref{SelfInj}).  
Rappelons (fait \ref{QuasiInverse}) que $\V$ est isomorphe à 
l'image de l'application naturelle de $\mm\otimes_{\H}\C$ dans $\C^d$.  
Puisque $\mm$ est dans $\ModP$, celle-ci s'identifie à l'image de 
l'application de $\mm\otimes_{\HP}\CP$ dans $\CP^d$.
Puisque $\CP$ est plat sur $\HP$, 
cette dernière application est injective, \ie que la représentation $\V$ est 
iso\-morphe à $\mm\otimes_{\HP}\CP$, ainsi qu'à une 
sous-représentation de $\CP^d$.
Elle est donc dans $\Bb'$. 
\end{proof}

\section{Le cas de $\GL_2(\k)$}
\label{OurExamples}
\label{puiss}

Dans toute cette section, on suppose que $n=2$, et l'on reprend les notations 
des para\-gra\-phes \ref{Borel} et \ref{Borel2}. 
Pour tout $\chi\in\hat\T$, les facteurs irréductibles de $\C_\chi$ sont tous 
de multiplicité $1$. 
Ils sont décrits par Diamond dans \cite[Proposition 1.1]{Diamond} 
(voir aussi Jeyakumar \cite{Jeya}). 
On note $s$ l'élément non trivial de $\W_0$ 
et $q$ le cardinal de $\k$.
Une orbite $\{\chi,\chi^{s}\}\in\Ga$ est dite {\it régulière} si 
$\chi\neq\chi^s$. 

\medskip

\'Etant donnée une orbite $\g$, 
la structure de la $\RR$-algèbre $\H_\g$ est déterminée par 
exemple dans \cite[\S 4]{CL}. 
Si $\g$ est régulière, $\H_{\g}$ 
est engendrée par $\SS_\g:=\t_{s}\e_{\g}$ et $\X:=\e_{\chi}$, 
avec les relations $(\SS_\g)^2=0$ et $\SS_\g\X+\X\SS_\g=\SS_\g$.
Sinon, la torsion par $\chi$ 
permet de se ramener au cas où $\g=\{1\}$, 
et $\HP$ est engendrée par $\SS:=\t_{s}\e_{1}$ avec la 
relation $\SS^2=-\SS$.

\begin{prop}
\label{PlatGL21} 
Le $\HP$-module à gauche $\CP$ est plat.
\end{prop}

\begin{proof}
En tant que $\HP$-module à gauche, 
$\HP$ est la somme directe des modules
pro\-jec\-tifs in\-dé\-com\-po\-sa\-bles $\HP\SS=\HP\t_{s}$ et 
$\HP(\SS+\e_1)$.
Une base de $\CP$ comme $\RR$-espace vectoriel
est $\{e_{a}, a\in\k\cup\{\infty\}\}$, où $e_{\infty}$ 
est la fonction caractéristique de $\B$ et $e_{a}$, 
pour $a\in\k$, celle de~:
\begin{equation*}
\B s
\begin{pmatrix} 1&a\\ 0&1\end{pmatrix}.
\end{equation*}
L'image de $e_{a}$ par $\t _{s}$ est la somme des $e_{b}$ 
pour les  $b\in\k\cup\{\infty\}$ tels que $b\neq a$ 
(voir par exemple \cite[2.2.1]{Platitude}). 
On en déduit que l'application de $\HP\oplus(\HP\SS)^{q-1}$ 
dans $\CP$ définie par~:
\begin{equation*}
(h,(h_{a})_{a\in\k^\times})
\mapsto
he_{\infty}+\sum_{a\in\k^\times}h_{a}e_{a}
\end{equation*}
est un isomorphisme de $\HP$-modules à gauche. 
On a prouvé que $\CP$ est un $\HP$-module projectif, 
donc plat.
\end{proof}

\begin{prop}
\label{PlatGL22}
Le $\H$-module à gauche $\C$ est plat si et seulement si $q=p$.
\end{prop}

\begin{proof}
D'après la proposition \ref{CetSesCopains}, l'étude de la platitude du 
$\H$-module à gauche $\C$ se ramène à celle des $\H_{\g}$-modules 
à gauche $\C_{\g}$.
Si l'orbite $\g=\{\chi,\chi^s\}$ n'est pas régulière, on se ramène en tordant 
par $\chi$ au cas où $\g=\{1\}$, et la proposition \ref{PlatGL21} implique 
que $\C_{\g}$ est plat. 
On suppose maintenant que $\g$ est régulière.
En tant que $\H_\g$-module à droite, 
$\H_{\g}$ est la som\-me directe des
modules $\e _{\chi}\H_{\gamma}$ et  
$\e _{\cs}\H_{\g}$ et on a la suite exacte de $\H_{\g}$-modules à droite~: 
\begin{equation}
\label{s1}
0\to\t _{s}\e_{\cs}\H_{\g}\to 
\e_{\chi}\H_{\g}\to \t _{s}\e_{\chi}\H_{\g}\to0.
\end{equation}
D'après \cite[Théorème 7.1]{CL}, les représentations $\t _{s}\C_{\chi}$
et $\t_{s}\C_{\cs}$ sont ir\-ré\-duc\-ti\-bles.
Le calcul de leurs dimensions se trouve par
exemple dans \cite{Paskunas}, dont le lemme 4.9 assure de plus que la somme
de ces dimensions est égale à celle de $\C_{\chi}$ si et seulement si $q=p$.
Par conséquent, la suite~:
\begin{equation}
\label{s3}
0\to\t _{s}\C_{\cs}\to \C_{\chi}\to\t _{s}\C_{\chi}\to0
\end{equation}
de représentations de $\G$ est exacte si et seulement si $q=p$.
Si $q$ est différent de $p$, cela signifie que l'exactitude 
de (\ref{s1}) n'est pas préservée par le produit 
tensoriel par le $\H_{\g}$-module à gauche $\C_{\g}$, qui n'est donc 
pas plat.

Supposons maintenant que $q$ soit égal à $p$. 
Pour montrer que $\C_{\g}$ est plat, il suffit de montrer que, 
pour tout idéal à droite $\A\subseteq\H_{\g}$, l'homo\-morphisme 
naturel de $\A\otimes_{\H_\g}\C_{\g}$ dans $\C_{\g}$ est injectif
(voir \cite{Bki}, chapitre 1, \S2, n$^\circ$3, proposition 1).  
D'après (\ref{s1}), 
il suffit de le montrer pour les idéaux 
$\e_{\chi}\H_{\g}$ et $\t_{s}\e_{\chi}\H_{\g}$ 
et leurs analogues obtenus en substituant $\cs$ à $\chi$.
Puisque $\e_{\chi}$ et $\e_{\cs}$ sont des idem\-po\-tents orthogonaux, 
la seule vérification non triviale concerne $\t_{s}\e_{\chi}\H_{\g}$ 
et elle est assurée par l'exactitude de (\ref{s3}). 
Ceci met fin à la preuve de la proposition \ref{PlatGL22}.
\end{proof}

Dans le cas où $q\neq p$, on construit une 
représentation qui est dans $\Ee$ sans être dans $\Bb$
(voir la proposition \ref{MarquiseDeMerteuil}).

\begin{prop}
\label{eXaMpLe}
Soit $\chi$ un caractère de $\T$ d'orbite régulière, 
et soit $\K$ le noyau dans $\Rr$ de $\t_{s}:\C_{\cs}\to\C_{\chi}$.
On suppose que $q\neq p$.
Alors $\K^{\vee}$ est dans $\Ee$ sans être dans $\Bb$.
\end{prop}

\begin{proof}
On a une suite exacte dans $\Rr$~:
\begin{equation}
\label{ContreExemple}
0\to\K\to\C_{\cs}\to\t_{s}\C_{\cs}\to0
\end{equation}
et puisque $q\neq p$, le noyau $\K$ contient  strictement $\t_s\C_\chi$
d'après la preuve de la proposition \ref{PlatGL22}.
En passant aux $\U$-invariants, on a les inclusions~:
\begin{equation*}
(\t_s\C_\chi)^\U\subseteq\K^\U\subseteq\C_{\cs}^\U.
\end{equation*}
Puisque $\t_s\C_\chi$ est irréductible, le
terme de gauche est de dimension $1$ (voir la proposition \ref{CL1}), 
et l'on vérifie que celui de droite est de dimension $2$.
Si l'on avait $\K^\U=\C_{\cs}^\U$, la sous-représentation de
$\K$ engendrée par $\K^\U$ serait égale à $\C_{\cs}$.
On aurait $\K=\C_{\cs}$, ce qui contredirait le fait
que la restriction de $\t_s$ à $\C_{\cs}$ n'est pas nulle. 
Aussi $\K$ et $\t_s\C_{\chi}$ ont-elles le même espace de vecteurs 
$\U$-invariants
sans être égales, \ie que $\K$ est une sous-représentation de $\C$
mais n'est pas engendrée par ses vecteurs $\U$-invariants.
Sa contragrédiente est donc dans $\Ee$ sans être dans $\Bb$.
\end{proof}

\section{Les foncteurs paraboliques}
\label{FoncParaIRJ}

\def\VV{\V}

Dans toute cette section, $\G$ est le groupe $\GL_n(k)$ et $\M$ est un
sous-groupe de Levi de $\G$, que l'on supposera standard à partir du 
paragraphe \ref{MulhollandDrive}. 
Tous les résultats de la section \ref{Sec1}, énoncés pour $\G$, 
s'étendent na\-tu\-rellement à $\M$. 
On note $\Rr(\M)$ la catégorie des représentations (de dimension finie) de 
$\M$. 

\subsection{Définition des foncteurs paraboliques\label{LeChienDeLespace}}

Soit $\P$ un sous-groupe parabolique de $\G$ muni d'une décomposition
de Levi $\P=\M\N$.
On note $\ii_{\P}$ le foncteur d'induction parabolique de 
$\Rr(\M)$ dans $\Rr(\G)$ défini pour toute représentation $\VV$ de $\M$ par~: 
\begin{equation*}
\ii_{\P}(\VV)=\{f:\G\to\VV\ |\ f(mng)=m\cdot f(g),\ m\in\M,\ n\in\N,\ g\in\G\}
\end{equation*}
que l'on munit de l'action de $\G$ par translations à droite, et 
$\rr_{\P}$ le foncteur de restriction para\-bo\-lique de $\Rr(\G)$ dans 
$\Rr(\M)$ défini pour toute représentation $\V$ de $\G$ par~:
\begin{equation*}
\rr_{\P}(\V)=\V^{\N}=\{v\in\V\ |\ n\cdot v=v, \ n\in\N\}
\end{equation*}
que l'on munit de l'action de $\M$ par restriction.
On désigne par 
$\jj_\P$ le foncteur de Jacquet de $\Rr(\G)$ dans $\Rr(\M)$ 
défini pour toute représentation $\V$ de $\G$ par~:
\begin{equation*}
\jj_{\P}(\V)=\V_{\N}=\V/\V(\N)
\end{equation*}
(où $\V(\N)$ désigne le sous-espace de $\V$ engendré par les vecteurs de la 
forme $n\cdot v-v$, 
pour $v\in\V$ et $n\in\N$),
que l'on munit de l'action naturelle de $\M$.
Pour $g\in\G$ et $v\in\VV$, on note $[g,v]$ l'élément de 
$\ii_\P(\VV)$ de support $\P g^{-1}$ et prenant en $g^{-1}$ la 
valeur $v$.
Remarquons qu'on a $[g,v]=g\cdot[1,v]$.
Les deux résultats suivants sont classiques.

\begin{prop}
\label{adjir}
Le foncteur $\ii_\P$ est adjoint à gauche de $\rr_\P$,
et le foncteur $\jj_\P$ est adjoint à gauche de $\ii_\P$.
\end{prop}

On rappelle que $\VV^\vee$ désigne la représentation contragrédiente 
de $\VV$.

\begin{prop}
Pour toute représentation $\VV$ de $\M$, 
les représentations $\ii_\P(\VV^\vee)$ et $\ii_\P(\VV)^\vee$
sont isomorphes. 
\end{prop}

On en déduit le résultat suivant. 

\begin{coro}
\label{Adjgg}
Pour toute représentation $\V$ de $\G$, les représentations 
$\rr_\P(\V)^{\vee}$ et $\jj_\P(\V^\vee)$ de $\M$ sont isomorphes.
\end{coro}

\begin{proof}
En appliquant la proposition \ref{adjir}, on voit que le foncteur~:
\begin{equation}
\V\mapsto(\rr_\P(\V^\vee))^\vee
\end{equation}
est adjoint à gauche de $\ii_\P$.
\end{proof}

\begin{rema}
Le foncteur $\ii_\P$ est exact de $\Rr(\M)$ dans $\Rr(\G)$.
Il suffit en effet de vérifier que, si $f:\V_1\to\V_2$ est un homomorphisme 
surjectif de représentations de $\M$, alors, pour tous les 
$g\in\G$ et $v_2\in\V_2$, 
la fonction $[g,v_2]\in\ii_{\P}(\V_2)$ se relève en $[g,v_1]\in\ii_{\P}(\V_1)$, 
où $v_1\in\V_1$ est un re\-lè\-ve\-ment de $v_2$.
\end{rema}

\subsection{Un système de représentants} 
\label{MulhollandDrive}

On suppose dorénavant 
que $\M$ est un sous-groupe de Levi standard de $\G$,
et que $\P=\M\N$ est le sous-groupe parabolique standard (constitué
de matrices triangulaires supérieures par blocs) lui correspondant. 

\medskip

On note $\U_\M=\U\cap\M$ le sous-groupe unipotent maximal standard de
$\M$, et $\C_\M$
la représenta\-tion de $\M$ induite à
partir du caractère trivial de $\U_\M$. 
On note~: 
\begin{equation}
\label{DefIM}
\boldsymbol{i}_\M:\C_\M\to\C
\end{equation}
l'unique homomorphisme de représentations de $\M$ envoyant $\1_{\U_\M}$ 
(la fonction caractéristique de $\U_\M$ dans $\M$) sur $\1_{\U}$. 
Son image est contenue dans l'es\-pa\-ce des vecteurs $\N$-invariants 
de $\C$.
On note $\H_\M$ l'algèbre des $\M$-endomorphismes de $\C_\M$, qu'on 
identifie au sous-espace de ses vecteurs $\U_\M$-in\-va\-riants. 
On pose $\W_\M=\W\cap\M$.

\medskip

D'après \cite[Proposition 2.3.3]{Carter}, pour tout élément $w\in\W$, 
la classe $w\W_\M$ contient un unique élément de $\W_0$ de longueur 
mini\-ma\-le, que l'on note $d(w)$. 
On pose~:
\begin{equation}
\label{Virgo}
\Dd_{\M}=\{d(w)\ |\ w\in\W\}.
\end{equation}
C'est un système de repré\-sen\-tants de $\W/\W_{\M}$. 
On note $\Root_{\M}\subseteq\Root$ l'ensemble des racines 
associées à une réflexion de $\W_{\M}$. 
Alors un élément $d\in\W_0$ appartient à $\Dd_\M$ si et seulement si 
$d\cdot\a\in \Root^+$ 
pour tout $\a\in\Root_{\M}\cap\Root^+$, ou encore, 
si et seulement si $d\cdot\a\in \Root^-$  pour toute racine $\a\in\Root_{\M}\cap\Root^-$. 
L'ensemble $\Dd_\M$ possède la propriété suivante, démontrée dans 
\cite[Proposition 2.2]{OParab}.

\begin{lemm}
\label{kitten}
Soient $s\in \SS_0$ et $d\in\Dd_\M$.
\begin{itemize}
\item Si $\ell(sd)=\ell(d)-1$, alors $sd\in \Dd_\M$.
\item Si $\ell(sd)=\ell(d)+1$, alors  ou bien $sd\in \Dd_\M$ ou bien 
$sd\in d\W_\M$.
\end{itemize}
\end{lemm}

La propriété suivante vient de \cite[Proposition 2.3.3]{Carter}.

\begin{lemm}
Pour tous $d\in\Dd_{\M}$ et $w\in\W_\M$, on a $\U d\U w\U=\U dw\U$.
\end{lemm}

Autrement dit, on a $\t_{d}\t_{w}=\t_{dw}$ pour tous 
$d\in\Dd_\M$ et $w\in\W_\M$.
On en déduit le résultat suivant.

\begin{prop}
\label{Liberte}
L'algèbre $\H$ est un $\H_{\M}$-module à droite (respectivement à gauche) 
libre de base $\{\t_{d}\}_{d\in\Dd_{\M}}$ (respectivement 
$\{\t_{d^{-1}}\}_{d\in\Dd_{\M}}$).  
\end{prop}

On note $\U^-$ le sous-groupe des matrices unipotentes 
triangulaires inférieures de $\G$.  

\begin{lemm} 
\label{TriLemme}
Pour tout $d\in\Dd_{\M}$, on a~:
\begin{enumerate}
\item 
$d\U_{\M} d^{-1}\subseteq\U$ et $d(\U^-\cap\M)d^{-1}\subseteq\U^-$~;
\item 
$d^{-1}\U d\cap\P\subseteq\U$~;
\item 
$d^{-1}\U d\N\cap\M=\U_{\M}$.
\end{enumerate}
\end{lemm}

\begin{proof} 
La propriété (1) est une conséquence de la caractérisation de $\Dd_\M$ 
en termes de racines.  
Ensuite, d'après \cite[2.5.12]{Carter}, tout élément $u\in\U$ s'écrit de fa\c con 
unique sous la forme $u=xy$, avec $x\in\U\cap d \U d^{-1}$ et 
$y\in\U\cap d\U^- d^{-1}$.  
Si $d^{-1}ud$ appartient à $\P$, alors $d^{-1}yd\in\U^-\cap\P$,
donc $y\in\U^-$.
On en déduit que $y=1$ et que $d^{-1}ud=d^{-1}xd\in\U$, ce qui prouve (2).

Soient $n\in\N$ et $u\in\U$ tels qu'on ait $d^{-1}udn\in\M$. 
D'après (2), l'élément $d^{-1}ud$ appartient à l'inter\-sec\-tion 
$d^{-1}\U d\cap\P$ qui est incluse dans $\U$. 
On en déduit que $d^{-1}udn\in\U\cap\M=\U_{\M}$, ce qui prouve que
$d^{-1}\U d\N\cap\M\subseteq\U_{\M}$. 
L'inclusion réciproque est une conséquence de (1).
\end{proof}

Notons que $\Dd_{\M}$ est aussi un système de représentants des doubles 
classes de $\U\backslash\G/\P$.

\begin{lemm}
\label{Lemme26}
L'application $\jp_\M:\H_\M\to\H$ définie par~:
\begin{equation*}
\jp_\M(\1_{\U_\M m \U_\M})=\1_{\U m\U},
\quad m\in\M,
\end{equation*}
est un homomorphisme injectif de $\R$-algèbres, 
égal à la restriction de $\boldsymbol{i}_\M$ à $\H_{\M}$.
\end{lemm}

\begin{proof} 
Pour $w\in\W_{\M}$, la fonction caractéristique de la double classe~: 
\begin{equation*}
\U_{\M}w\U_{\M}=\coprod_{x\in(\U_{\M}\cap w^{-1}\U_{\M}w)\backslash\U_{\M}} 
\U_{\M}wx
\end{equation*}
est envoyée par $\boldsymbol{i}_\M$ sur la somme des $\1_{\U wx}$. 
Comme $\U w\U=\U w\N\U_\M$, comme $\N$ est normalisé par $w$ 
et comme $wx\in\U w$ équivaut à $wxw^{-1}\in\U_\M$, 
la double classe $\U w\U$ est la réunion disjointe des $\U wx$. 
Ainsi la restriction de $\boldsymbol{i}_\M$ à $\H_{\M}$ est 
égale à $\jp_\M$.
Enfin, comme on a $\U w\U w'\U=\U w\U_{\M} w'\U$ pour tout 
$w,w'\in\W_{\M}$, on vérifie que $\jp_{\M}$ est un 
homo\-mor\-phisme de $\RR$-algèbres.
\end{proof}

On identifiera dorénavant $\H_\M$ à son image dans $\H$.

\subsection{Calcul de $\rr_\P(\C)$}
\label{neige}

On considère l'homomorphisme de représentations de $\M$~:
\begin{equation}
\label{Fraises}
\xi_{\P}:
\H\otimes_{\H_{\M}}\C_{\M}\to\rr_\P(\C),
\quad 
h\otimes c\mapsto h*\boldsymbol{i}_\M(c),
\end{equation} 
où $*$ désigne l'action à gauche de $\H$ sur $\C$.
L'objet de ce para\-graphe est de prouver la proposition suivante. 

\begin{prop}
\label{Pelleas}
L'application $\xi_{\P}$ est un isomorphisme à la fois 
de représentations de $\M$ et de $\H$-modules à gauche. 
\end{prop}

\begin{proof}
Tout élément de $\rr_\P(\C)$ est une combinaison linéaire de fonctions de la 
forme $\1_{\U g\N}$, où $g\in\G$ peut être choisi de la forme $g=dm$, avec 
$d\in\Dd_\M$ et $m\in\M$.
On fixe $d\in\Dd_\M$, et on considère l'application de $\M$ dans 
$\U\backslash \G/\N$ définie par~:
\begin{equation*}
m\mapsto\U dm\N=\U d\N m.
\end{equation*}
Elle a pour image l'ensemble des doubles classes modulo $(\U,\N)$ qui sont 
contenues dans $\U d \P$. 
D'après le lemme \ref{TriLemme}(1),
pour tout $u\in\U_\M$, les 
éléments $m$ et $um$ ont la même image.
Inversement, si $m,m'\in\M$ ont la même image par cette application, alors, 
comme $\M$ normalise $\N$, on trouve $\U dmm'^{-1}\N=\U d\N$, et donc 
$mm'^{-1}$ appartient à $d^{-1}\U d\N\cap\M$.
D'après le lemme \ref{TriLemme}(3),
on en déduit que $\U_{\M}m=\U_{\M}m'$. 
En d'autres termes, l'application~: 
\begin{equation*}
\1_{\U_\M m}\mapsto\1_{\U dm\N}
\end{equation*}
est injective et $\M$-équivariante de $\C_\M$ dans $\C$, et son image est le
sous-espace des fonctions de $\C$ supportées dans $\U d\P$.
On voit maintenant que la réciproque de $\xi_\P$ est donnée par~:
\begin{equation}
\label{SPT}
\1_{\U dm\N}\mapsto\t_{d}\otimes\1_{\U_\M m}.
\end{equation}
Enfin, on vérifie immédiatement que $\xi_\P$ est un homomorphisme 
de représentations de $\M$ et de $\H$-modules à gauche. 
\end{proof}

\begin{rema}
\label{CetteVieilleHistoireDeChambres} 
On a vu au passage dans cette preuve (voir (\ref{SPT})) 
que la fonction carac\-téristique de $\U dm \N$ est 
égale à $\t_{d}({\bf 1}_{\U m})$.
\end{rema}

On en déduit le résultat suivant. 
Soit $\Bb(\M)$ la sous-catégorie de $\Rr(\M)$ déterminée par la 
dé\-fi\-nition \ref{CatBiz}. 

\begin{prop}
\label{IndParBiz}
Soit $\VV$ une représentation dans $\Bb(\M)$.
Alors $\ii_{\P}(\VV)$ est dans $\Bb(\G)$.
\end{prop}

\begin{proof}
C'est une conséquence de l'exactitude de $\ii_\P$.
Si $\VV$ est l'image d'un endomorphisme 
$u\in\End_\G((\C_\M)^b)$ pour un entier 
$b\>1$ convenable, alors $\ii_{\P}(\VV)$ est l'image de $\ii_\P(u)$, 
qui est dans $\End_\G(\C^b)$.
\end{proof}

Le résultat suivant est un analogue de la proposition \ref{IndParBiz}
pour les foncteurs $\rr_\P$ et $\jj_\P$.

\begin{prop}
\label{ResJacParBiz}
Soit $\V$ une représentation dans $\Bb(\G)$.
Alors les représentations 
$\rr_{\P}(\V)^{\dag}$ et $\jj_\P(\V)^{\ddag}$ 
sont dans $\Bb(\M)$ (voir le paragraphe \ref{CatRep} pour les
définitions de $\dag$ et $\ddag$).
\end{prop}

\begin{proof}
D'après le corollaire \ref{Adjgg}, 
il suffit de le prouver pour $\rr_\P$.
Soit un entier $b\>1$ tel que $\V$ se 
plonge dans $\C^b$.
Puisque $\rr_\P$ est exact à gauche, 
$\rr_\P(\V)$ se plonge dans $(\rr_\P(\C))^b$. 
D'après la proposition \ref{Pelleas},
et comme $\H$ est libre comme $\H_\M$-module
à droite, $(\rr_\P(\C))^b$ est isomorphe 
à une somme directe de copies de $\C_\M$.
Ainsi $\rr_{\P}(\V)^{\dag}$ est dans $\Bb(\M)$.
\end{proof}

\subsection{Diagrammes commutatifs}
\label{Diag1}


Puisque $\U$ se décompose sous la forme $\U_\M\cdot\N$, on a~: 
\begin{equation}
\label{Nemesis}
\V^\U=\rr_\P(\V)^{\U_\M}
\end{equation}
pour toute représentation $\V$ de $\G$, qui est une égalité de 
$\H_\M$-modules à droite. 

Par adjonction, on obtient le résultat suivant. 

\begin{prop}
\label{Middlemarch}
Soit $\mm$ un $\H_\M$-module à droite de type fini.
Il existe un unique homo\-mor\-phis\-me de re\-pré\-sen\-ta\-tions de $\G$~: 
\begin{equation}
\label{EQ1}
\mm\otimes_{\H_\M}\C\to\ii_\P(\mm\otimes_{\H_\M}\C_\M)
\end{equation}
envoyant $x\otimes\1_{\U g^{-1}}$ sur $[g,x\otimes\1_{\U_\M}]$ 
pour tous $x\in\mm$ et $g\in\G$, et c'est un iso\-mor\-phisme. 
\end{prop}

\begin{proof}
Étant donné un élément $x$ de $\mm$, on note $f_x$ l'unique homomorphisme 
de $\C$ dans $\ii_\P(\mm\otimes_{\H_\M}\C_\M)$ envoyant $\1_\U$ sur 
$[1,x\otimes\1_{\U_\M}]$, qui est bien défini puisque cette dernière 
est invariante par $\U$ dans $\ii_\P(\mm\otimes_{\H_\M}\C_\M)$. 
On vérifie que l'application $x\mapsto f_x$ est $\H_\M$-linéaire. 
Par adjonction, il lui correspond l'homomorphisme 
(\ref{EQ1}), que l'on note ${\rm\Psi_{\mm}}$. 
À partir de (\ref{Nemesis}), on a~:
\begin{equation}
\label{Slan}
\Hom_{\H_\M}(\mm,\V^\U)=\Hom_{\H_\M}(\mm,\rr_\P(\V)^{\U_\M})
\end{equation}
pour toute représentation $\V$ de $\G$. 
Par une succession d'ajonctions, le membre de gauche de (\ref{Slan}) 
est $\RR$-isomorphe à~:
\begin{equation*}
\Hom_{\H}(\mm\otimes_{\H_\M}\H,\V^\U)\simeq
\Hom_{\G}(\mm\otimes_{\H_\M}\C,\V)
\end{equation*}
et le membre de droite à~:
\begin{eqnarray*}
\Hom_{\M}(\mm\otimes_{\H_\M}\C_\M,\rr_\P(\V))\simeq
\Hom_{\G}(\ii_\P(\mm\otimes_{\H_\M}\C_\M),\V)
\end{eqnarray*}
pour toute représentation $\V$ de $\G$,
ce qui prouve que les deux membres de (\ref{EQ1}) sont iso\-mor\-phes 
en tant que représentations de $\G$.
Il suffit donc de prouver que l'homomorphisme ${\rm\Psi_{\mm}}$ 
est surjectif.
Or l'induite parabolique $\ii_\P(\mm\otimes_{\H_\M}\C_\M)$ 
est engendrée comme représentation de $\G$ 
par les fonc\-tions $[1,x\otimes\1_{\U_\M}]$, avec $x\in\mm$,
qui sont dans l'image de ${\rm\Psi_\mm}$ par construction.
\end{proof}

\begin{rema}
En particulier, lorsque $\mm$ est libre de rang $1$, on obtient un 
isomorphisme de re\-pré\-sen\-ta\-tions entre $\C$ et $\ii_\P(\C_\M)$,
qui n'est autre que l'isomorphisme naturel provenant de la transitivité 
de l'induction. 
\end{rema}


\begin{prop}
\label{Mehta}
Pour toute représentation $\V$ de $\M$,
on a un isomorphisme de $\H$-modules à droite entre $\ii_\P(\V)^\U$ 
et $\V^{\U_\M}\otimes_{\H_\M}\H$.
\end{prop}

\begin{proof}
On désigne par $\Dd_{\M}'$ l'ensemble des $d^{-1}$ avec $d\in\Dd_\M$
(voir (\ref{Virgo})).
Pour $d\in \Dd_\M'$ et $x\in \V^{\U_\M}$, on désigne par $\psi_{d,x}$ la 
fonction $\U$-invariante de $\ii_\P(\V)$ de support $\P d \U$ et de valeur $x$ 
en $d$. 
Une base de $\ii_\P(\V)^\U$ est donnée par l'ensemble des $\psi_{d,x}$ pour
$d\in \Dd'_\M$ et $x$ parcourant une base de $\V^{\U_\M}$ 
(voir par exemple \cite[I.5.6]{VigBook}, en utilisant le lemme \ref{TriLemme} (2)). 
Soient $x\in \V^{\U_\M}$, $w\in \W_\M$ et $d\in \Dd'_\M$.  
Les égalités suivantes sont vérifiées~:
\begin{equation}
\label{hummus}
\psi_{1,x}\t_d=\psi_{d,x} \end{equation}
\begin{equation}
\label{bear}
\psi_{1,x}\t_w=\psi_{1,(x\,\t_w)}.
\end{equation}
Pour la première égalité, on note d'abord que la fonction $\psi_{1,x}\tau_d$ 
est $\U$-invariante de support $\P d\U$. Pour démontrer que sa valeur en $d$ 
est $x$, il suffit de remarquer que pour $u\in \U$, on a $\P d u= \P d$ si et 
seulement si $\U d u =\U d$, ce qui est donné par le lemme \ref{TriLemme} 
(2).  La seconde s'obtient aisément grâce à la décomposition de la double 
classe $\U w \U$ décrite dans la preuve du lemme \ref{Lemme26}.  

L'égalité (\ref{bear}) assure que l'on a un morphisme de $\H_\M$-modules à droite
$\V^{\U_\M}\rightarrow \ii_\P(\V)^{\U}$ bien défini par $x\mapsto\psi_{1, x}$.
Il induit un morphisme $\H$-équivariant~:
\begin{equation}
\label{EQ8f}
\V^{\U_\M}\otimes_{\H_\M}\H \longrightarrow \ii_\P(\V)^{\U}.\end{equation}
L'égalité (\ref{hummus}) assure que (\ref{EQ8f}) est surjective. 
Par la proposition \ref{Liberte}, les  espaces en question ont même dimension.
Donc (\ref{EQ8f}) est bijective.  
\end{proof}

On déduit de la proposition \ref{Mehta} le résultat suivant.

\begin{prop}
\label{Griffin}
Pour tout $\H$-module $\mm$ de type fini, 
on a un isomorphisme de 
repré\-sen\-ta\-tions de $\M$ entre $\mm\otimes_{\H_\M}\C_\M$ 
et $\jj_\P(\mm\otimes_\H\C)$.
\end{prop}

\begin{proof}
Par une succession d'adjonctions, on a~:
\begin{equation*}
\Hom_{\H}(\mm,\ii_\P(\V)^{\U})\simeq\Hom_{\G}(\mm\otimes_{\H}\C,\ii_\P(\V))
\simeq\Hom_{\M}(\jj_{\P}(\mm\otimes_{\H}\C),\V)
\end{equation*}
pour toute représentation $\V$ de $\M$, et~:
\begin{eqnarray*}
\Hom_{\H}(\mm,\Hom_{\H_\M}(\H,\V^{\U_\M}))&\simeq&
\Hom_{\H}(\mm,\Hom_{\M}(\H\otimes_{\H_\M}\C_\M,\V))\\
&\simeq&\Hom_{\M}(\mm\otimes_{\H_\M}\C_\M,\V).
\end{eqnarray*}
On en déduit que les représentations 
$\jj_\P(\mm\otimes_{\H}\C)$ et $\mm\otimes_{\H_\M}\C_\M$
sont isomorphes.
\end{proof}

\subsection{Un cas particulier} 
\label{Rentiers}

Dans ce paragraphe, on suppose que $\M=\T$ et $\N=\U$.  
Le $\H_\T$-module à gauche $\C_\T$ est libre de rang $1$, 
et l'on identifie les catégories $\Mod(\H_\T)$ et $\Rr(\T)$.  

\begin{prop}
Soit $\mm$ un $\H$-module à droite de dimension $1$.
Alors $\mm\otimes_{\H}\C$ est irréduc\-tible si et 
seulement si elle est dans $\Bb(\G)$.  
\end{prop}

\begin{proof}
D'après la proposition \ref{Griffin}, le module 
$\jj_\P(\mm\otimes_{\H}\C)$ est de dimen\-sion $1$.
C'est encore le cas de $\rr_\P((\mm\otimes_{\H}\C)^\vee)$, 
qui engendre donc dans $(\mm\otimes_{\H}\C)^{\vee}$ une 
sous-re\-pré\-sen\-ta\-tion ir\-ré\-duc\-ti\-ble 
puisque, pour toute représentation irréductible $\V$, le 
module $\V^\U$ est de dimension $1$.
On en déduit que $(\mm\otimes_{\H}\C)^{\vee}$ est irréductible 
si et seulement si elle est dans $\Bb(\G)$.
Par passage à la contragrédiente, il en est de même pour 
$\mm\otimes_{\H}\C$.
\end{proof}

Illustrons ceci au travers des exemples 
classiques du ca\-rac\-tère \emph{trivial} et du ca\-rac\-tère 
\emph{signe} de $\H$. 
On note $w_{0}$ l'élément de $\W_{0}$ de longueur maximale, 
égale à $m=n(n-1)/2$, le nombre d'éléments de $\Root^+$.
Pour chaque élément $s\in\SS_{0}$, il existe une écriture réduite 
de $w_{0}$ commen\c cant par $s$. 
On choisit une écriture réduite $w_{0}=s_{i_{m}}\dots s_{i_{1}}$. 
Avec les notations du paragraphe \ref{Borel}, pour tout $s\in\SS_0$, on pose~: 
\begin{equation*}
\t_{s}^*=\t_{s}^{}+\sum_{\chi^s=\chi}\e_{\chi}.
\end{equation*} 
C'est un élément de $\H$ vérifiant $\t_{s}^{}\t_{s}^*=\t_{s}^*\t_{s}^{}=0$. 
La fonction caractéristique de $\G$ s'écrit alors~:
\begin{equation*}
\1_{\G}=(-1)^m\t_{s_{i_{1}}}^*\dots\t_{s_{i_{m}}}^*\e_{1}.
\end{equation*}
C'est un idempotent de $\H$, donc le $\H$-module à droite 
$\1_{\G}\cdot\H$ est un facteur direct de $\H$.
C'est de plus un espace vectoriel de di\-men\-sion $1$ qui porte le caractère
trivial de $\H$, défini par $\e_{1}\mapsto 1$ et $\t_{s}\mapsto 0$ pour tout 
$s\in\SS_{0}$. 
L'application naturelle $\1_{\G}\otimes_{\H}\C\to\C$ est
injective et $\1_{\G}\otimes_{\H}\C$ est isomorphe à la représentation 
triviale de $\G$. 

\medskip

On considère maintenant le $\H$-module à droite 
$\t_{w_{0}}\e_1\cdot\H=\t_{w_{0}}\cdot\HP$, qu'on note $\mathfrak{st}$.
D'après la remarque \ref{Fourier} (2), on a
$\t_{w_{0}}\t_{s}\e_{1}=-\t_{w_{0}}\e_{1}$ pour $s\in\SS_{0}$.
On en déduit que $\mathfrak{st}$ est de dimension 
$1$ et porte le caractère signe de $\H$, défini par $\e_{1}\mapsto 1$ et 
$\t_{s}\mapsto -1$ pour tout $s\in\SS_{0}$.
De plus, $\t_{w_{0}}\e_1$ est, au signe près, un idempotent de $\H$, 
de sorte que l'application $\G$-équivariante~:
\begin{equation*}
\mathfrak{st}\otimes_{\H}\C\to\C
\end{equation*}
est injective et $\mathfrak{st}\otimes_{\H}\C$ est dans $\Bb(\G)$. 
La remarque faite au début de ce paragraphe s'applique bien au $\H$-module 
$\mathfrak{st}$ et la représentation $\mathfrak{st}\otimes_\H\C$, qui s'identifie à 
$\mathfrak{st}\otimes_{\HP}\CP$, est irréductible. 

\medskip

On va montrer 
qu'elle est isomorphe à la représentation de Steinberg ${\rm St}$ 
de $\G$ définie comme le quotient de l'induite $\ii_\B(1)$ par la 
somme de ses sous-représentations $\ii_\P(1)$, où $\P$ décrit 
l'ensemble des sous-groupes paraboliques propres de $\G$ contenant $\B$. 
Rappelons que $\ii_\P(1)$ est engendrée par la 
fonction caractéristique de $\P$.  
Or, pour $s\in\SS_{0}$, l'élément 
$\t_{s}^*\e_{1}^{}$ est la fonction caractéristique 
de $\B\cup \B s\B$, de sorte que ${\rm St}$ est isomorphe au quotient de 
$\CP$ par la somme des images $\t_{s}^*(\CP)$ pour 
$s$ parcourant $\SS_{0}$. 
Puisque $\t_{w_{0}}\t_{s}^*=0$ pour tout $s\in\SS_{0}$, on a une surjection
$\G$-équivariante~: 
\begin{equation}
\label{DinerA4HeuresDuMatin}
{\rm St}\to\mathfrak{st}\otimes_{\HP}\CP.
\end{equation} 
Il reste à voir que cet homomorphisme est injectif. 

\begin{lemm}
Pour tout $j\in\{1,\dots,m\}$, l'idempotent central $\e_1$ 
appartient à~: 
\begin{equation*}
(-1)^j\t_{s_{i_{j}}\dots s_{i_{1}}}^{}\e_1+\t_{s_{i_{j}}}^*\e_{1}
\H+\dots+\t_{s_{i_{1}}}^*\e_{1}\H.
\end{equation*}
\end{lemm}

\begin{proof} 
Pour $j=1$, on a en effet $\e_1(\t_{s_{i_{1}}}^*-\t_{s_{i_{1}}}^{})=\e_1$ d'après 
(\ref{Fourier2}). 
Supposons le lemme vrai au rang $j$, avec $1\<j\<m-1$. 
On écrit~:
\begin{equation*}
(-1)^j\t_{s_{i_{j}}\dots s_{i_{1}}}^{}\e_1
=(-1)^j(\t_{s_{i_{j+1}}}^*-\t_{s_{i_{j+1}}}^{})\t_{s_{i_{j}}\dots s_{i_{1}}}^{}\e_1.
\end{equation*}
Le résultat s'ensuit par récurrence.
\end{proof}

Le lemme au rang $j=m$ assure que pour $f\in\CP$, l'égalité
$\t_{w_{0}}f=0$ implique que
$f$ appartient à $\t_{s_{i_{m}}}^*\CP+\dots+\t_{s_{i_{1}}}^*\CP$. 
Autrement dit, l'application (\ref{DinerA4HeuresDuMatin}) est injective. 

\medskip

Remarquons que 
Cabanes et Enguehard définissent la représentation de Steinberg comme la 
représentation irréductible de $\G$ correspondant au caractère signe de $\H$ 
(\cite[Definition 6.13]{CE}) et proposent de retrouver la présente définition 
au travers de \cite[Chap. 6, Exercices 3 et 4]{CE}. 

\subsection{Condition nécessaire de platitude pour le $\H$-module $\C$}
\label{Abramovic}

On donne une condition reliant la platitude de $\C$ à celle de $\C_\M$.

\begin{prop}\label{ProspectPark}
Si $\C$ est un $\H$-module plat, alors $\C_{\M}$ est un 
$\H_{\M}$-module plat.
\end{prop}

\begin{proof}
Soit $\mm$ un idéal à droite de $\H_{\M}$. 
Notons $\K$ le noyau de l'application naturelle~: 
\begin{equation}
\label{one}
\mm\otimes_{\H_{\M}}\C_{\M}\to\C_{\M}.
\end{equation}
Le foncteur $\ii_\P$ est un foncteur exact (à gauche) de $\Rr(\M)$ dans 
$\Rr(\G)$, si bien que le noyau de~: 
\begin{equation}\label{two} 
\ii_{\P}(\mm\otimes_{\H_{\M}}\C_{\M})\to\ii_{\P}(\C_\M)
\end{equation} 
est isomorphe à $\ii_\P(\K)$. 
Les isomorphismes fournis par la proposition \ref{Middlemarch} assurent 
alors que le noyau de l'application naturelle~: 
\begin{equation}
\label{three} 
\mm\otimes_{\H_{\M}}\C\to\C
\end{equation} 
est lui aussi isomorphe à $\ii_\P(\K)$.
La proposition \ref{Liberte} dit 
que le $\H_{\M}$-module à gauche $\H$ est libre, de sorte que 
$\mm\otimes_{\H_{\M}}\H$ est isomorphe à l'idéal à droite de $\H$ engendré par 
$\mm$. 
Par conséquent,
si $\C$ est un $\H$-module plat, (\ref{three}) est injective, $\I_{\P}(\K)$ 
est la représentation nulle, et le noyau $\K$ de (\ref{one}) est trivial. 
Nous 
avons prouvé que si $\C$ est un $\H$-module plat, alors l'application 
(\ref{one}) est injective pour tout idéal à droite $\mm$ de $\H_{\M}$, 
c'est-à-dire que $\C_{\M}$ est un $\H_{\M}$-module plat. 
\end{proof}

On note $\C^{(n)}$ et $\H^{(n)}$ les quantités $\C$ et $\H$ associées à 
$\G=\GL_n(\k)$ pour $n\>1$.  

\begin{coro}
\label{ThomasMarvel}
S'il existe un entier 
$n_0\>1$ tel que le $\H^{(n_0)}$-module $\C^{(n_0)}$ ne soit pas 
plat, alors, pour tout $n\>n_0$, le $\H^{(n)}$-module $\C^{(n)}$ n'est pas 
plat. 
\end{coro}

\begin{proof}
Il suffit d'appliquer la proposition \ref{ProspectPark} avec 
$\M=\GL_{n_0}(k)\times\GL_{n-n_0}(k)$, et de remarquer que 
$\C_\M=\C^{(n_0)}\oplus\C^{(n-n_0)}$ n'est pas plat sur
$\H_\M=\H^{(n_0)}\oplus\H^{(n-n_0)}$ 
puisque $\C^{(n_0)}$ n'est pas plat sur $\H^{(n_0)}$.
\end{proof}

On déduit de ce corollaire et de la proposition \ref{PlatGL22} le 
résultat suivant. 

\begin{coro}
\label{Athos}
On suppose que $q\neq p$.
Alors, pour tout $n\>2$, le $\H^{(n)}$-module $\C^{(n)}$ n'est pas 
plat. 
\end{coro}

Voici deux résultats corollaires de la proposition \ref{Liberte}.

\begin{lemm}
\label{intents}
Tout $\H_\M$-module à gauche $\mm$ est un facteur direct de la 
restriction de $\H\otimes_{\H_\M} \mm$ à $\H_\M$.
\end{lemm}

\begin{proof} 
D'après la proposition \ref{Liberte}, l'espace vectoriel $\H\otimes_{\H_\M} \mm$ 
s'identifie à la somme directe des $ \t_{d}\mm$ pour $d\in\Dd_\M$.
On note $\nn$ la somme directe des $ \t_{d}\mm$ pour $d\in\Dd_\M$, 
$d\neq 1$.  
C'est un sous-espace vectoriel de $\H\otimes_{\H_\M}\mm$.  
Pour prouver le lemme, il suffit de s'assurer que ce sous-espace est stable 
sous l'action de $\H_\M$.  
On rappelle que l'algèbre $\H_\M$ est engendrée par 
$\{\t_{s},\,\t_{t}\ |\ s\in \SS_{0}\cap \W_\M, t\in \T\}$.

Puisqu'un élément $t\in\T$ est de longueur nulle dans $\W$, les
relations (\ref{pawf}) assurent que, pour tout $d\in\Dd_\M$, on a l'égalité~: 
\begin{equation*}
\t_t \t_d=\t_d\t_{d^{-1}t d},
\end{equation*}
et l'on remarque que ${d^{-1}td}\in\T$. 
Ainsi, un sous-espace de la forme $\t_d\mm$ avec 
$d\in\Dd_\M$ est stable sous l'action de $\t_t$ pour $t\in\T$.
On en déduit que $\nn$ est stable sous l'action de $\t_t$ pour tout $t\in \T$.

Soit $s\in \SS_{0}\cap \W_\M$. 
Vérifions que l'espace $\nn$ est stable sous l'action de $\t_{s}$.  
Soit $d\in\Dd$ tel que $d\neq 1$. 
Si $\ell(sd)= \ell(d)-1$, 
alors $\t_{d}=\t_{s}\t_{sd}$ et $sd\in\Dd_\M$ d'après le lemme \ref{kitten}.  
D'après la remarque \ref{Fourier}, on a~: 
\begin{equation*}
\t_{s}^2=-\t_{s}\sum_{\chi^s=\chi}\e_{\chi}.
\end{equation*}
Ainsi l'élément~:
\begin{equation*}
a=-\sum_{\chi^s=\chi}\e_{\chi}
\end{equation*}
vérifie $\t_{s}\t_{d}=\t_{s} a\t_{sd}$.
On en déduit que $\t_{s}\t_{d}\mm=\t_{s} a\t_{sd}\mm$, 
qui est inclus dans $ \t_s \t_{sd}\mm=\t_{d}\mm$ d'après 
l'argument précédent. 
Si $\ell(sd)=\ell(d)+1$, alors $\t _{s}\t _{d}=\t _{sd}$.  
Si l'élément $sd$ appartient à $\Dd_\M$, alors 
$\t _{s}\t _{d}\mm=\t _{sd}\mm$ est inclus dans $\nn$. 
Sinon, $sd\in d\W_\M$ d'après le lemme \ref{kitten}.  
Il existe donc $w\in \W_\M$ tel que $\t _{sd}=\t _{d}\t _{w}$, 
de sorte que $\t _{s}\t _{d}\mm=\t_{d}\t_{w}\mm$ est inclus dans 
$\t _{d}\mm\subseteq\nn$. 
Ainsi, $\nn$ est bien un $\H_\M$-module. 
\end{proof}

\begin{lemm}
\label{purposes}
Un $\H_\M$-module à gauche de type fini $\mm$ est plat si et seulement si 
le $\H$-module $\H\otimes_{\H_\M}\mm$ est plat, 
et dans ce cas, ils sont tous deux projectifs. 
\end{lemm}

\begin{proof} 
D'après le lemme \ref{lam}, un module de type fini sur $\H$ 
ou $\H_\M$ est plat si et seulement s'il est projectif, ce qui 
prouve la dernière assertion.
Par ailleurs, un module (sur $\H$ ou $\H_\M$) est 
projectif si et seulement s'il est facteur direct d'un module 
libre.  
On en déduit que si $\mm$ est un $\H_\M$-module à gauche 
projectif, alors $\H\otimes_{\H_\M}\mm$ est un 
$\H$-module à gauche projectif. 

Inversement, si le $\H$-module à gauche $\H\otimes_{\H_\M}\mm$ 
est projectif, c'est un facteur direct d'un $\H$-module libre. 
D'après le lemme \ref{intents}, le $\H_\M$-module à gauche $\mm$ 
est alors un facteur direct d'une somme de copies de $\H$, qui est 
un $\H_\M$-module à gauche libre donc projectif. 
\end{proof}

On déduit de ce qui précède la proposition suivante.  

\begin{prop}
\label{Monsters0} 
Le $\H$-module $\C^\N$ est plat si et seulement si le $\H_{\M}$-module 
$\C_{\M}$ est plat. 
Dans ce cas ils sont tous deux projectifs et le $\H$-module $\C^\N$ 
est un facteur direct de $\C$. 
\end{prop}

\begin{proof} 
Rappelons que l'homomorphisme $\ip_{\M}$ défini en (\ref{DefIM})
est un homomorphisme injectif de $\C_{\M}$ dans $\C$, 
et que, d'après la proposition \ref{Pelleas}, l'homomorphisme
(\ref{Fraises}) est un isomorphisme de $\H$-modules entre 
$\H\otimes _{\H_{\M}}\C_{\M}$ et $\C^\N$. 
Le résultat est alors une conséquence des lemmes \ref{purposes} et 
\ref{Ohio}. 
\end{proof}

\section{Le cas de $\GL_3(\k)$}
\label{CasGL3}

Dans ce paragraphe, on suppose que $n=3$.
L'objectif de cette section est de démontrer les 
propositions \ref{OldStuff1} et \ref{OldStuff2}
énoncées plus loin.

\subsection{Représentations algébriques de $\GL_n(\mathbf{F}_{p^r})$}

Dans ce paragraphe, on rappel\-le quelques résultats de \cite{Jantzen}
et \cite{Herzig}.

\medskip

Soit $\TB(\kb)$ le sous-groupe des matrices diagonales de $\GB(\kb)$, 
soit $\BB(\kb)$ le sous-groupe des matrices triangulaires supérieures 
de $\GB(\kb)$, soit $\UB(\kb)$ son radical unipotent et soit $\UB^-(\kb)$ 
le radical unipotent du sous-groupe de Borel opposé.
On note $\X$ le groupe des caractères algé\-bri\-ques de $\TB(\kb)$,
que l'on iden\-tifie à $\ZZ^n$, 
et $\X_+$ l'ensemble des $n$-uplets $(a_1,\dots,a_n)\in\X$ tels que 
$a_1\>\dots\>a_n$. 
Pour $\l\in\X_+$, on note $\Ww(\l)$ l'espace de fonc\-tions rationnelles~:
\begin{equation*}
\{f:\GL_n(\kb)\to\kb\ |\ f(gtu)=\l(t)^{-1}f(g),\ g\in\GB(\kb),\ t\in\TB(\kb),\ u\in\UB^-(\kb)\}
\end{equation*}
qu'on munit de l'action de $\GB(\kb)$ par translations à gauche. 
C'est une re\-pré\-sen\-ta\-tion algébri\-que de $\GL_n(\kb)$.
On note $\Ll(\l)$ son socle.  
Pour tout entier $r\>0$, on pose~: 
\begin{equation*}
\X_r=\{(a_1,\dots,a_n)\in\X_+\ |\ 0\<a_i-a_{i+1}\<p ^r-1,\ i\in\{1,\dots,n\}\}.
\end{equation*}
Pour tout $r\>1$, on note $\mathbf{F}_{p^r}$ l'unique sous-corps de $\fpb$ 
de cardinal $p^r$ et, étant donnés $\l\in\X_r$ et $i\in\{0,\dots,r-1\}$, on note 
$\Ll_{r}(\l)$ la restriction de $\Ll(\l)$ à $\GB(\mathbf{F}_{p^r})$ 
et $\Ll_{r}(\l)^{(i)}$ la composée de $\Ll_{r}(\l)$ avec 
l'automorphisme de $\GB(\mathbf{F}_{p^r})$ induit par 
$x\mapsto x^{p^i}$. 
On a les résultats suivants.

\begin{prop}[\cite{Jantzen}, II.3] 
\label{Jan2}
On fixe un entier $r\>1$.
\begin{enumerate}
\item 
Pour tout $\l\in\X_{r}$, la re\-pré\-sen\-ta\-tion $\Ll_{r}(\l)$ de 
$\GB(\mathbf{F}_{p^r})$ est irréductible.
\item 
L'application $\l\mapsto\Ll_{r}(\l)$ induit une bijection entre 
$\X_{r}/(p^r-1)\X_0$ et l'ensemble des classes d'isomorphisme de 
représentations irréductibles de $\GB(\mathbf{F}_{p^r})$.
\item 
L'espace des vecteurs $\UB(\mathbf{F}_{p^r})$-invariants de 
$\Ll_{r}(\l)$ est de dimension $1$, 
et la re\-pré\-sen\-ta\-tion de $\TB(\mathbf{F}_{p^r})$ sur 
cet espace est égale à $\l$. 
\item
Si $\l=(a_1,\dots,a_n)\in\X_{r}$, alors 
$\l^*=(-a_n,\dots,-a_1)\in\X_{r}$ et $\Ll_{r}(\l^*)$ 
est isomorphe à la représentation contra\-grédiente de $\Ll_{r}(\l)$.
\item 
Soit $\l\in\X_{r}$, qu'on écrit sous la forme~:
\begin{equation}
\label{PadreDinis}
\l=\l_0+\l_1p+\dots+\l_{r-1}p^{r-1},
\quad
\l_i\in\X_1,
\quad
i\in\{0,\dots,r-1\}.
\end{equation}
Alors on a un isomorphisme de représentations de $\GB(\mathbf{F}_{p^r})$~:
\begin{equation}
\label{SebastiaoDeMelo}
\Ll_{r}(\l)\simeq\Ll_{r}(\l_0)\otimes\Ll_{r}(\l_1)^{(1)}\otimes
\dots\otimes\Ll_{r}(\l_{r-1})^{(r-1)}.
\end{equation}
\end{enumerate}
\end{prop}

En comparant les propositions \ref{Jan2} et \ref{CL1}, on obtient le résultat suivant. 

\begin{coro}
\label{Restr}
Soit $\l\in\X_{r}$ et soit $\chi$ un caractère de $\TB(\mathbf{F}_{p^r})$.
La représentation $\Ll_{r}(\l)$ est un quotient ir\-ré\-ductible de $\C_\chi$ 
si et seulement si la restriction de $\l$ à $\TB(\mathbf{F}_{p^r})$ 
est égale à $\chi$.
\end{coro}

\begin{exem}
On suppose que $n$ est égal à $2$.
Soit $\l=(a,b)\in\X_{r}$, et écrivons $a-b$ sous la forme 
$e_0+e_1 p+\dots+e_{r-1}p^{r-1}$
avec $e_i\in\{0,\dots,p-1\}$.
Alors $\Ll_{r}(\l)$ est la représentation irréductible~:
\begin{equation*}
{\rm Sym}^{e_0}(\overline{\mathbf{F}}{}^2_p)\otimes\dots\otimes
{\rm Sym}^{e_{r-1}}(\overline{\mathbf{F}}{}^2_p)^{(r-1)}\otimes\det{}^{b},
\end{equation*}
qui est de dimension $(e_0+1)\dots(e_{r-1}+1)$.
\end{exem}

On suppose maintenant que $n$ est égal à $3$. 
On rappelle quelques résultats de Herzig \cite{Herzig} sur la 
semi-simplification de $\C_\chi$ pour 
$\G={\rm GL}_3(\mathbf F_p)$. 

\begin{prop}[\cite{Herzig}, Proposition 4.9]
\label{HenryJames}
Soit $(a,b,c)\in\X_1$.
Si~:
\begin{equation}
\label{RamPasquier}
0\<a-b,b-c<p-1
\quad\text{et}\quad
p-1\<a-c,
\end{equation}
alors la représentation $\Ww(a,b,c)$ est de longueur $2$.
Sinon, $\Ww(a,b,c)$ est irréductible.
\end{prop}

On commence par étudier $\C_\chi$ lorsque $\chi=1$.
Pour tout $\l\in\X_1$, on note $\Ww_{1}(\l)$ la restriction de $\Ww(\l)$ à 
${\rm GL}_3(\mathbf F_p)$, qui est de longueur $\<2$ d'après la proposition 
\ref{HenryJames}.

\begin{prop}
\label{AngelaDeLima}
Dans le groupe de Grothendieck des représentations de longueur finie de
$\GL_3(\fp)$, la semi-simplification de $\CP$ est égale à~:
\begin{equation}
\label{DecSSC1Weyl}
\Ww_{1}(0,0,0)+2\Ww_{1}(p-1,0,0)+2\Ww_{1}(p-1,p-1,0)+\Ww_{1}(2p-2,p-1,0).
\end{equation}
Chacune des représentations $\Ww_1(\l)$ apparaissant ci-dessus 
est irréductible, et on a~:
\begin{eqnarray}
\dim \Ww_{1}(0,0,0)&=&1,\\
\dim \Ww_{1}(p-1,0,0)&=&{p(p+1)}/{2},\\
\dim \Ww_{1}(p-1,p-1,0)&=&{p(p+1)}/{2},\\
\dim \Ww_{1}(2p-2,p-1,0)&=&p^3.
\end{eqnarray}
\end{prop}

\begin{proof}
Le théorème \cite[5.1]{Herzig} donne la décomposition de 
la semi-simplification de $\CP$ 
en une somme de représentations $\Ww_1(\l)$ et la proposition \ref{HenryJames} 
montre que ces représentations sont ir\-ré\-duc\-ti\-bles.
Plus précisément, ce théorème décrit la décomposition 
de la semi-simplification de la réduction modulo $p$ 
de chacun des facteurs irréductibles de l'induite du 
$\mathbf{Z}_p$-ca\-rac\-tère trivial de $\BB(\fp)$ à $\GL_3(\fp)$~:
\begin{itemize}
\item[$\bullet$]
la réduction modulo $p$ du $\mathbf{Z}_p$-caractère trivial de $\GL_3(\fp)$ est
isomorphe à la représenta\-tion irréductible $\Ww_{1}(0,0,0)$, qui est de dimension 
$1$~;
\item[$\bullet$]
la réduction modulo $p$ de la $\mathbf{Z}_p$-représentation de Steinberg 
est isomorphe à la représenta\-tion irréductible
$\Ww_{1}(2p-2,p-1,0)$, qui est de dimension $p^3$~;
\item[$\bullet$]
la $\mathbf{Z}_p$-représentation irréductible de $\GL_3(\fp)$ apparaissant avec 
multiplicité $2$ dans l'induite à $\GL_3(\fp)$ du $\mathbf{Z}_p$-ca\-rac\-tère 
trivial de $\BB(\fp)$ est de dimension $p^2+p$. 
Sa réduction modulo $p$ est isomorphe à la somme des deux représentations 
irréductibles $\Ww_{1}(p-1,0,0)$ et $\Ww_{1}(p-1,p-1,0)$.
Il suffit donc de montrer que ces deux-là ont la même dimension,
ce qui découle du fait qu'elles sont duales d'après la proposition 
\ref{Jan2}(4). 
\end{itemize}
On en déduit le résultat annoncé.
\end{proof}

On étudie maintenant $\C_\chi$ avec $\chi$ régulier, 
\ie dont l'orbite sous l'action de $\W_0$ est de cardinal $6$.

\begin{prop}
\label{BartonFink}
Soit $\chi$ un caractère régulier de $\T(\fp)$.
Alors $\C_\chi$ est de longueur strictement supérieure à 
$6$ dans la ca\-té\-gorie $\Rr(\GL_3(\fp))$.
\end{prop}

\begin{proof}
On choisit $(a,b,c)\in\X_1$ tel que~:
\begin{equation*}
\chi:
\begin{pmatrix}
x&&\\
&y&\\
&&z\\
\end{pmatrix}
\mapsto x^a y^b z^c.
\end{equation*}
Puisque $\chi$ est régulier, on a $p>2$ et l'on peut supposer que $1\<a-b,b-c<p-1$.
D'après la formule \cite[(5.2)]{Herzig}, la semi-simplification de $\C_\chi$ 
est~: 
\begin{equation*}
\begin{split}
\Ww_{1}(a,b,c)\; +&\; \Ww_{1}(p-1+b,p-1+c,a)+\Ww_{1}(p-1+c,a,b)\\
&+\Ww_{1}(2p-2+c,p-1+b,a)+\Ww_{1}(p-1+a,p-1+c,b)
+\Ww_{1}(p-1+b,a,c).
\end{split}
\end{equation*}
On vérifie que l'un des deux triplets $(a,b,c)$ et $(p-1+a,p-1+c,b)$ 
satisfait à la condition \eqref{RamPasquier} 
de la proposition \ref{HenryJames}, de sorte que 
l'une ou l'autre des représentations~:
\begin{equation*}
\Ww_{1}(a,b,c),
\quad
\Ww_{1}(p-1+a,p-1+c,b)
\end{equation*}
est de longueur $2$, ce qui prouve l'assertion. 
\end{proof}

\begin{rema}
On peut montrer de la même façon que, si l'orbite de $\chi$ sous l'action de 
$\W_0$ est de cardinal $3$, alors $\C_\chi$ est de longueur $>6$. 
Compte tenu de la proposition \ref{AngelaDeLima}, on en déduit que 
$\C_\chi$ est de longueur $6$ si et seulement $\chi$ est invariant par $\W_0$.
\end{rema}

On en déduit le résultat suivant.
Soit $r\>1$ un entier.

\begin{prop}
\label{Swift}
Soit $(a,b,c)\in\X_{r}$.
La représentation irréductible $\Ll_{r}(a,b,c)$ est isomorphe à un 
quotient de $\CP$ si et seulement si $(a,b,c)$ est congru à l'un 
des poids~: 
\begin{equation}
\label{ListePoids}
(0,0,0),
(p^r-1,0,0),
(p^r-1,p^r-1,0),
(2p^r-2,p^r-1,0),
\end{equation}
modulo $(p^r-1)\X_0$.
En outre, on a~:
\begin{eqnarray}
\dim \Ll_{r}(0,0,0)&=&1,\\
\dim \Ll_{r}(p^r-1,0,0)&=&\left({p(p+1)}/{2}\right)^{r},\\
\dim \Ll_{r}(p^r-1,p^r-1,0)&=&\left({p(p+1)}/{2}\right)^{r},\\
\dim \Ll_{r}(2 p^r-2,p^r-1,0)&=&p^{3r}.
\end{eqnarray}
\end{prop}

\begin{proof}
La première partie de la proposition est une conséquence du corollaire 
\ref{Restr}.  
Ensuite, puisque $p^r-1=(p-1)(1+p+\dots+p^{r-1})$, 
chaque poids $\l$ dans (\ref{ListePoids}) se décom\-pose sous 
la forme (\ref{PadreDinis}) avec des $\l_i\in\X_1$ indépendants de $i$ et 
respectivement égaux, suivant $\l$, à~:
\begin{equation}
\label{ListePoidsp}
(0,0,0),
(p-1,0,0),
(p-1,p-1,0),
(2p-2,p-1,0).
\end{equation}
Compte tenu de la formule (\ref{SebastiaoDeMelo}) et 
de la proposition \ref{AngelaDeLima}, on trouve
les formules annoncées. 
\end{proof}

\subsection{L'algèbre de Hecke $\HP$}
\label{PremH}

D'après \cite[\S4]{CL}, la $\RR$-algèbre $\HP$ est 
en\-gen\-drée par $\Sb_1=\t_{s_{1}}\e_{1}$ 
et $\Sb_2=\t_{s_{2}}\e_{1}$ avec les relations~: 
\begin{equation*}
\Sb_1 \Sb_{2} \Sb_{1}=\Sb_{2}\Sb_{1}\Sb_{2},
\quad
\Sb_{1}^{2}+\Sb_{1}^{}=\Sb_{2}^{2}+\Sb_{2}^{}=0.
\end{equation*}
Remarquons que $\Sb_{1}$ et $\Sb_{2}$ sont les fonctions caractéristiques 
respectives de $\B s_{1}\B$ et $\B s_{2}\B$. 
Dans la suite, on pose $\Sb_{1}^*=\Sb_{1}^{}+\e_1^{}$ et
$\Sb_{2}^*=\Sb_{2}^{}+\e_1^{}$.
On pose~: 
\begin{equation*}
\Xb=-\Sb_{1}\Sb_{2}\Sb_{1},
\quad
\Yb = -\Sb_{1}^{} \Sb_{2}^* \Sb_{1}^{},
\quad
\Zb=-\Sb_{1}^*\Sb_{2}^{}\Sb_{1}^*,
\quad
\Tb=\Sb_{1}^*\Sb_{2}^*\Sb_{1}^*.
\end{equation*}
On vérifie que ce sont des idempotents deux à deux orthogonaux de $\HP$ 
qui décomposent l'unité.  
Notons que $\Xb $ et $\Tb$ sont centraux et qu'on a les relations~:
\begin{equation*}
\Yb\Sb_{2}^{}=\Sb_{2}^*\Zb,
\quad 
\Sb_{2}^{}\Yb =\Zb\Sb_{2}^*,
\end{equation*} 
qui permettent 
en particulier de s'assurer que la somme $\Yb+\Zb $ est un idempotent central.  
Ainsi, les idéaux à droite $\Xb\HP,\Yb\HP,\Zb\HP$ et $\Tb\HP$ 
sont des $\HP$-modules projectifs indécomposables.

\medskip

On sait (voir \cite[Theorem 1.25]{CE}) 
que l'application $\mm\mapsto{\rm soc}(\mm)$ qui 
à un $\HP$-module à droite associe son socle (\ie 
son plus grand sous-module semi-simple) induit 
une bijection entre les classes d'isomorphisme de 
$\HP$-modules projectifs indécomposables et les 
classes d'isomorphisme de $\HP$-modules simples.
On va expliciter cette bijection.  

\begin{defi}
Pour  $\alpha_{1},\alpha_{2}\in \{0,-1\}\subseteq\RR$, on désigne par
$\chi_{\alpha_{1},\alpha_{2}}$  le caractère de $\HP $ défini par
$\chi_{\alpha_{1},\alpha_{2}}(\Sb_1)= \alpha _{1}$ et
$\chi_{\alpha_{1},\alpha_{2}}(\Sb_2)= \alpha _{2}$. 
\end{defi}

L'application $(\alpha_{1},\alpha_{2})\mapsto\chi_{\alpha_{1},\alpha_{2}}$ définit 
une bijection de $\{0,-1\}\times\{0,-1\}$ sur l'ensemble des classes 
d'isomorphisme de $\HP$-modules simples.

\begin{prop}
On a~: 
\begin{equation*}
{\rm soc}(\Xb\HP)=\chi_{-1,-1},
\quad
{\rm soc}(\Yb\HP)=\chi_{0,-1},
\quad
{\rm soc}(\Zb\HP)=\chi_{-1,0},
\quad
{\rm soc}(\Tb\HP)=\chi_{0,0}.
\end{equation*}
\end{prop}

\begin{proof}
On vérifie d'abord que 
les idéaux bilatères $\Xb \HP$ et $\Tb\HP$ sont de dimension $1$ 
et correspondent respectivement aux caractères $\chi_{-1, -1}$ et $\chi_{0,0}$.  
Ensuite, 
en utilisant la relation $\Sb_{1}^{}\Sb_{2}^*=-\Yb \Sb_{2}^*$, on vérifie que 
$\Yb \HP =\Sb_{1}^{}\Sb_{2}^*\HP^{}$ est un $\RR$-espace vectoriel de dimension $2$ 
de base $\{ \Yb \Sb_{2}^{},\Sb_{1}^{}\Sb_{2}^*\}$.  
On a la suite exacte non scindée de $\HP$-modules~: 
\begin{equation}
\label{h1}
0\to\Yb \Sb_{2}\HP = \Sb_{2}^*\Zb \HP^{} \to \Yb \HP
\overset{\Sb_{2}}\longrightarrow \Sb_{2}\Yb \HP\to 0,
\end{equation} 
\ie que le $\HP$-module à droite $\Yb \HP $ est l'enveloppe projective de 
$\chi_{0,-1}$, et c'est une ex\-tension non scindée de $\chi_{-1,0}$ par 
$\chi_{0, -1}$. 
Le $\HP$-module à droite $\Yb \HP $ est l'enveloppe projective de $\chi_{0,-1}$.
De même, on a la suite exacte de $\HP$-modules~: 
\begin{equation}
\label{h2}
0\to\Zb \Sb_{2}^* \HP^{} = \Sb_{2}\Yb \HP \to \Zb \HP
\overset{\Sb_{2}^*}\longrightarrow \Sb_{2}^*\Zb \HP^{} \to 0.
\end{equation} 
C'est une extension non scindée de $\chi_{0,-1}$ par $\chi_{-1, 0}$.  
Le $\HP$-module $\Zb \HP $ à droite est l'enveloppe projective de $\chi_{-1,0}$.  
\end{proof}

Ainsi les idéaux à droite de $\HP$ sont, à isomorphisme près, les 
$\H'$-modules projectifs indé\-com\-posables $\Xb \HP$, $\Yb \HP$, $\Zb \HP$ et $\Tb\HP$, 
auxquels s'ajoutent les idéaux non projectifs $\Sb_{ 2} \Yb \HP$ et $\Yb\Sb_{ 2}\HP$.

\begin{rema}
\label{Gizzi}
Notons que $\Sb_2\Yb\C'= \Sb_1^* \Sb_2\Sb_1\C'$ et 
$\Yb\Sb_2\C'= \Sb_2^*\Sb_1\Sb_2\C'$ de sorte que la 
classification de Carter et Lusztig (\cite[Theorem 7.4]{CL})
 assure que les représentations irréductibles de $\G$ 
possédant un vecteur $\B$-invariant sont, à isomorphisme près~:
\begin{equation*}
\Xb\C', \Sb_2\Yb \C', \Yb \Sb_2\C', \Omega \C'
\end{equation*}
et leurs espaces $\U$-invariants (donc $\B$-invariants)
respectifs portent les caractères $\chi_{-1, -1}$, $\chi_{-1, 0}$, $\chi_{0, -1}$, $\chi_{0, 0}$ de $\HP$.
\end{rema}

\subsection{Platitude de $\C'$} 
\label{Prem1}

On a le résultat suivant.

\begin{prop}
\label{OldStuff1}
Le $\HP$-module $\CP$ est plat si et seulement si $q=p$. 
\end{prop}

\begin{proof}
On note $q$ le cardinal de $\k$.
Rappelons que $\CP$ est de dimension~: 
\begin{equation*}
(\G:\B)=(1+q)(1+q+q^2).
\end{equation*} 
En tant que $\HP$-module, c'est la somme directe de 
$\Xb\CP$, $(\Yb +\Zb ) \CP$ et $\Tb \CP$. 
D'après le paragraphe \ref{Rentiers}, $\Xb\CP$ est la représentation de 
Steinberg.  
Elle est irréductible et de dimension $q^3$. 
On déduit du paragraphe précédent que le $\HP$-module $\Xb\CP$ est 
isomorphe à une somme directe de $q^3$ copies de $\HP\Xb$.  
C'est un $\HP$-module projectif. 

L'élément $\Tb$ s'identifie dans $\CP$ à la fonction caractéristique de $\G$.  
On en déduit que $\Tb\CP$ est un espace vectoriel de dimension $1$ sur 
lequel $\G$ agit trivialement (voir le paragraphe \ref{Rentiers}).  
En tant que $\HP$-module, $\Tb\CP$ est isomorphe à $\HP\Tb$ et c'est un 
$\HP$-module projectif.

La sous-représentation $\Sb_{1}^*\CP^{}\subseteq\CP$ est engendrée par 
la fonction caractéristique du
sous-groupe parabolique $\P_{1}=\B\cup \B s_{1}\B$ de $\G$.
Elle est donc isomorphe à l'induite 
$\ii_{\P_1}(1)$ qui est de dimension $(\G:\P_{1})=1+q+q^2$.  
Or $\Sb_{1}^*\CP^{}$ est la somme directe de $\Zb \CP$ et $\Tb\CP$, 
donc $\Zb\CP$ est de dimension $q+q^2$.  
Ainsi, la dimension de $\Yb \CP$ est aussi $q+q^2$.

Le noyau de la restriction de $\Sb_{2}$ à $\Yb \CP$ contient 
$\Yb \Sb_2\CP= \Sb_{2}^*\Zb \CP^{}$, 
et le noyau de la restriction de $\Sb_{2}^*$ à $\Zb \CP$ contient 
$\Sb_{2}\Yb \CP=\Zb  \Sb_{2}^*\CP^{}$.
On a les complexes~: 
\begin{equation}
\label{c1}
0\to\Yb\Sb_2\CP\to\Yb\CP\overset{\Sb_{2}}
\longrightarrow\Sb_{2}\Yb\CP\to0
\end{equation}
et~:
\begin{equation}
\label{c2}
0\to\Zb\Sb_{2}^*\CP\to\Zb\CP\overset{\Sb_{2}^*}
\longrightarrow\Sb_{2}^*\Zb\CP\to0
\end{equation}
de représentations de $\G$, dont nous discutons l'exactitude. 
D'après la remarque \ref{Gizzi}, les repré\-sentations $\Yb\Sb_2\CP$ et 
$\Sb_{2}\Yb \CP$ sont ir\-ré\-duc\-ti\-bles.  
Comme elles ne sont isomorphes ni au caractère trivial, ni à la représentation 
de Steinberg, elles sont donc (d'après la proposition \ref{Swift}) de dimension~: 
\begin{equation*}
(p(p+1)/2)^r, 
\end{equation*} 
où l'on a posé $q=p^r$.
Par conséquent, chacun des complexes est exact 
si et seulement si $q=p$.  

Le $\HP$-module à gauche $(\Yb +\Zb )\CP$ est projectif si et seulement si, 
pour tout idéal à droite $\A\subseteq\HP$, l'application~:
\begin{equation*}
\A\otimes_{\HP} (\Yb +\Zb )\CP\rightarrow \CP
\end{equation*}
est injective.  
Il suffit de tester cette propriété sur les idéaux indécomposables de $\HP$, 
et, puisque $\Xb $, $\Yb $, $\Zb $ et $\Tb$ sont des idempotents orthogonaux, 
seuls les cas des idéaux $\Sb_{ 2} \Yb \HP$ et $\Sb_{ 2}^* \Zb \HP^{}$ 
néces\-si\-tent une vérification. 
Traitons le cas de $\Sb_{ 2} \Yb \HP$ en utilisant les complexes (\ref{h1}) 
et (\ref{c1}). 
Le cas de $\Sb_{ 2}^* \Zb \HP^{}$ s'obtient de façon analogue en utilisant 
les complexes (\ref{h2}) et (\ref{c2}). 

D'après \cite[Chapitre 1, \S2, n°11]{Bki}, un élément 
$\Sb_{2}\Yb \otimes c\in \Sb_{2}\Yb \otimes_{\HP}( \Yb +\Zb )\CP$ est 
nul si et seulement s'il existe une famille finie $(h_{i})_{i}$ de $\HP$ et
une famille finie $(c_{i})_{i}$ de $(\Yb +\Zb )\CP$ telles que 
$c=\sum_i{h_{i}} c_{i}$ et $\Sb_{2}\Yb  h_{i}=0$ pour tout $i$, c'est-à-dire, d'après
l'exactitude de (\ref{h1}), si et seulement si 
$c\in \Yb \Sb_{2}\CP+\Zb \CP$.  
Si $q=p$, le complexe (\ref{c1}) est exact, donc cette condition est 
équivalente à 
$\Sb_{2}\Yb c=0$ et l'homomorphisme $\Sb_{2}\Yb \otimes_{\HP} (\Yb +\Zb 
)\CP\rightarrow \CP$ est injectif.  
Si $q\neq p$, tensoriser (\ref{h1}) par
le $\HP$-module $(\Yb +\Zb )\CP$ donne le complexe $(\ref{c1})$ qui n'est 
pas exact.  
Donc $(\Yb +\Zb )\CP$ n'est pas plat. 
La proposition \ref{OldStuff1} est démontrée. 
\end{proof}

\subsection{Platitude de $\C$}
\label{Prem2}

On a le résultat suivant.

\begin{prop}
\label{OldStuff2}
Le $\H$-module $\C$ est plat si et seulement si $q=p=2$.
\end{prop}

\begin{proof}
Le fait que $\C$ n'est pas plat sur $\H$ lorsque $q$ est différent de $p$ est 
donné par le corollaire \ref{Athos}. 
On suppose maintenant que $q=p$.  

Supposons que $p=2$.  
Dans ce cas, le groupe $\hat\T$ est réduit au caractère trivial, et $\C$ est 
égal à $\CP$. 
D'après la proposition \ref{CetSesCopains} et la 
proposition \ref{OldStuff1}, le $\H$-module $\C$ est plat.  

Supposons que $p>2$, et soit $\chi\in\hat\T$ un caractère de $\T$.  
Le $\H$-module correspondant à $\C_{\chi}$ est $\e_{\chi}\H$, 
qui est de dimension $6$ et de base $\{\e_{\chi}\t_{w}\ \vert\ w\in\W_{0}\}$ 
en tant que $\RR$-espace vectoriel. 
Puisque tous les $\H$-modules simples sont de dimension $1$ 
(voir par exemple \cite[Theorem 6.10 (iii)]{CE}), ce mo\-du\-le 
est de longueur $6$ dans $\Mod$.  
Si $\C$ était un $\H$-module plat, le foncteur $\FF$ des 
$\U$-invariants fournirait, d'après la proposition \ref{DS2019}, 
une équivalence entre $\Ee$ et $\Mod$. 
Pour montrer que $\C$ n'est pas un module plat, il suffit donc 
de trouver un caractère $\chi$ tel que la représentation $\C_{\chi}$ de $\G$
soit de longueur strictement supérieure à $6$ dans $\Ee$. 
Remarquons que, 
puisque toute représentation non nulle de $\G$ admet un vecteur $\U$-invariant
non  trivial,  les  éléments  irréductibles  des  catégories  $\Ee$  et  $\Rr$
coïncident.  
Il suffit donc de trouver $\chi$ tel que $\C_{\chi}$ est de
longueur strictement supérieure à $6$ dans $\Rr$, ce qui 
a été fait à la proposition \ref{BartonFink}. 
\end{proof}

Associé au corollaire \ref{ThomasMarvel}, ce résultat fournit le corollaire 
suivant. 

\begin{coro}
\label{Negatif}
On suppose que $q\neq2$.
Alors, pour tout $n\>3$, le $\H^{(n)}$-module $\C^{(n)}$ n'est pas plat. 
\end{coro}

On peut traiter le cas de $\CP$ de façon analogue.
On note $\CP^{(n)}$ et $\HP^{(n)}$ les quantités $\CP$ et $\HP$ associées à 
$\G=\GL_n(\k)$ pour $n\>1$.  
On déduit de la proposition \ref{OldStuff1} le résultat 
suivant. 

\begin{coro}
\label{NegatifP}
On suppose que $q\neq p$.
Alors, pour tout $n\>3$, le $\HP^{(n)}$-module $\CP^{(n)}$ n'est pas 
plat. 
\end{coro}

\section{Représentations de $\GL_{n}$ sur un corps $p$-adique}
\label{RepPadic}

\subsection{Préliminaires\label{Preli}}

\subsubsection{\label{Dorian}}

Soit $\F$ un corps localement compact non archimédien de corps résiduel $\k$.  
On désigne par $\Oo$ l'anneau des entiers de $\F$ et par $\p$ l'idéal 
maximal de $\Oo$.
On note $q$ le cardinal de $k$. On fixe une uniformisante $\varpi$ de $\F$. 

Étant donné un entier $n\>1$, on pose $\Gp=\GL_{n}(\F)$ et 
$\K=\GL_{n}(\Oo)$.  L'image de $\K$ par réduction modulo $\p$ est le groupe 
$\Gf=\GL_{n}(\k)$ des sections précédentes.  On note $\I$ le sous-groupe de 
$\K$ constitué 
des matrices dont la réduction modulo $\p$ est triangulaire su\-pé\-rieu\-re, 
\ie égale à $\B$.  C'est le sous-groupe d'Iwahori (supérieur) standard. Son 
unique pro-$p$ sous groupe de Sylow est noté $\I_1$, c'est l'ensemble des 
matrices dont la réduction modulo $\p$ est unipotente supérieure, \ie égale à 
$\U$. On l'appelle le pro-$p$-sous-groupe d'Iwahori de $\Gp$. 
Pour $m\>1$, on pose $\K_m=1+\p^{m}\M_n(\Oo)$. 

Par \textit{représentation} de $\Gp$ ou de l'un de ses sous-groupes fermés, 
on entendra
re\-pré\-sen\-ta\-tion lisse à coef\-fi\-cients dans $\RR$.

On désigne par $\Bp$ le sous-groupe de Borel de $\Gp$ des matrices
triangulaires  supérieures, de décomposition de Levi
$\Bp=\Tp\Up$ où $\Tp$ désigne le tore diagonal et $\Up$ le sous-groupe unipotent supérieur de $\Gp$.
On identifie le tore fini $\Tf$ de $\Gf$ avec un sous-groupe de $\Tp$ grâce au relèvement de Teichmüller
$\k^\times\rightarrow \Oo^\times$.

On désigne par
$\L$ le sous-groupe de $\Gp$ des matrices diagonales dont les c\oe fficients
non nuls sont des puissances de $\varpi$, et par
$\L^{(1)}$ le sous-groupe de $\Gp$ engendré par $\L$ et $\T$. 
On note $\Wp$ le produit  semi-direct de $\W_0$ par $\L^{(1)}$. Il s'identifie
à un sous-groupe de $\Gp$ qui constitue un système de représentants des 
doubles classes de $\Gp$ modulo $\I_1$.  

Notons que le groupe de Weyl affine étendu, 
système de représentants des doubles classes de $\Gp$ modulo le sous-groupe 
d'Iwahori $\I$, est le sous-groupe de $\Wp$ engendré par $\W_0$ et $\L$. On le 
notera $\Wp'$.

\medskip

On considère la donnée radicielle affine associée à $(\Gp, \Bp, \Tp)$. On se 
réfère à \cite[1]{Lu} par exemple. 
Les racines s'identifient avec celles de la donnée radicielle finie décrite au
paragraphe \ref{HeckeFinie} et l'on note à nouveau $\Root=\Root^+\cup \Root^-$ 
leur ensemble.  
Le sous-ensemble des racines positives $\Root^+$ est le semi-groupe engendré 
par l'ensemble des racines simples ${\rm \Pi}=\{\root_1, \dots,\root_{n-1}\}$. 
On considère la racine $\root_i$ comme le morphisme $\L\rightarrow \ZZ $ 
suivant~: 
$$\root_i: {\rm diag}(\varpi^{x_1},\varpi^{x_2}, \dots, \varpi^{x_n})\longmapsto
{x_{i+1}-x_{i}}.$$ 
Il s'étend  par inflation  en un morphisme  $\root_i:\L^{(1)}\rightarrow \ZZ$.
L'action naturelle de $\W_0$ sur $\L^{(1)}$ induit une action de $\W_0$ sur 
l'ensemble des racines que l'on note 
$(w_0, \root)\mapsto w_0\root$. 

\medskip

On définit comme dans \cite{Lu} l'ensemble des racines affines par 
$\tilde\Root:=\Root\times \ZZ$. On considère $\Root$ comme un sous-ensemble de 
$\tilde \Root$ en identifiant $\root\in\Root$ avec $(\root,0)\in \tilde 
\Root$. L'action de $\Wp$ sur l'ensemble des racines affines est définie comme 
suit~: $$w_0\lambda :\: (\root, k)\mapsto (w_0\root, k-\root(\lambda)).$$ Les 
ensembles des racines affines positives et négatives sont respectivement~: 
$$\tilde \Root^+:=\{(\root, k), \, \root\in \Root, \, k>0\}\cup\{(\root, 0),\, 
\root\in \Root^+\},$$ $$\tilde \Root^-:=\{(\root, k), \, \root\in \Root, \, 
k<0\}\cup\{(\root, 0),\, \root\in \Root^-\}.$$ 
Le groupe $\Wp$ est muni d'une longueur $\ell$ 
qui prolonge la longueur sur $\W$ décrite au paragraphe \ref{HeckeFinie} (on 
se réfère à \cite[\S1.4]{Lu} qui définit la longueur sur $\Wp'$.  
Pour 
l'extension au cas de $\Wp$, voir \cite[\S1.2]{VigProp}). La longueur d'un 
élément $w\in \Wp$ est le nombre de racines affines 
positives rendues négatives sous l'action de $w$. 

\medskip

Le résultat suivant provient de \cite[Proposition 2.5]{OParab}.

\begin{prop}\label{AdagioFaure}
Il existe un système de représentants $\Dd$ des classes à gauche 
$\Wp'/\W_0$ tel que pour tout $d\in\Dd$ et $w\in\Wp'$, on a~:
\begin{equation*}
\ell(dw)=\ell(d)+\ell(w). \end{equation*} 
Chaque $d\in\Dd$ est l'unique élément de longueur minimale dans $d\W_0$.
\end{prop}

\begin{rema} 
\begin{enumerate}
\item 
Puisque les éléments de $\T$ sont de longueur nulle dans $\Wp$,
l'ensemble $\Dd$ est aussi un système de représentants des classes à gauche 
$\Wp/\W$. Chaque $d\in \Dd$ est de longueur minimale dans $d\W$ mais il y a 
d'autres éléments de même longueur dans cet ensemble. 
\item 
L'ensemble $\Dd$ est un système de représentants des doubles classes 
$\I_{1}\backslash \Gp/\K$. 
\item 
D'après la preuve de la proposition \emph{loc.cit.}, 
l'ensemble $\Dd\subseteq \Wp'$ est exactement l'ensemble des $d\in \Wp'$ 
vérifiant $d\Root^+\subseteq{\tilde \Root}^+$. 
Par définition de l'action de $\Wp$ sur les racines affines,  cela signifie que $d\in\Dd$ si et
seulement si $d\in \Wp'$ et $d(\I_1\cap \Up)d^{-1}\subseteq \I_1$.
\end{enumerate}
\label{really?}
\end{rema}

La propriété suivante est prouvée par \cite[Proposition 2.7]{OParab}. 

\begin{lemm} \label{cocoon}
Soient $d\in \Dd$ et $s\in \SS_0$. Si $\ell(sd)=\ell(d)-1$ alors $sd\in\Dd$. Si 
$\ell(sd)=\ell(d)+1$ alors $sd\in \Dd$ ou bien $sd\W_0=d\W_0$.
\end{lemm}

\subsubsection{\label{Shou}}

On note $\Cc=\ind_{\I_1}^{\Gp}(1)$ la représentation de $\Gp$ obtenue 
par induction compacte à partir du ca\-rac\-tè\-re trivial de $\I_1$.
Elle s'identifie à l'espace des fonctions à support fini sur les classes à droite $\I_1\backslash\Gp$ muni de 
l'action de $\Gp$ par translation à droite.
On note $\Hh=\End_\Gp(\Cc)$ la $\RR$-algèbre de ses $\Gp$-endomorphismes. 
Par réciprocité de Frobenius, $\Hh$ s'identifie canoniquement à l'espace 
$\Cc^{\I_1}$ des fonctions de $\Gp$ dans $\RR$ qui sont invariantes par 
$\I_1$ par translations à droite et à gauche, muni du produit de 
convolution d'unité $\1_{\I_1}$.

\medskip

Le module universel fini $\C$ attaché au couple $(\Gf,\Uf)$, 
dont on a étudié les propriétés comme représentation de $\Gf$ 
et comme $\H$-module dans les sections 
\ref{Sec1} et \ref{FoncParaIRJ}, est maintenant vu comme une 
représentation de $\K$ qui se factorise par l'action 
triviale de $\K_{1}$. 
Comme telle, on l'identifie au sous-espace $\ind_{\I_{1}}^\K(1)$ 
des fonctions de $\Cc$ à support dans $\K$.
En considérant les espaces des vecteurs 
$\I_{1}$-invariants de ces représentations, cette 
identification fournit une injection naturelle de l'algèbre de Hecke finie 
$\H$ dans $\Hh$. 
On assimilera désormais $\H$ à son image dans $\Hh$~: 
c'est la sous-algèbre de $\Hh$ engendrée par les fonctions caractéristiques 
des doubles classes de la forme $\I_{1}w\I_{1}$ pour $w\in\W\subseteq\Wp$. 
Ainsi, on étend en toute légitimité les notations du paragraphe 
\ref{HeckeFinie}~: les fonctions caractéristiques~:
\begin{equation*}
\t_{w}=\1_{\I_1 w\I_1},
\quad 
w\in\Wp,
\end{equation*}
forment une base de $\Hh=\Cc^{\I_1}$ comme $\RR$-espace vectoriel.

\medskip

On a \begin{equation}\label{paw}\t_{w}\t_{w'}=\t_{ww'}\textrm{ pour tous les éléments }w,w'\in\Wp
\textrm{ vé\-ri\-fiant }\ell(w)+\ell(w')=\ell(ww').\end{equation} 
Avec celles de la remarque \ref{Fourier}, ce sont les seules relations dans 
$\Hh$ dont nous ferons usage. 
On trouve dans \cite{VigProp} une présentation de $\Hh$ 
par générateurs et relations.

\medskip

On note $\X$ l'immeuble de Bruhat-Tits affine réduit de $\Gp$ sur $\F$. Pour
$k\in\{0,\dots,n-1\}$  on  note  $\X_k$  l'ensemble  des  simplexes,  ou  encore
\emph{facettes},  de  $\X$ de  dimension  $k$. Il  est  muni  d'une action  de
$\Gp$. Les chambres de $\X$ sont 
les facettes de dimension maximale $n-1$, les sommets de $\X$ sont les 
facettes de dimension $0$.  
Ces derniers sont en bijection avec les classes d'homothétie des $\Oo$-réseaux de $\F^n$.
On note $\sigma_0$ le sommet  correspondant à la classe d'homothétie du réseau
$\Oo e_1\oplus\dots\oplus \Oo e_n$ engendré  par la base canonique de $\F^n$ et
l'on   choisit   $\sigma_0$   pour   origine  de   l'immeuble.    On   appelle
\emph{appartement standard}  de $\X$ le  complexe simplicial dont  les sommets
correspondent  aux réseaux  de la  forme $\varpi^{k_1}\Oo  e_1\oplus\dots\oplus
\varpi^{k_n}\Oo  e_n$,  avec  $k_1,  \dots,  k_n\in \ZZ$,  formés  sur  la  base
canonique. Un  appartement de $\X$  est un conjugué de  l'appartement standard
sous l'action d'un élément de $\Gp$.  

Rappelons que deux sommets $\sigma$ et $\sigma'$ 
sont  voisins  s'il existe  des  réseaux  $\Lambda$  et $\Lambda'$  de  $\F^n$
correspondant    respectivement   à   $\sigma$    et   $\sigma'$    tels   que
$\varpi\Lambda'\subseteq \Lambda\subseteq \Lambda'$. 
On dira que deux sommets voisins sont à distance $1$ l'un de l'autre et 
la distance combinatoire entre deux sommets quelconques est alors définie par
récurrence comme dans \cite[2.1.3]{BO}.
 Lorsque l'on dira qu'un sommet  $\sigma\in \X_0$ est à distance $m\in\NN$, il
 sera sous-entendu qu'il est à distance  $m$ de l'origine. De même, on parlera
 la boule de  rayon $m$ sans préciser qu'il s'agit de  la boule fermée centrée
 en $\sigma_0$, qui contient exactement tous les sommets à distance $\< m$. 

\medskip

L'immeuble   $\X$   est   étiquetable,   comme  rappelé   par   exemple   dans
\cite[2]{Broussous}~: il existe une application simpliciale 
$\lambda: \X\rightarrow \Delta_{n-1}$ qui  respecte la dimension des simplexes,
où $\Delta_{n-1}$ est le simplexe standard construit sur 
$\{0, 1,\dots, n-1\}$.  
Un étiquetage de $\X$ est unique à un automorphisme de $\Delta_{n-1}$ près et 
l'on en fixe un pour la suite.  
Pour  $i\in\{1,  \dots,  n-1\}$,  on  considère un  simplexe  $\tau=\{t_0,  \dots,
t_{i-1}\}$  de dimension $i-1$  contenu dans  un simplexe  $\sigma=\{s_0, \dots,
s_i\}$ de dimension $i$.  
Le nombre d'incidence $[\sigma, \tau]$ de $\tau$ dans $\sigma$ est défini 
comme suit~: 
$$[\sigma, \tau]=(-1)^j\textrm{ si }\{\lambda(s_0), \dots., 
\lambda(s_i)\}-\{\lambda(t_0), \dots, \lambda(t_i)\}=\{j\}.$$ 

À toute facette $\s$ de $\X$ correspond un pro-$p$-sous-groupe 
$\U_{\s}$ de $\Gp$ (\cite[\S1]{SSReso}), qui est le pro-$p$-radical du sous-groupe 
parahorique compact de $\Gp$ fixant $\s$ point par point.
Pour tout élément $g\in\Gp$, on a $\U_{g\cdot\s}=g\U_{\s}g^{-1}$. 
Suivant \cite[\S 1]{SSReso},
pour chaque entier $i\in\{0,\dots,n-1\}$ et toute représentation 
$\V$ de $\Gp$, on pose~:
\begin{equation}
\label{GulliversTravels}
\Ff_{i}(\X,\V)=\bigoplus_{\s\in X_i}\V^{\U_{\s}},
\end{equation}
espace naturellement muni d'une structure de représentation de $\Gp$. 
Plus précisément, si $g\in\Gp$ et si $f\in\Ff_{i}(\X,\V)$, 
alors $gf$ est la fonction définie par 
$\s\mapsto g\cdot f(g^{-1}\cdot\s)$.
Pour $i\in\{1,\dots,n-1\}$, on dé\-fi\-nit une application de transition~: 
\begin{equation}
\begin{array}{ccc}
\partial_{i}:\Ff_{i}(\X,\V)&\to&\Ff_{i-1}(\X,\V),\\
f&\mapsto&\Big(\tau\mapsto\sum\limits_{\tau\subseteq\s
\atop{\dim(\s)=i}}[\s:\tau]\, f(\s)\Big).
\label{transition}
\end{array}
\end{equation}
Cette application est équivariante sous
l'action  du sous-groupe de $\Gp$ des éléments à déterminant inversible dans $\Oo$.
On a ainsi un complexe augmenté~: 
\begin{equation}
\label{CoefficientSystems}
0\to\Ff_{n-1}(\X,\V)\to\dots\to\Ff_{0}(\X,\V)\xrightarrow{\alpha}\V\to0, 
\end{equation}
où $\alpha$ désigne l'homomorphisme d'augmentation défini par~:
\begin{equation*}
\alpha(f)=\sum\limits_{\dim(\s)=0}f(\s).
\end{equation*}
Ce complexe est appelé le \emph{système de c\oe fficients} associé à la représentation $\V$ de $\Gp$.

\medskip

Si l'on choisit $\V=\Cc$, l'espace (\ref{GulliversTravels}) est 
un $\Hh$-module à gauche ($\Hh$ opère sur 
(\ref{GulliversTravels}) composante par composante), les 
applications de transition (\ref{transition}) sont $\Hh$-linéaires
et 
\begin{equation}
\label{Fizzy}
0\to\Ff_{n-1}(\X,\Cc)\to\dots\to\Ff_{0}(\X,\Cc)\xrightarrow{\alpha}\Cc\to0, 
\end{equation}
est un complexe de $\Hh$-modules qui
 est exact d'après la preuve du théorème de 
\cite[\S3]{SSReso}.

\medskip

On note $\CcI=\ind_{\I}^{\Gp}(1)$ l'induite compacte du caractère trivial de 
$\I$. 
Elle s'identifie à l'espace des fonctions à support fini dans
$\I\backslash\Gp$ muni de l'action de $\Gp$ par translation à droite.  Soit
 $\HhI=\End_\Gp(\CcI)$ la $\RR$-algèbre de ses 
$\Gp$-endo\-mor\-phis\-mes. 
Par réciprocité de Frobenius, $\HhI$ s'identifie canoniquement à l'espace 
$\CcI^{\I}$ des fonctions de $\Gp$ dans $\RR$ qui sont invariantes par 
$\I$ par translations à droite et à gauche, muni du produit de 
convolution d'unité $\1_{\I}$.

\begin{rema} Au paragraphe \S \ref{Borel}, nous avons introduit l'idempotent central $\e_1$ de $\H$. Vu comme élément de
$\Hh$ c'est encore un idempotent central et il correspond  à
la projection naturelle $\Cc\rightarrow \CcI$. Ainsi, $\e_1\Hh$ est un facteur direct de $\Hh$ comme $\Hh$-module à droite et   $\CcI$ est un facteur
direct  de $\Cc$ comme $\Hh$-module à gauche. 
\label{Wsquare}
\end{rema}

Le lemme suivant  est immédiat.

\begin{lemm}
\label{Fox} Soit $\Aa$ un idéal à droite  de  $\Hh$.
Pour tout $i\in\{0, \dots, n-1\}$, on a \\
$\Aa\Im(\partial_{i})=\partial_{i}(\Aa\Ff_{i}(\X,\Cc))$,  $\Aa\Ff_i(\X,\Cc)=\Ff_i(\X,\Aa\Cc)$ et 
$\Ker(\partial_i)\cap\e_1 (\Ff_i(\X,\Cc))=\e_1(\Ker(\partial_i))$.
\end{lemm}

 En 
appliquant ce lemme
 avec l'idéal à droite de $\Hh$ engendré par $\e_\1$,  le complexe (\ref{Fizzy}), donne 
un complexe augmenté de $\HhI$-modules~:
\begin{equation}
\label{CoefficientSystemsII}
0\to\Ff_{n-1}(\X,\CcI)\to\dots\to\Ff_{0}(\X,\CcI)\xrightarrow{\alpha}\CcI\to0 
\end{equation}
qui est également  exact. 

\subsection{L'espace des vecteurs $\K_{1}$-invariants de $\Cc$}

Dans ce paragraphe, on donne une description de 
l'espace des vecteurs $\K_{1}$-invariants de $\Cc$. 

\begin{prop}
\label{DistanceUn}
L'espace des vecteurs $\K_{1}$-invariants de $\Cc$ est engendré 
en tant que $\Hh$-mo\-du\-le par les $x\cdot\1_{\I_1}$, avec $x\in\Gp$
vérifiant $x^{-1}\K_1x\subseteq\I_1$, \ie par les $\1_{\I_1 x}$, 
avec $x\in\Gp$ vérifiant $\I_{1}x\K_1=\I_1x$.
\end{prop}

\begin{rema}
Ce résultat, ainsi que la preuve que nous en donnons, est valable 
 sur un corps  algébriquement clos de caractéristique 
quelconque, non nécessairement $p$.
\end{rema}

\begin{proof}
Comme espace vectoriel, l'espace des vecteurs $\K_{1}$-invariants 
de $\Cc$ est engen\-dré par les fonctions caractéristiques de 
la forme $\1_{\I_1 g\K_{1}}$, avec $g\in\Gp$. 
Si $f=\1_{\I_1 g\K_{1}}$ 
est une telle fonc\-tion, il suffit de montrer qu'il existe 
une fonction $\K_{1}$-invariante de la forme 
$e=\1_{\I_1 x}$, avec $x\in\Gp$ vé\-ri\-fiant $\I_{1}x\K_1=\I_1x$, 
telle qu'on ait $f\in\Hh\cdot e$.
Puisque $\Hh\cdot e$ est l'espace des vecteurs de $\Cc$
invariants par $x^{-1}\I_1x$, cela revient à montrer qu'il y a 
un élément $x\in\Gp$ tel que $x^{-1}\I_1x$ fixe $f$ et contienne 
$\K_{1}$ \ie tel que
$$\K_1\subseteq x^{-1}\I_1 x\subseteq g^{-1}\I_1 g\K_1.$$
Par  la décomposition de Bruhat, 
l'élément $g$ se décompose dans $\Gp=\I_1\Wp\I_1$. De plus, tout élément
de $\Wp$ s'écrit comme un produit d'un élément
de $\Dd$ défini par la proposition \ref{AdagioFaure} et d'un élément de $\W$ 
qui est inclus dans $\K$. Puisque $\K$ normalise $\K_1$, 
 on peut se ramener au cas où $g$ 
est un élément de $\Dd$. D'après la remarque \ref{really?},  on a  alors \begin{equation} g (\I_1\cap \Up) g^{-1}\subset \I_1\label{Oren}.\end{equation}
On note $\Up^{-}$ le sous-groupe 
unipotent opposé à $\Up$ relativement à $\Tp$ et l'on écrit la dé\-com\-po\-si\-tion d'Iwahori~:
\begin{displaymath}
\I_1=(\I_1\cap\Up)\cdot(\I_1\cap\Tp)\cdot(\I_1\cap\Up^-).
\end{displaymath}
Les deux derniers facteurs sont contenus dans $\K_1$
 de sorte que (\ref{Oren}) implique
 $$\I_1\subseteq g^{-1}\I_1 g\K_{1}.$$ 
Autrement dit, $x=1$ convient.
\end{proof}

\subsection{Du cas fini au cas $p$-adique}
\label{Fairway}

Le résultat suivant est une conséquence immédiate de la proposition \ref{AdagioFaure}. 

\begin{prop}
\label{BrocoliRabe}
L'algèbre de Hecke $\Hh$ est un module à droite (respectivement  à gauche) libre de base 
$\{\tau_{d}\}_{d\in\Dd}$  (respectivement  $\{\tau_{d^{-1}}\}_{d\in\Dd}$)  sur
l'algèbre  de  Hecke  finie  $\H$,  où  l'ensemble $\Dd$  est  défini  par  la
proposition \ref{AdagioFaure}. 
\end{prop}

\begin{rema}
\label{Brooklyn}
Le sous-espace $\t_{d}\H$ de $\Hh$ est exactement celui des fonctions 
$\I_{1}$-invariantes à support dans $\I_{1}d\K$. 
\end{rema}

\subsubsection{} On déduit de la proposition \ref{BrocoliRabe} les  deux lemmes suivants.

\begin{lemm}
Tout $\H$-module à gauche $\mm$ est un facteur direct de la restriction de 
$\Hh\otimes_\H \mm$ à $\H$.
\label{Honey}
\end{lemm}

\begin{proof} Comme le lemme \ref{intents} a été prouvé grâce au lemme
\ref{kitten}, on prouve le présent résultat en utilisant le lemme \ref{cocoon}.
Par la proposition  \ref{BrocoliRabe}, l'espace vectoriel $\Hh\otimes_{\H} \mm$ s'identifie à la somme directe des
  $ \t_{d}\mm$ pour $d\in \Dd$.
On définit le sous-espace vectoriel $\nn$ de  $\Hh\otimes_{\H} \mm$
comme la somme
directe des  $ \t_{d}\mm$ pour $d\neq 1$ et l'on montre qu'il est stable sous 
l'action de l'algèbre $\H$,  engendrée par 
$\{\t_{s},\,\t_{t},\:s\in \SS_{0},\,t\in \T\}$.

\noindent\textbf{a/}  Comme dans la preuve du lemme \ref{intents}, un sous-espace de la forme 
$\t_d \mm$  avec $d\in\Dd$ est stabilisé par l'action de $\t_t$, $t\in \T$.

\noindent\textbf{b/} Le point a/ assure que $\nn$ est stable sous l'action de $\t_t$ pour tout $t\in \T$.

\noindent\textbf{c/}  Soit $s\in  \SS_{0}$. Vérifions que
$\t_{s}$ stabilise l'espace $\nn$. Soit   $d\in\Dd$, $d\neq 1$.\\
Si $\ell(sd)= \ell(d)-1$ alors 
$\t_{d}=\t_{s}\t_{sd}$ et $sd$ appartient à $\Dd$ d'après le lemme \ref{cocoon}.  D'après la remarque
\ref{Fourier}, la combinaison linéaire $a=-\sum_{\chi\in\hat \T,\, \cs=\chi}\e_{\chi}$ d'éléments de $\{\t_{t}, t\in \T\}$ vérifie  $\t_{s}\t_{d}=\t_{s} a\t_{sd}$ de sorte que 
$\t_{s}\t_{d}\mm=\t_{s} a\t_{sd}\mm\subset \t_s \t_{sd}\mm=\t_{d}\mm$ d'après a/.\\
Si  $\ell(sd)= \ell(d)+1$ alors  $\t _{s}\t _{d}=\t _{sd}$. Si 
l'élément $sd$ appartient à $\Dd$, alors 
$\t _{s}\t _{d}\mm=\t _{sd}\mm\subset \nn$. 
Sinon, $sd\in d\W$ d'après le lemme \ref{cocoon} 
donc il existe $w\in \W$ 
tel que $\t _{sd}=\t _{d}\t _{w}$ de sorte que
$\t _{s}\t _{d}\mm=\t_{d}\t_{w}\mm\subset \t _{d}\mm\subset \nn$. 

Ainsi, $\nn$ est bien stable sous l'action de $\H$.
\end{proof}

\begin{lemm}
\label{Fur}
Un $\H$-module à gauche de type fini 
$\mm$ est plat si et seulement si le $\Hh$-module $\Hh\otimes_\H\mm $ est 
plat, et dans ce cas, ils sont même tous deux projectifs. 
\end{lemm}

\begin{proof} 
 Soient $\A$ un idéal à droite de $\H$ et $\EuScript{A}$ l'idéal 
à droite de $\Hh$ qu'il engendre, \ie que 
$\EuScript{A}=\A\Hh$, qui est isomorphe à $\A\otimes_{\H}\Hh$ par la proposition \ref{BrocoliRabe}. 
Sous l'hypothèse que $\Hh\otimes_\H\mm $ est un $\Hh$-module à gauche plat,
 l'application linéaire naturelle
\begin{equation*}
\Aa\otimes _\H\mm\simeq \A\otimes _\H\Hh\otimes_\H\mm
\longrightarrow \Hh\otimes _\H\mm
\end{equation*}  est injective. 
D'après le lemme \ref{Honey},  l'espace
 $\A\otimes_\H\mm$ s'injecte dans $\A\otimes _\H\Hh\otimes_\H\mm$ et la 
 restriction de l'application précédente à cet espace n'est autre que 
 l'application linéaire naturelle
$\A\otimes _\H\mm
\longrightarrow \mm$, qui est donc également injective.  
D'après \cite[chap. I, \S2, n$^\circ$3, Proposition 1 a)]{Bki}, 
cela suffit à assurer la platitude de $\mm$.

Supposons que $\mm$ est un $\H$-module à gauche plat, donc projectif par le lemme \ref{lam}~: c'est  un
facteur direct d'un $\H$-module libre.  Par la proposition
\ref{BrocoliRabe}
tensoriser par $\Hh$
 montre alors que $\Hh\otimes _\H\mm$ est un facteur direct d'un $\Hh$-module libre. 
C'est un $\Hh$-module projectif donc plat.
\end{proof}

\subsubsection{}

\begin{coro}
\label{Alex}
L'application naturelle de $\Hh\otimes_{\H}\C$ dans $\Cc^{\K_1}$ est un 
isomorphisme de $\Hh$-modules à gauche et de représentations de $\K$.  En particulier, la $\K$-représentation 
$\Cc^{\K_1}$ se décompose en la somme directe des $\t_d \C$, $d\in \Dd$.
\end{coro}

\begin{proof}
L'injectivité est une con\-sé\-quen\-ce de la proposition \ref{BrocoliRabe}
par les arguments suivants.
L'espace
$\Hh\otimes_{\H}\C$ se décompose comme représentation de $\K$ en la somme 
directe des $\t_{d}\otimes\C$ pour $d\in \Dd$. 
Les espaces images des $\t_{d}\otimes\C$ dans $\Cc^{\K_1}$ sont en somme directe puisque
$\t_d\C$ est un ensemble de fonctions à support dans $\I_1d \K$. Il suffit donc de vérifier que $\t_d\otimes \C\rightarrow \t_d\C$ est injective.  Pour cela, on remarque que  chaque application $\K$-équivariante $\t_{d}: \C\rightarrow \Cc$ est injective
puisque sa restriction à l'espace $\I_{1}$-invariant $\H$ l'est.

En vertu de la proposition \ref{DistanceUn}, et puisque $\K$ normalise $\K_1$,
la surjectivité sera prouvée lorsque l'on  aura  vérifié que si la fonction $f=\1_{\I_1 d}$ avec $d\in \Dd$ est $\K_1$-invariante, alors elle appartient à l'image de 
$\Hh\otimes _\H\C\rightarrow \Cc^{\K_1}$.

L'hypothèse de $\K_1$-invariance se traduit par $\I_1d \K_1= \I_1d$. Ainsi, par la décomposition d'Iwahori, $\I_1d \I_1= \I_1 d\, (\I_1\cap \Up)$.
D'après la remarque \ref{really?}, on a donc $\I_1 d\I_1= \I_1 d$ de sorte que  $f$ est égale à la fonction caractéristique de  
$\I_1d \I_1$~: on a prouvé que  $f=\t_d \1_{\I_1}\in \t_d \C$. 
\end{proof}

\begin{lemm} 
\label{Wyoming2}
Le $\Hh$-module $\Cc^{\Ip}$ est un facteur direct de $\Cc^{\K_{1}}$.
\end{lemm}

\begin{proof} 
Le corollaire \ref{Wyoming} dit que le $\H$-module $\C^{\U}$ est un
facteur direct de $\C$.  On conclut en notant que dans l'isomophisme 
$\Hh\otimes_{\H}\C\simeq \Cc^{\K_{1}}$ le sous-espace 
$\Hh\otimes_{\H}\C^{\U}$ de $\Hh\otimes_{\H}\C$ s'identifie à $\Cc^{\I_{1}}$. 
\end{proof}

Les lemmes \ref{Honey}  et \ref{Fur} appliqués à $\mm=\C$ donnent respectivement~:
\begin{prop}
\label{facteurdirect} 
Le $\H$-module à gauche $\Cf$ est un facteur direct de $\Cc^{\K_{1}}$. 
\end{prop}

\begin{prop}
\label{LockedOut}
Le $\Hh$-module  à gauche $\Cc^{\K_{1}}$ est plat si et seulement si le $\H$-module 
$\Cf$ est plat. 
\end{prop}
\begin{rema}\label{Summer42}
On a donc prouvé que si $\Cc^{\K_{1}}$ est un $\Hh$-module  plat, il est même projectif.
\end{rema}

\begin{exem} \label{Yellowstone} Par la remarque \ref{Wsquare},
le travail précédent est valable en rempla\c cant
$\Cc$ par $\CcI$, $\Hh$ par $\HhI$,
$\C$ par $\C'$ et $\H$ par $\H'$ avec les notations du paragraphe \ref{Borel}.
De cette remarque, des propositions \ref{PlatGL21}, \ref{PlatGL22} 
et des corollaires \ref{Negatif}, \ref{NegatifP}, 
on déduit les assertions suivantes.
\begin{enumerate}
\item Si $n=1$, le $\Hh$-module $\Cc^{\K_{1}}$ est projectif.
\item Si $n=2$, le $\Hh$-module $\CcI^{\K_{1}}$ est projectif.
\item Si $n=2$, le $\Hh$-module $\Cc^{\K_{1}}$ est projectif si et seulement si $q=p$.
\item 
Si $n=3$, le $\Hh_{1}$-module $\CcI^{\K_{1}}$ est projectif si et seulement si $q=p$.
\item Si $n=3$, le $\Hh$-module $\Cc^{\K_{1}}$ est plat si et seulement si 
$q=p=2$, auquel cas il est même projectif.
\item Si $n\> 4$ et $q\neq2$, le $\Hh$-module $\Cc^{\K_{1}}$ n'est pas 
  plat. 
\item Si $n\> 4$ et $q\neq p$, le $\Hh$-module $\CcI^{\K_{1}}$ n'est pas 
  plat.
 
\end{enumerate}
\end{exem}

\subsubsection{\label{samba}}
\def\Np{{\Nn}}

\def\Mp{\Mm}
On reprend les notations du paragraphe \ref{Diag1}.  En particulier $\P=\M\N$ 
est un sous-groupe parabolique supérieur de décomposition de Levi 
$\P=\M\N$. On désigne par $\Np$ le sous-groupe de $\K$ image réciproque de 
$\N$ par la réduction modulo $\p$, par $\Mp$ celle de $\M$.  Ce sont des 
sous-groupes ouverts et compacts de $\Gp$.  

\begin{prop}\label{DuelAtDawn}
L'isomorphisme $\Hh\otimes_{\H}\C\simeq \Cc^{\K_{1}}$ du corollaire \ref{Alex} induit par restriction un isomorphisme   de $\Hh$-modules à gauche et de représentations de $\Mp$~: $$\Hh\otimes_{\H}\C^{\N}\simeq \Cc^\Np.$$En particulier, le $\H$-module $\C^{\N}$ s'identifie à un facteur direct de $\Cc^\Np$.

\end{prop}
 
\begin{proof}  
D'après le corollaire \ref{Alex}, 
La $\K$-représentation $\Cc^{\K_1}$ s'identifie à la somme directe des 
$\t_d\C$, $d\in \Dd$. 
L'application $\t_d :\C\rightarrow \t_d\C$ est injective comme nous l'avons 
noté dans la preuve du corollaire \ref{Alex}, de sorte que l'espace des 
$\Nn$-invariants de $\t_d\C$ est égal à $\t_d\C^\N$.  Ainsi, on a bien 
l'isomorphisme annoncé. La dernière assertion provient du lemme 
\ref{Honey}.
\end{proof}

\begin{coro} Si $\C_\M$ est un $\H_\M$-module à gauche plat,  alors
 le $\Hh$-module à gauche
$\Cc ^\Np$  est plat et même projectif,  et est un facteur direct de $\Cc^{\K_1}$.
\label{Monsters}
\end{coro}

\begin{proof}
Si $\C_\M$ est un $\H_\M$-module à gauche plat, alors d'après la proposition 
\ref{Monsters0}, le $\H$-module $\C^\N$ est projectif et est un facteur direct 
de $\C$.  On conclut en utilisant l'isomorphisme de $\Hh$-module à gauche 
$\Cc^\Nn\simeq \Hh\otimes_\H\C^\N$.  
\end{proof}

\begin{rema} \label{bloom}Le lemme \ref{Wyoming} assure que $\C^\U$ est un facteur direct de $\C^\N$ comme $\H$-module.
On déduit de la proposition \ref{DuelAtDawn} que $\Cc^{\I_1}$ est un facteur direct de $\Cc^\Np$ comme $\Hh$-module.

\end{rema}

\begin{exem}\label{Sweetgrass}
Nous ferons usage de l'exemple suivant dans la section \ref{IAS}. 
On choisit $n=3$. 
Soit $\M$ le sous-groupe de Lévi standard de ${\rm GL}_{3}(k)$ isomorphe à 
${\rm GL}_{2}(k)\times {\rm GL}_{1}(k)$. On désigne par $\P$ le sous-groupe 
parabolique supérieur de ${\rm GL}_{3}(k)$ associé, de radical unipotent 
$\N=\begin{pmatrix}1&0&k\cr 0&1&k\cr 0&0&1\cr\end{pmatrix}.$ 
Dans ce cas, $\Np=1+
\begin{pmatrix} 
\p &\p&\Oo\\ 
\p&\p&\Oo\\
\p&\p&\p
\end{pmatrix}.
$

Si $q=p$, on déduit de ce qui précède et de la proposition \ref{PlatGL22} 
que $\Cc^\Np$ est un $\Hh$-module projectif facteur direct de 
$\Cc^{\K_{1}}$. Si $q\neq p$, le $\Hh$-module $\Cc^\Np$ n'est pas plat. 
\end{exem}

\begin{rema}\label{parme}
La proposition \ref{DuelAtDawn} et le corollaire \ref{Monsters} sont valables 
en rempla\c cant 
$\Cc$ par $\CcI$, $\Hh$ par  $\HhI$,
$\C$ par $\C'$ et $\H$ par $\H'$ avec les notations du paragraphe \ref{Borel} 
et par la remarque \ref{Wsquare}.  
On déduit alors de la proposition \ref{PlatGL21} que 
$\CcI^\Np$ est un $\Hh_1$-module projectif facteur direct de $\CcI^{\K_{1}}$ 
que $q$ soit égal à $p$ ou non. 
\end{rema}

\section{Le cas de $\GL_2(\F)$}
\label{sec6}

\subsection{Filtration de $\Cc$}

\subsubsection{}

En tant que $\Hh$-module et en tant que représentation de $\K$, 
l'espace $\Cc$ est la limite inductive des $\Cc^{\K_{m}}$ 
pour $m\> 1$.

\begin{prop}
\label{DistanceEm}
L'espace des vecteurs $\K_{m}$-invariants de $\Cc$, pour $m\>1$, 
est engendré en tant que $\H$-mo\-du\-le par les $x\cdot\1_{\I_1}$, 
avec $x\in\G$ vérifiant $x^{-1}\K_mx\subseteq\I_1$, \ie par les 
$\1_{\I_1 x}$, avec $x\in\G$ vérifiant $\I_{1}x\K_m=\I_1x$.
\end{prop}

\begin{proof}

De fa\c con analogue au cas où $m=1$ traité par la proposition 
\ref{DistanceUn}, on va chercher, pour tout $g\in\Gp$, un élément 
$x\in\Gp$ tel que $x^{-1}\I_1x$ fixe $f=1_{\I_1g\K_{m}}$ et contienne 
$\K_{m}$. 
En appliquant la décomposition de Bruhat et puisque $\K$ normalise $\K_m$, on peut se ramener au cas où $g$ est une matrice diagonale de la forme
\begin{equation*}
g={\rm diag}(\varpi^{a_1},\varpi^{a_2})=
\begin{pmatrix}
\varpi^{a_1}&0\\0&\varpi^{a_2}
\end{pmatrix},
\quad
a_1,a_2\in\ZZ,
\end{equation*}
et l'on va chercher $x$ sous la forme diagonale 
${\rm diag}(\varpi^{b_1},\varpi^{b_2})$, avec $b_1,b_2\in\ZZ$.
La condition pour que $x^{-1}\I_1x$ contienne $\K_m$ s'écrit~:
\begin{equation}
\label{eq:lmo}
1-m\<b_2-b_1\<m,
\end{equation}
et la condition pour que $x^{-1}\I_1x$ soit contenu dans 
$g^{-1}\I_1g\K_{m}$ s'écrit~:
\begin{equation}
\label{Condm}
{\rm min}\{m,a_2-a_1\}\<b_2-b_1\<{\rm max}\{1-m,a_2-a_1\}.
\end{equation}
Si $\1_{\I_1g}$ est invariant par $\K_m$, \ie si 
$g^{-1}\I_1g$ contient $\K_m$, on a~:
\begin{equation*}
\label{eq:lmo2}
1-m\<a_2-a_1\<m
\end{equation*}
et il suffit de choisir $x=g$.
Sinon, 
les conditions 
(\ref{eq:lmo}) et (\ref{Condm}) imposent le choix~:
\begin{equation*}
b_2-b_1=
\left\{
\begin{array}{ll}
m & \text{si } a_2-a_1\>m+1,\\
1-m & \text{si } a_2-a_1\<-m.\\
\end{array}
\right.
\end{equation*}
Dans tous les cas, il y a un seul choix possible pour $x$ à un scalaire près. 
\end{proof}

\subsubsection{\label{Merapi}} Dans le cas de ${\rm GL}_2(\F)$, l'immeuble $\X$ a une structure d'arbre. On choisit l'étiquetage des sommets de sorte que l'étiquette de $\s_0$ est $0$.
On appelle  \emph{arête} plutôt que  \emph{chambre} les facettes  de dimension
$1$. Rappelons que deux arêtes quelconques de l'arbre sont 
conjuguées sous l'action d'un élément de $\Gp$.
Pour $m\in \NN$, $m{\> 1}$, on dira  d'une arête qu'elle est à distance $m$ si
l'un des deux sommets qui la constituent est 
à distance $m$ et l'autre à distance $m-1$ (sous-entendu de l'origine $\s_0$).

\medskip

La décomposition de Cartan pour $\Gp$ s'écrit~:
\begin{equation*}
\Gp=\K\begin{pmatrix}\varpi^{\NN}&0\cr 0&1 \cr\end{pmatrix}\K\Zp,
\end{equation*}
où $\Zp$ désigne le centre de $\Gp$, que l'on identifie à $\mult\F$.
Pour $m\in\ZZ$, on note $\l_{m}$ la matrice diagonale 
${\rm diag}(\varpi^m,1)$, et l'on fixe un système de représentants 
$\EuScript{K}_{m}$ de $\K/\K\cap\l_{m}^{}\K\l_{m}^{-1}$.
Alors l'ensemble des $k\l_{m}$, avec $m\>0$ et $k\in\EuScript{K}_{m}$,
est un système de représentants de $\Gp/\K\Zp$.
Pour $m\>0$, les sommets à distance $m$ de l'arbre $\X$ sont les 
$k\l_{m}\ori$ pour $k\in\EuScript{K}_{m}$.

\medskip

Pour tout sommet $\sigma$ à distance $m\> 1$, on désigne par $e(\s)$ 
l'unique arête à distance $m$ contenant le sommet $\s$.  Elle relie $\sigma$ à l'unique
 voisin de $\sigma$ à distance $m-1$. On pose $e(\sigma_0):=\{\sigma_{0},\l_{-1}\sigma_{0}\}$.

\begin{rema}
\label{Incubus}
\begin{enumerate}
\item 
Le fixateur de l'arête $e(\sigma_0)$ est le sous-groupe d'Iwahori $\I$.  
Le pro-$p$-groupe $\U_{e({\s_0})}$ associé à l'arête 
$e({\s_0})$ est le pro-$p$-sous-groupe d'Iwahori $\I_{1}$. 
\item Soit $g\in\Gp$ et $m\>1$. 
Remarquons que le sommet $g\ori$ est à distance $\< m$ si et seulement si 
$\K_m\subset  g \K  g^{-1}$ ou  encore si  et seulement  si $\K_{m+1}\subseteq
g\K_{1}g^{-1}$~; l'arête $ge({\s_0})$ est à distance 
$\< m$ si et seulement si c'est le cas de chacun de ses sommets et l'on 
vérifie que cela signifie que $\K_{m}\subseteq g\I_{1}g^{-1}$, en remarquant 
que $\K\cap \lambda_{-1}\K\lambda_1 =\I$.  
\end{enumerate}
\end{rema}

Un calcul explicite donne le lemme suivant.

\begin{lemm}
\label{Hurt} 
Pour tout $m\> 1$, on a 
$\l_{m}\Cc^{\K_{1}}\subseteq\Cc^{\K_{m+1}}$ et 
$\l_{m}\Cc^{\I_{1}}=(\l_{m}\Cc^{\K_{1}})\cap \Cc^{\K_{m}}$. 
\end{lemm}

\begin{prop}
\label{PSHQ}
Pour $m\> 0$, on a l'égalité suivante~:
\begin{equation*}
\Cc^{\K_{m+1}}=
\sum_{\substack{ {\sigma\in\X_{0}}\\ {\text{\rm à distance } m}} }\Cc^{\U_{\sigma}}. 
\end{equation*}
\end{prop}

\begin{proof} 
L'inclusion indirecte est assurée par le lemme \ref{Hurt}
car $\K_{m+1}$ est distingué dans $\K$.
Montrons l'inclusion directe par récurrence sur $m$. 
Elle est immédiate pour $m=0$ puisque $\U_{\sigma_0}=\K_1$.
Supposons-la vraie à un certain rang $m\> 0$ et montrons-la au rang suivant. 
Puisque les espaces en présence sont stables sous l'action de $\Hh$, il 
suffit, d'après la proposition 
\ref{DistanceEm}, de montrer que toute fonction $\K_{m+1}$-invariante de la 
forme 
$f=\1_{\I_{1}g}$, avec $g\in \Gp$, appartient à l'espace de droite. 
Dire que $f$ est $\K_{m+1}$-invariante signifie que 
$\K_{m+1}$ est inclus dans $g^{-1} \Ip g$. 
Si cette fonction est même $\K_{m}$-invariante, alors on conclut par 
récurrence. 
Sinon, cela  signifie que $\K_{m}$  n'est pas inclus  dans $g^{-1} \Ip  g$ et,
d'après la remarque \ref{Incubus}, que l'arête $g^{-1}e({\s_0})$ est à 
distance exactement $m+1$. L'un des  deux sommets $g^{-1}\ori$ ou 
$g^{-1}\lambda_{-1}\ori$ est donc à distance $m+1$ et l'on note $\sigma$ le sommet en 
question. Dans chacun des deux cas,
on vérifie que $f$ appartient à $\Cc^{\U_{\sigma}}$ 
\end{proof}

\begin{lemm}
Pour tout sommet $\sigma$ à distance $m\> 1$,
on a une injection $\Hh$-équivariante
\begin{equation}
\label{Einziger}
\Cc^{\U_{\sigma}}/\Cc^{\U_{e_{\sigma}}}\hookrightarrow\Cc^{\K_{m+1}}/\Cc^{\K_m}.
\end{equation} 
\end{lemm}

\begin{proof} 
Le sommet $\sigma$ est de la forme $k\lambda_{m}\sigma_{0}$ avec $k\in \K$ et l'on 
vérifie que $e_{\sigma}$ est alors l'image par translation par $k\lambda_{m}$ de 
l'arête $e_{\s_0}$. 
Ainsi, 
$\U_{\sigma}={k\lambda_{m}}\K_{1}{(k\lambda_{m})}^{-1}$ et 
$\U_{e_\sigma}={k\lambda_{m}}\Ip{(k\lambda_{m})}^{-1}$.
Puisque les sous-groupes de congruence $\K_{m}$ et $\K_{m+1}$ sont normalisés 
par $k$, on se ramène au cas de $k=1$ qui est prouvé par le lemme \ref{Hurt}.
\end{proof}

D'après la proposition \ref{PSHQ}, le $\Hh$-module $\Cc^{\K_{m+1}}/\Cc^{\K_m}$ 
est égal à la somme des images des applications (\ref{Einziger}) lorsque 
$\sigma$ parcourt l'ensemble des sommets à distance $m$. On a un morphisme 
$\Hh$-équivariant surjectif 
\begin{equation}\label{Key}
\bigoplus_{\substack{ {\sigma\in\X_{0}}\\ {\textrm{ à distance } m}} 
}\Cc^{\U_{\sigma}}/\Cc^{\U_{e_{\sigma}}}\longrightarrow 
\Cc^{\K_{m+1}}/\Cc^{\K_m}. 
\end{equation}

\subsubsection{} \begin{prop} Le complexe  \begin{equation}
\label{FCS-2}
0\to\Ff_{1}(\X,\Cc)
\xrightarrow{\partial_{1}}
\Ff_{0}(\X,\Cc)\xrightarrow{\alpha}\Cc\to0, 
\end{equation}  défini au paragraphe \ref{Shou} est exact si et seulement si l'application (\ref{Key})
est injective pour tout $m\> 1$.  
\end{prop}

\begin{proof}

Supposons que l'application (\ref{Key}) est injective pour tout $m\> 1$.

Soit $f\in\Ff_{0}(\X,\Cc)$~:  
c'est une fonction de support  ${\rm supp}(f)$ fini sur l'ensemble  $\X_0$ des sommets de l'arbre et   
de valeurs $f({\sigma})\in\Cc^{\U_{\sigma}}$ pour tout $\s\in {\rm supp}(f)$.
Soit $m\> 0$ le rayon minimal tel que le support de $f$
est contenu dans la boule de rayon $ m$.

On suppose que $f$ appartient au noyau de $\alpha$ c'est à dire que 
$\sum_{\sigma\in \X_{0}}f({\sigma})=0$et l'on montre par récurrence sur $m$ que 
$f$ est dans l'image de $\partial_{1}$. 
Remarquons que $m$ ne saurait être égal à $0$ puisque
la restriction de $\alpha$ aux fonctions de support le sommet origine $\ori$
n'est autre que l'inclusion $\Cc^{\K_{1}}\subset \Cc$.

Supposons que $m=1$. 
Alors 
$$\sum_{\substack{{\sigma\in \X_{0}} \\{ \textrm{ à distance } 1}}} 
f({\sigma})=-f({\sigma_{0}})\in \Cc^{\K_{1}}$$ 
donc $f({\sigma})\in \Cc^{\U_{e_{\sigma}}}$  pour tout sommet $\sigma$ à distance $1$
par injectivité de (\ref{Key}). 
On définit  alors  $\varphi\in \Ff_{1}(\X,\Cc)$ de support l'ensemble des arêtes  à distance $1$  de la fa\c con suivante~: pour tout sommet $\sigma$  à distance $1$, la fonction $\varphi$ prend  en $\{\sigma_{0}, \sigma\}$ la valeur
$f({{\sigma}})$.
Suivant la définition de l'application de transition $\partial_{1}$,  on vérifie  que $\partial_{1}(\varphi)\in  \Ff_{0}(\X,\Cc)$ est la fonction de support contenu dans la boule de rayon $1$ donnée par~:
\begin{equation*}\partial_{1} (\varphi)(\sigma_{0})=-\sum_{\substack{{\sigma\in \X_{0}} \\{ \textrm{ à distance } 1}}} f({{\sigma}})\:=\:f({\sigma_{0}})
\end{equation*}
et  $\partial_{1} (\varphi)(\sigma)=f({\sigma})$  pour tout sommet  $\sigma$ à distance $1$. Autrement dit, $\partial_{1}(\varphi)=f$. 

\medskip
Supposons  que $m\> 2$ et que la propriété est vraie aux rangs $\< m-1$. D'après le lemme \ref{Hurt}, on a l'inclusion $\Cc^{\U_{\sigma}}\subset \Cc^{\K_{m}}$ pour tout sommet $\s$ à distance $\< m-1$. Par conséquent,  on a $$\sum_{\substack{{\sigma\in \X_{0}} \\{ \textrm{ à distance } m}}} f({{\sigma}})\in  \Cc^{\K_{m}}.$$ Par injectivité de (\ref{Key}), on en déduit que $f({\sigma})\in \Cc^{\U_{e_\sigma}}$ pour tout sommet $\sigma$ à distance $m$. 
On  peut  alors  définir la fonction $\varphi\in \Ff_{1}(\X,\Cc)$ de support l'ensemble des arêtes  à distance $m$  de la fa\c con suivante~: pour tout sommet $\sigma$  à distance $m$, la fonction $\varphi$ prend  en $e_{\sigma}$ la valeur
$f({{\sigma}})$.
On vérifie  que $\partial_{1}(\varphi)\in  \Ff_{0}(\X,\Cc)$ 
est une fonction de support contenu dans la boule de rayon $m$ 
telle que, pour tout sommet $\sigma$ à distance $m$,
\begin{equation*}\partial_{1} (\varphi)(\sigma)=-(-1)^{m} f({{\sigma}}).
\end{equation*} La fonction $f+(-1)^m\partial_{1}(\varphi)$, de support inclus dans la boule de rayon $m-1$,  appartient au noyau de $\alpha$. Par hypothèse de récurrence, $f$ appartient à l'image de $\partial_{1}$.

\bigskip

Supposons que le 
complexe (\ref{FCS-2}) est exact. Soit $m\> 1$. 
Montrons que (\ref{Key}) est injective. 
Soit $(v_{\sigma})_{\sigma}$ une famille non nulle d'éléments de $\Cc$ indexée 
par les sommets à distance $m$ avec $v_{\sigma}\in\Cc^{\U_{\sigma}}$. 
Supposons que $$\underset{\substack{{\sigma\in \X_{0}} \\{ \textrm{ à distance 
      } m}}} \sum v_{\sigma}\in\Cc^{\K_{m}}.$$ 
Par la proposition 
\ref{PSHQ}, il existe une famille $(v_{\sigma})_{\sigma}$ d'éléments de $\Cc$ 
indexée par les sommets $\sigma$ à distance $\< m-1$ avec 
$v_{\sigma}\in\Cc^{\U_{\sigma}}$ et 
\begin{equation}\label{dance}\sum_{\substack{{\sigma\in \X_{0}} \\{ \textrm{ à 
          distance } m}}} v_{\sigma}+\sum_{\substack{{\sigma\in \X_{0}} \\{ 
        \textrm{ à distance } \< m-1}}}v_{\sigma}=0.\end{equation} 

On définit la fonction $\v \in\Ff_{0}(\X,\Cc)$ à support dans la boule de 
rayon $m$ en posant $\v(\sigma):=v_{\sigma}$ pour tout sommet à distance $\< 
m$. 
Par (\ref{dance}), la fonction $\v$ appartient au noyau de $\alpha$ et par exactitude de 
(\ref{FCS-2}), à l'image de $\partial _{1}$~: il existe une fonction $\w \in 
\Ff_{1 }(\X,\Cc)$ telle que 
$\partial_{1}(\w)=\v$. 
Soit $r\> 1$ le rayon minimal tel que le support de $\w$ est contenu dans la 
boule de rayon $r$ et soit $\tau$ une arête à distance $r$ dans ce support. 
L'arête $\tau$ possède un unique sommet à distance $r$ et il n'y a pas d'autre 
arête dans le support de $\w$ le contenant. On en déduit que 
$\partial_{1}(\w)$ prend une valeur non nulle en ce sommet, puis que $r=m$. 
La valeur de $\partial_{1}(\w)$ en chaque sommet ${\sigma}$ à distance $m$ 
est égale à $-(-1)^m \w(e_{\sigma})\in \Cc^{\U_{e_{\sigma}}}$. Mais elle est 
aussi égale à $v_{\sigma}$. D'où $v_{\sigma}\in\Cc^{\U_{e_{\sigma}}}$. 
\end{proof}

Comme rappelé au paragraphe \ref{Shou}, la proposition suivante est démontrée dans \cite{SSReso} et l'on en déduit le corollaire
\ref{TheProp}.

\begin{prop}
Le complexe (\ref{FCS-2}) est exact.
\end{prop}
\begin{coro}\label{TheProp}
Pour tout $m\> 1$, l'application (\ref{Key}) est un isomorphisme.
\end{coro}

\subsection{Platitude de $\Cc$}
Nous allons montrer la proposition \ref{Montana} suivante.  Notons que les arguments utilisés ici sont latents dans   \cite{Platitude}, où un résultat similaire est prouvé, mais seulement pour le module 
$ \Cc/\varpi \Cc$ sur l'algèbre de Hecke 
$\Hh$. Remarquons d'ailleurs
que pour ce module universel, l'analogue du  crucial corollaire \ref{TheProp} est démontré directement  dans   \cite{Platitude} sans utiliser 
\cite{SSReso}. Cette preuve directe est  adaptable au cas du module universel $\Cc$.

\begin{prop}
\begin{enumerate}
\item Le $\Hh$-module $\Cc$ est plat si et seulement si le $\H$-module fini $\C$ est plat, autrement dit si et seulement si $q=p$.
Dans ce cas, ce sont même des modules projectifs.
\item  Le $\Hh$-module $\Cc_1$ est projectif.
\end{enumerate}
\label{Montana}
\end{prop}

Le critère de platitude suivant est tiré de
\cite[chap. I, \S2, n$^\circ$5, Proposition 5]{Bki}. 
\begin{lemm} \label{Bki1}Soit~:
\begin{equation*}
0\to\E'\to\E\to\E''\to0
\end{equation*}
une suite exacte de $\Hh$-modules à gauche. 
On suppose que $\E''$ est plat. 
Alors $\E$ est plat si et seulement si $\E'$ est plat. 
\end{lemm}

\begin{proof}[Preuve de la proposition \ref{Montana}]
\begin{enumerate}
\item
Supposons que $\Cc$ est un $\Hh$-module plat. 
En utilisant le lemme \ref{Bki1} et l'exactitude du complexe (\ref{FCS-2}), la platitude du $\Hh$-module $\Ff_{0}(\X,\Cc)$ est  équivalente à celle de
$\Ff_{1}(\X,\Cc)$. Ce dernier étant une somme directe de copies de
$\Hh$, il est plat, et l'on en déduit que $\Ff_{0}(\X,\Cc)$ est plat.
Puisque $\Ff_{0}(\X,\Cc)$  est une somme directe de copies de $\Cc^{\K_{1}}$, la proposition  2 \cite[chap. I, \S2, n$^\circ$2]{Bki} dit que $\Cc^{\K_{1}}$ est un 
$\Hh$-module plat, ce qui, par la proposition \ref{LockedOut} implique que
$\C$ est un $\H$-module plat.

Pour montrer l'implication réciproque, rappelons que $\Cc^{\I_{1}}$ est un facteur direct de $\Cc^{\K_{1}}$  (lemme \ref{Wyoming2}).
Supposons que $\C$ est un $\H$-module plat, c'est-à-dire que 
$\Cc^{\K_{1}}$ est un 
$\Hh$-module plat par la proposition \ref{LockedOut}, et même projectif  par la remarque qui la suit.
Alors,  pour tout $m\> 1$ et tout sommet à distance $m$, le $\Hh$-module
$\Cc^{\U_{\sigma}}/\Cc^{\U_{e_{\sigma}}}$ est projectif, puisqu'il est isomorphe.
à $\Cc^{\K_{1}}/\Cc^{\I_{1}}$ par la preuve du lemme \ref{Einziger}.
Par la proposition \ref{TheProp}, le $\Hh$-module $\Cc^{\K_{m+1}}/\Cc^{\K_{m}}$ est projectif pour tout $m\> 1$. 
Ainsi, $\Cc$ est  projectif, comme somme directe de modules projectifs.
\item 
En appliquant la projection $\e_1:\Cc\rightarrow \Cc_1$ qui est un idempotent central de $\Hh$, on obtient un isomorphisme de $\Hh_1$-modules analogue à $(\ref{Key})$ en remplaçant $\Cc$ par $\Cc_1$.
Or on sait que $\C_1$ est un $\H$-module projectif. Les arguments du point précédent s'appliquent alors et l'on en déduit que $\Cc_1$ est un $\Hh_1$-module projectif.
\end{enumerate}
\end{proof}

\subsection{Remarques sur les présentations des représentations de ${\rm GL}_{2}(\F)$}
\def\Mo{\mathfrak M}
\subsubsection{\label{sabin}}

On rappelle que le  complexe
 \begin{equation}
\tag{\ref{FCS-2}}
0\to\Ff_{1}(\X,\Cc)
\xrightarrow{\partial_{1}}
\Ff_{0}(\X,\Cc)\xrightarrow{\alpha}\Cc\to0, 
\end{equation} est un complexe exact  de 
$\Hh$-modules. On note $\Gp^\circ$ le sous-groupe de $\Gp$ des éléments de déterminant inversible dans $\Oo$.
Le complexe ci-dessus n'est pas $\Gp$-équivariant, mais simplement $\Gp^\circ$  équivariant.
Comme dans \cite{SSReso}, on  le modifie   pour en faire un complexe de représentations de $\Hh$ et de $\Gp$ qui lui est isomorphe en tant  que complexe de représentations de  $\Hh$ et de $\Gp^\circ$, en considérant les espaces de chaînes orientées comme suit.  Une arête orientée est un couple de la forme $(\sigma, \sigma')$ où
$\{\sigma, \sigma'\}$ est une arête de $\X$. On note $\X_{(1)}$ leur ensemble. On définit alors
$\Ff_{(1)}(\X,\Cc)$  comme l'ensemble des fonctions $f$ à support fini dans l'ensemble  des arêtes orientées   telles que $f(\s, \s')=-f(\s',\s)\in \Cc^{\U_{\{\s, \s'\}}}$ pour toute arête $\{\s,\s'\}$. L'espace $\Ff_{(0)}(\X,\Cc)$ est identique à $\Ff_{0}(\X,\Cc)$ et l'application de transition $ \partial_{(1)} :\Ff_{(1)}(\X,\Cc)\rightarrow \Ff_{(0)}(\X,\Cc)$ associe à la fonction de support $\{(\s, \s'), (\s', \s)\}$ et de valeur $v$ en $(\s, \s')$ la fonction  $\s\mapsto -v$,  $\s'\mapsto v$.
Le complexe  de représentations de $\Gp$ et de $\Hh$ \begin{equation}
\label{new}
0\to\Ff_{(1)}(\X,\Cc)
\xrightarrow{\partial_{(1)}}
\Ff_{(0)}(\X,\Cc)\xrightarrow{\alpha}\Cc\to0, 
\end{equation}
est exact.

On se donne un $\Hh$-module à droite $\Mo$.
Désormais dans ce paragraphe, on suppose que $q=p$ de sorte que $\Cc$ est un
$\Hh$-module projectif par la proposition \ref{Montana}.  
La  suite exacte (\ref{new}) admet  donc
une section $\Hh$-équivariante et le complexe suivant de représentations de 
$\Gp$ est exact 
 \begin{equation}
\label{FCS-M}
0\to\Mo\otimes_{\Hh}\Ff_{(1)}(\X,\Cc)
\xrightarrow{}
\Mo\otimes_{\Hh}\Ff_{(0)}(\X,\Cc)\xrightarrow{}\Mo\otimes_{\Hh}\Cc\to0.
\end{equation}
On obtient ainsi une résolution pour toute représentation de $\Gp$ de la forme
$\Mo\otimes_{\Hh}\Cc$. 

\medskip

\begin{rema}
Le complexe (\ref{new}) admet une section $\Gp$-équivariante naturelle $\mathscr S: \Cc\rightarrow \Ff_{(0)}(\X,\Cc)$ définie par $\1_{\I_{1}}\rightarrow f_{0}$ où $f_{0}$ est la fonction de support $\ori$ et de valeur $\1_{\I_{1}}\in\Cc^{\K_{1}}$. Vérifions que cette section n'est pas $\Hh$-équivariante. 
On pose  $\omega:=\begin{pmatrix}0&1\cr \varpi&0\cr\end{pmatrix}$. Cet élément normalise $\I_{1}$ et est de longueur nulle dans $\Wp$. L'image de $\1_{\I_{1}}$ par $\t_{\omega}$ est la fonction caractéristique  $\1_{\I_{1}\omega}=\omega^{-1}.\1_{\I_{1}}$ qui est envoyée par $\mathscr S$ sur la fonction de support $\omega^{-1}\ori$ et de valeur  $\omega^{-1}\1_{\I_{1}}$. L'image par $\t_{\omega}$ de $f_{0}$ est, quant à elle, la fonction de support $\ori$ et de valeur $\t_{\omega}\in\Cc^{\K_1}$.

\end{rema}

\subsubsection{} Dans ce paragraphe, on suppose de plus
 que $\F=\mathbb Q_{p}$. Soit $\pi$ une représentation de $\Gp$ ayant un caractère central. D'après \cite{Inv}, elle est isomophe à la représentation
  $$\pi^{\I_{1}}\otimes_{\Hh}\Cc$$ de sorte que (\ref{FCS-M}) fournit une résolution pour la représentation $\pi$. 
Dans les termes de Breuil-Paskunas (\cite{Paskunas}, \cite{BP}), 
  la représentation $\pi$ est l'homologie en degré $0$ du diagramme donné par l'injection $\I_{1}Z$-équivariante
\begin{equation}\label{diagramme}\pi^{\I_{1}}\hookrightarrow \pi^{\I_{1}}\otimes_{\Hh} \Cc^{\K_{1}}.\end{equation} 

Dans le cas où $\pi$ est admissible, en particulier si $\pi$ est irréductible 
(\cite{BL}, \cite{B}), l'espace vectoriel $\pi^{\I_{1}}\otimes_{\Hh}\Cc^{\K_{1}}$ est de dimension  finie puisque $\Cc^{\K_{1}}$ est un $\Hh$-module de type fini (proposition \ref{DistanceUn}). On obtient ainsi une \emph{présentation standard} pour toute représentation admissible
de ${\rm GL}_{2}(\mathbb Q_{p})$ ayant un caractère central.
Ce résultat a été démontré pour la première fois par
 \cite{C}.  
On en trouve d'autres preuves dans \cite{BP}, \cite{VigCrit}, \cite{Hu}, \cite{GK2}.
   
\begin{rema}
\begin{enumerate}

\item L'espace $\pi^{\I_{1}}\otimes_{\Hh} \Cc^{\K_{1}}$ est la sous-$\K$-représentation de $\pi$ engendrée par $\pi^{\I_{1}}$.  C'est un sous-espace $\K_{1}$-invariant.
Cet espace ne coïncide pas en général avec l'espace $\pi^{\K_{1}}$ des $\K_{1}$-invariants de $\pi$. Par exemple, dans \cite[Théorème 8.6]{Columbia}, Breuil 
donne un exemple où $\pi$ est une représentation supersingulière, $\pi^{\K_{1}}$ n'est pas une représentation semi-simple de $\K$, et $\pi^{\I_{1}}\otimes_{\Hh} \Cc^{\K_{1}}$ est son socle.

\item 
Comme dans  \cite[\S3]{SSReso}, on peut construire le complexe de représentations de $\Gp$ suivant~:
\begin{equation}
\label{Compl}
0\to\Ff_{(1)}(\X,\pi)
\xrightarrow{\partial_{(1)}}
\Ff_{(0)}(\X,\pi)\xrightarrow{\alpha}\pi\to0 
\end{equation}
où 
$\pi$ est une représentation lisse de $\Gp$ à c\oe fficients dans un corps 
algébriquement clos $\mathfrak K$ de caractéristique quelconque. 
Les espaces $\Ff_{(0)}(\X,\pi)$ et $\Ff_{(1)}(\X,\pi)$, l'application de 
transition $\partial_{(1)}$ et le morphisme d'augmentation sont définis comme 
au début du paragraphe \ref{sabin} en remplaçant $\Cc$ par $\pi$.

Dans le cas où $\mathfrak K$ est le corps des nombres complexes et $\pi$ est 
une représentation de $\Gp$ engendrée par son espace des vecteurs 
$\I_1$-invariants, il est prouvé dans \cite{SSReso} que le complexe 
(\ref{Compl}) est exact, donc en particulier, $\pi$ est isomorphe 
à  l'homologie en degré $0$
du diagramme donné par l'injection $\I_{1}$-équivariante
$\pi^{\I_{1}}\hookrightarrow \pi^{\K_{1}}.$

Le point (1) de la présente remarque dit que le complexe (\ref{Compl}) n'est pas exact en général si $\mathfrak K$ est de caractéristique $p$.
\end{enumerate}
\end{rema}   

\section{Le cas de $\GL_3(\F)$}
\label{IAS}

Dans toute cette section, on suppose que $n=3$.

\subsection{Notations et préliminaires\label{twin}}
\def\Zp{Z}

On note $\Zp$ le centre de $\Gp$. On reprend les notations de l'exemple 
\ref{Sweetgrass}.  
On note $\Pp$ le sous-groupe de $\K$ image réciproque du sous-groupe 
parabolique standard de $\G$ de  facteur de Levi $\M={\rm GL}_2(k)\times{\rm GL}_1(k)$. 
C'est un sous-groupe parahorique de $\Gp$ de pro-$p$-radical égal au
sous-groupe $\Nn$ défini dans l'exemple \ref{Sweetgrass}.
On pose  $$\omega:=\begin{pmatrix}0&1&0\cr 0&0&1\cr \varpi &0&0\cr\end{pmatrix}.$$ Remarquons que $\omega^3=\varpi{ \rm Id}$ où ${ \rm Id}$ est la matrice identité de $\Gp$.
Distinguons et notons $\sigma_2$ la chambre standard de l'immeuble $\X$ :  les sommets $\sigma_0$, 
 $\omega\sigma_0$ et  $\omega^2\sigma_0$ forment une chambre dont le
 stabilisateur sous l'action de $\Gp$ est le sous-groupe engendré par
$\omega$ et le sous-groupe d'Iwahori $\I$.
On note $\s_1$ l'arête $\{\sigma_0, \omega\sigma_0\}$.  Elle a pour 
stabilisateur le sous-groupe de $\Gp$ engendré par $\Zp$ et $\Pp$. 

Les facettes $\s_0,\s_1$ et $\s_2$ sont telles que les pro-$p$-sous-groupes de 
$\Gp$ associés (voir \S \ref{Shou}) sont respectivement 
 $\U_{\s_0}=\K_1$ et $\U_{\s_2}=\I_1$ et~:
\begin{equation*}
\U_{\s_{1}}=\Nn=1+
\begin{pmatrix} 
\p &\p&\Oo\\ 
\p&\p&\Oo\\
\p&\p&\p
\end{pmatrix}.
\end{equation*}

Commençons par un résultat  relatif à la structure de l'immeuble de 
${\rm GL}_3$.  

Soit $m\>1$ et $\sigma$ un simplexe de $\X$. On dira que $\sigma$ est à distance $m$ si $m$ est le plus petit rayon d'une boule
de centre $\sigma_0$ contenant le simplexe $\sigma$.

Soit $\sigma$ 
une chambre de $\X$ à distance $m$. 
Elle possède un sommet $x$ à distance $m$, un sommet $z$ à distance $m-1$ et 
un sommet $y$ à distance $m-1$ ou $m$. 
On dira qu'elle est de type $(a)$ si $y$ est à distance $m$, de type $(b)$ s'il est à distance $m-1$.

\begin{lemm}
\label{room} L'ensemble des chambres de $\X$ à distance $m$ contenant l'arête $\{x,y\}$ est
\begin{itemize}
\item réduit à $\{\sigma\}$ si $y$ est à distance $m$ c'est-à-dire si $\sigma$ est de type $(a)$,
\item constitué de $\sigma$ et de  chambres de type $(a)$
si $y$ est à distance $m-1$, c'est-à-dire si $\sigma$ est de type $(b)$.

\end{itemize}
\end{lemm}

\begin{proof}
Quitte à conjuguer par un élément de ${\Gp}$, on peut supposer
que $\sigma$ appartient à l'appartement standard.
D'après \cite[Lemmes 2.3.4 et 2.3.5]{BO}, tous les voisins de $x$ à distance $m-1$  appartiennent à chacun des appartements contenant à la fois $\s_{0}$ et $x$ (ils sont donc dans l'appartement standard) et sont au nombre de $1$ ou $2$.   De plus, si $x$ a deux voisins à distance $m-1$, il forme avec eux une chambre de l'immeuble.

\medskip

Soit $\sigma'=\{x,y,t\}$ une autre chambre à distance $m$ contenant $\{x,y\}$. 
\begin{itemize}
\item Supposons que $y$ est à distance $m$ et que $\sigma\neq \sigma'$. Alors $t$ est nécessairement à distance $m-1$ et est un voisin de $x$ à distance $m-1$ dans l'appartement standard. 
Si  $x$ a un seul  voisin à distance $m-1$  alors $t=z$  et $\sigma=\sigma'$ ce qui est exclu. Donc $x$ possède deux voisins à distance $m-1$, ce sont $t$ et $z$. Mais alors  $\{x, t,z\}$ forme une autre chambre de l'appartement standard, ce qui contredit le fait que  $y$ et $t$ sont voisins.

\item  
Supposons que  $y$ est à distance $m-1$ et que $\sigma'$ est, comme $\sigma$, de type $(b)$ c'est-à-dire que $t$ est également à distance $m-1$.
Alors $y$, $z$ et $t$ sont des voisins à distance $m-1$ de $x$ ce qui 
conduit à la seule possible conclusion que $z=t$ et $\sigma'=\sigma$.
\end{itemize}
\end{proof}

Pour ${\rm GL}_3(\F)$, le complexe augmenté de $\Hh$-modules à gauche (\ref{Fizzy}), dont on rappelle qu'il est exact,  s'écrit~: 
\begin{equation}
\label{ComplexeGL3}
0\to\Ff_{2}(\X,\Cc)\xrightarrow{\partial_{2}}\Ff_{1}(\X,\Cc) 
\xrightarrow{\partial_{1}}\Ff_{0}(\X,\Cc)\xrightarrow{\alpha}\Cc\to0.
\end{equation}

\begin{prop}
\label{morningside}
Soit $f\in \Ff_{2}(\X,\Cc)$. On suppose que toute arête dans le support de $\partial_2(f)$ est à  distance $\< m$. Alors les chambres du support de $f$ sont à distance $\<m$.

\end{prop}

\begin{proof} Soit $m'\>1$ le plus petit entier tel que  le support de $f$ est inclus dans la boule (fermée de centre $\sigma_0$) de rayon $m'$. On suppose que $m'>m$.
 Soit $\sigma$ une chambre du support de $f$ à distance $m'$~: 
 notons $x$ un sommet de $\sigma$ à distance $m'$ et
 $z$ un sommet de $\sigma$ à distance $m'-1$. Le troisième sommet noté $y$ est
 à distance $m'$ ou $m'-1$ selon que $\sigma$ est de type (a) ou (b).

Il est  impossible que  la chambre $\sigma$  soit la  seule du support  de $f$
contenant l'arête $\{x,y\}$. En effet, si c'était le cas, 
  la valeur de $\partial_{2}(f)$ en $\{x,y\}$ serait 
$f(\s)$ ou $-f(\s)$, qui est non nulle, ce qui 
est impossible puisque $\{x,y\}$ n'est pas dans le support de $\partial_2(f)$. 
D'après le lemme \ref{room}, la chambre  $\sigma$ ne saurait donc être de type
(a).  On a prouvé que $y$ est à distance $m'-1$ et qu'il n'y pas de chambre de 
type (a) à distance $m'$ dans le support de $f$. Une chambre du support de $f$ 
contenant $\{x,y\}$ est donc de type (b) à distance $m'$ et 
$\sigma$ est la seule qui
    satisfait ces critères  par le lemme \ref{room}. On obtient  une contradiction, donc $m'\<m$.\end{proof}

\begin{lemm}\label{sky}
Pour tout $i\in\{0,1,2\}$, le $\Hh$-module à gauche,  
$\Ff_{i}(\X,\Cc)$
est une somme directe de copies de $\Cc^{\U_{\s_i}}$.
\end{lemm}

\begin{proof}
Dans l'immeuble de $\Gp=\GL_3(\F)$, les facettes 
de dimension $i$ sont toutes conjuguées à $\s_i$ sous l'action de $\Gp$.  
Par ailleurs, pour $g\in\Gp$, on a 
$\U_{g\cdot\s_i}=g\U_{\s_i}g^{-1}$  (\cite[\S1]{SSReso}) et 
$\Cc^{\U_{g\cdot\s_i}}=g\cdot\Cc^{\U_{\s_i}}$. 
Le résultat s'ensuit.
\end{proof}

\begin{lemm}Soit $\Aa$ un idéal à droite de $\Hh$. 
On a~: 
\begin{equation}
\label{Hamstring}
\Aa\Im(\partial_{2})=\Aa\Ff_1(\X,\Cc)\cap\Im(\partial_{2}).
\end{equation}
\end{lemm}

\begin{proof}
L'inclusion du membre de gauche dans celui de droite est 
immédiate.
Il s'agit donc de prouver que l'on a l'inclusion contraire, 
ce que l'on fait en raisonnant par l'absurde. 
Soit donc $f\in\Ff_{2}(\X,\Cc)$ une fonction non nulle 
telle que $\partial_{2}(f)$ appartient à $\Aa\Ff_{1}(\X,\Cc)$ 
mais n'appartient pas à $\Aa\Im(\partial_{2})$.
S'il existe une chambre $\s$ dans le support de $f$ telle que 
$f(\s)\in\Aa\Cc^{\U_{\s}}$, on note $f'$ l'élément de 
$\Ff_2(\X,\Cc)$ de support $\s$ et prenant 
la même valeur que $f$ en $\s$. 
On a donc $\partial_{2}(f')
\in\Aa\Im(\partial_{2})$ d'après le lemme \ref{Fox}
et, quitte à remplacer $f$ par $f-f'$, on peut  supposer que~: 
\begin{equation}
\label{Eq*}
f(\s)\notin\Aa\Cc^{\U_{\s}}
\end{equation}
pour toute chambre $\s$ du support de $f$. 
Nous allons aboutir à une contradiction en suivant un raisonnement analogue à 
celui de la preuve de la proposition \ref{morningside}.  

\medskip

Soit $m$ le plus petit entier $\>1$ tel que 
toutes les chambres dans le support de $f$ soient contenues 
dans la boule (fermée de centre $\s_{0}$) et de rayon $m$. 
Soit $\s$ une chambre à distance $m$ dans le support de $f$ : elle possède un sommet $x$ à 
distance $m$ et ses autres sommets $y$ et $z$ sont à distance 
$\<m$. Au moins l'un des deux, disons $z$, est à distance $m-1$.  
Remarquons qu'il est impossible que la chambre $\sigma$ soit la seule du support de $f$ contenant  l'arête $\{x,y\}$. En  effet, si c'était le cas,
regardons $\partial_{2}(f)$~: 
c'est une fonction dont la valeur en $\{x,y\}$ est 
$f(\s)$ ou $-f(\s)$, qui est non nulle.
Puisque l'on a supposé que $\partial_{2}(f)\in\Aa\Ff_{1}(X,\Cc)$, 
c'est que cette valeur appartient à 
$\Aa\Cc^{\U_{\{x,y\}}}$.
On a donc~:
\begin{equation*}
f(\s)\in\Aa\Cc^{\U_{\{x,y\}}}\cap\Cc^{\U_{\s}}=\Aa\Cc^{\U_{\s}}
\end{equation*}
car $\Cc^{\I_{1}}$ est un facteur direct de $\Cc^{\U_{\sigma_1}}$ comme 
$\Hh$-module d'après la remarque \ref{bloom}, et en utilisant le lemme 
\ref{sky}.  
Ceci est en contradiction avec (\ref{Eq*}), donc la chambre $\sigma$ n'est pas 
la seule du support de $f$ contenant l'arête $\{x,y\}$.

\medskip

D'après le  lemme \ref{room},  la chambre  $\sigma$ est donc  de type  (b), le
sommet $y$ est à distance $m-1$. Ce raisonnement étant valable pour toutes les
chambres du  support de $f$  à distance $m$,  on sait qu'elles sont  toutes de
type (b). 
Mais alors, appliquant à nouveau le lemme \ref{room}, la chambre $\sigma$
est l'unique du support  de $f$ contenant $\{x,y\}$ et l'on aboutit  à une contradiction.

On a ainsi prouvé l'égalité (\ref{Hamstring})\end{proof}

Notons respectivement $\Bb_1(0, m)\subset \Ff_1(\X,\Cc)$ et $\Bb_2(0, m)\subset \Ff_2(\X,\Cc)$
l'ensemble  des fonctions  de  support  inclus dans  l'ensemble  des arêtes  à
distance $\leq m$, respectivement des chambres à distance $\leq m$. 
Ce sont des sous-espaces stables sous l'action de $\Hh$. 

\begin{lemm}\label{notch}
Soit $\Aa$ un idéal à droite de $\Hh$.  
On a~: \begin{equation} \partial_2(\Bb_2(0, m))\cap \Aa \Bb_1(0, m)= 
  \Aa \partial_2(\Bb_2(0, m)). \end{equation} 
\end{lemm}

\begin{proof}Seule l'inclusion du membre de gauche dans celui de droite  mérite une vérification.
Soit $f$ une fonction de $\Bb_2(0, m)$ telle que $\partial _2(f)\in \Aa \Bb_1(0, m)\subset \Aa\Ff_1(\X,\Cc) $. D'après l'égalité (\ref{Hamstring}), 
l'élément $\partial _2(f)$ appartient à $\Aa\Im\partial_2$ ce qui signifie qu'il existe une fonction $g\in \Aa\Ff_2(\X,\Cc)$ telle que 
$\partial_2(f)=\partial_2(g)$. Mais $\partial_2$ est injective, donc 
 $f\in \Aa\Ff_2(\X,\Cc)\cap \Bb_2(0, m)= \Aa\Bb_2(0, m)$.
\end{proof}

 L'espace $\Ff_{0}(\X,\Cc)$ s'identifie comme représentation de $\Gp$ et comme $\Hh$-module à gauche à l'induite compacte
$\ind_{\K Z}^\Gp  \Cc^{\U_{\sigma_0}} $   de la représentation  $\Cc^{\U_{\sigma_0}}$ 
de $\K \Zp$ sur laquelle on fait agir $\Zp$ trivialement.
 D'après le corollaire \ref{Alex},  
 l'espace $\Cc^{\U_{\sigma_0}}$ s'identifie comme $\Hh$-module à gauche et comme représentation de $\K$ à $\Hh\otimes_\H\C$ et $\C$ en est un facteur direct comme $\H$-module d'après la proposition \ref{facteurdirect}.  
Le $\Hh$-module à gauche  $\Ff_{0}(\X,\Cc)$ s'identifie  donc   au produit tensoriel
 $\Hh\otimes_\H \ind_{\K \Zp}^\Gp \, \C $ et l'espace 
 $\ind_{\K \Zp}^\Gp \, \C$  en est un facteur direct comme $\H$-module.

\medskip

L'espace $\Ff_{1}(\X,\Cc)$ s'identifie comme  représentation de $\Gp$ et comme
$\Hh$-module à gauche à l'induite compacte 
$\ind_{\Pp    \Zp}^\Gp   \Cc^{\U_{\sigma_1}}    $    de   la    représentation
$\Cc^{\U_{\sigma_1}}$   de  $\Pp\Zp$   sur   laquelle  on   fait  agir   $\Zp$
trivialement.  L'action  de $\Pp$  sur $\Cc^{\U_{\sigma_1}}$   se factorise
par  le  quotient $\Pp/\Nn\simeq  \M$.  Par  la proposition  \ref{DuelAtDawn},
l'espace  $\Cc^{\U_{\sigma_1}}$  s'identifie comme  $\Hh$-module  à gauche  et
comme représentation de $\Pp$ 
à $\Hh\otimes_\H\C^{\N}$ et $\C^{\N}$ en est un facteur direct comme
$\H$-module.  
Le $\Hh$-module à gauche  $\Ff_{1}(\X,\Cc)$ s'identifie donc  au produit tensoriel
 $\Hh\otimes_\H \ind_{\Pp\Zp}^\Gp  \,\C^{\N}$ et l'espace 
 $\ind_{\Pp \Zp}^\Gp \, \C^{\N}$  en est un facteur direct comme $\H$-module.

\medskip
On considère l'application   de transition 
${\partial_1}$
 comme  une application $\ind_{\Pp Z}^\Gp  \Cc^{\U_{\sigma_1}} \rightarrow \ind_{\K \Zp}^\Gp  \Cc^{\U_{\sigma_0}} $ et l'on note  $\partial_1^0$ sa restriction à
 $ \ind_{\Pp\Zp}^\Gp  \C^{\N}$. D'après la définition de $\partial _1$, l'application $\partial_1^0$ est à valeurs dans
  $\ind_{\K \Zp}^\Gp  \C $.

 \begin{rema}\label{david}Puisque $\partial_1$ est $\Hh$-équivariante, elle s'identifie à l'application 
   $$\id_\Hh\otimes \partial _1^0: \Hh\otimes_\H \ind_{\Pp \Zp}^\Gp  \,\C^\N\longrightarrow \Hh\otimes_\H \ind_{\K \Zp}^\Gp  \,\C$$
 de sorte que $\Im \partial _1$ est isomorphe à  $\Hh\otimes _\H\Im \partial _1^0$ comme $\Hh$-module à gauche, car $\Hh$ est un $\H$-module libre d'après la proposition \ref{BrocoliRabe}. 
\end{rema}

\subsection{Filtration de $\Im \partial_1$\label{anya}}  Soit $m\>1$.
On rappelle les critères 
suivants, tirés du corollaire à la 
pro\-po\-sition 7 de \cite[chap. I, \S2, n$^\circ$6]{Bki}. 

\begin{lemm}
\label{Bki2}
Soit~:
$0\to\E'\to\E\to\E''\to0$
une suite exacte de $\Hh$-modules à gauche. 
\begin{enumerate}
\item 
On suppose que le $\Hh$-module $\E''$ est plat.
Alors, pour tout idéal à droite $\Aa$ de $\Hh$, on a~:
\begin{equation} 
\label{DownwardFacingDog}
\Aa\E'=\E'\cap\Aa\E.
\end{equation}
\item 
On suppose $\E$ est plat et que, pour tout idéal à droite 
$\Aa$ de $\Hh$, on a l'égalité (\ref{DownwardFacingDog}).
Alors le $\Hh$-module $\E''$ est plat.
\end{enumerate}
\end{lemm}

\begin{prop}
\label{Fierce}
Le $\Hh$-module à gauche
$\Bb_1(0, m)$ est plat si et seulement si  $\Cc^{\U_{\sigma_1}}$ est un $\Hh$-module à gauche plat, auquel cas
$\partial _1(\Bb_1(0, m))$ est également plat.

\end{prop}

\begin{proof} La première assertion provient du fait que le
$\Hh$-module à gauche $\Bb_1(0, m)$ s'identifie à une somme directe (finie) de copies de $\Cc^{\U_{\sigma_1}}$. Supposons maintenant que
le $\Hh$-module à gauche
$\Bb_1(0, m)$ est plat.
Le $\Hh$-module $\partial _1(\Bb_1(0, m))$ est isomorphe au quotient
de $\Bb_1(0, m)$ par l'intersection de $\Bb_1(0, m)$ avec l'image de $\partial _2$ par exactitude du complexe (\ref{ComplexeGL3}). D'après  la proposition \ref{morningside}, cette intersection est égale à $\partial_2(\Bb_2(0, m))$.   Les lemmes \ref{Bki2} (2) et \ref{notch} assurent alors la platitude du module quotient  $\partial _1(\Bb_1(0, m))$.
\end{proof}

Puisque le $\Hh$-module $\Im \partial_1$ est la limite inductive  des $\partial _1(\Bb_1(0, m))$, on obtient le corollaire suivant en appliquant \cite[Chapitre 1, \S2, n$^\circ$3, Proposition 2 (ii)]{Bki}.

\begin{coro}\label{cock}
Si $\Cc^{\U_{\sigma_1}}$ est un $\Hh$-module plat,  alors le $\Hh$-module à gauche $\Im \partial_1$ est plat.\end{coro}

\begin{coro} Si $\Cc^{\U_{\sigma_1}}$ est un $\Hh$-module plat,
 le $\Hh$-module à gauche $\Im\partial_1$ est un facteur direct de $\Ff_{0}(\X,\Cc)$.
\label{Fiery}
\end{coro}

\begin{proof}[Preuve du corollaire \ref{Fiery}]
Par la remarque \ref{david}, le $\Hh$-module $\Im \partial _1$ s'identifie au produit tensoriel $\Hh\otimes_\H\Im \partial _1^0$.   Le corollaire \ref{Fiery} sera prouvé lorsque l'on aura vérifié que le  $\H$-module à gauche $\Im \partial _1^0$ est un facteur direct de $ \ind_{\K Z}^\Gp  \,\C$. Ce résultat est donné par le lemme suivant.
 
 \begin{lemm}\label{Fairy} Si $\Cc^{\U_{\sigma_1}}$ est un $\Hh$-module plat, le  $\H$-module à gauche $\Im \partial _1^0$ est injectif.
 
 \end{lemm}

\begin{proof}   On appelle  $\Bb^0_1(0, m)$ l'intersection $  \ind_{\Pp\Zp}^\Gp  \C^{\N}\cap \Bb_1(0, m)$. Le  $\H$-module à gauche $\Im \partial _1^0$ est la limite inductive  des $\partial_1^0(\Bb^0_1(0, m))$. Puisque $\H$ est un anneau noetherien, le théorème
 de Bass-Papp \cite[(3.46)]{LMR} assure que le lemme sera démontré
 si l'on prouve que chaque $\H$-module  $\partial_1^0(\Bb^0_1(0, m))$ est injectif. Par auto-injectivité de $\H$, il suffit  de prouver que
 $\partial_1^0(\Bb^0_1(0, m))$  est un $\H$-module projectif, ou encore qu'il est plat  par le lemme \ref{lam}.
  On remarque alors que $\Hh\otimes _\H \partial_1^0(\Bb^0_1(0, m))$
  s'identifie à $\partial_1(\Bb_1(0, m))$ qui est un $\Hh$-module plat   d'après la proposition  \ref{Fierce}, sous l'hypothèse que $\Cc^{\U_{\sigma_1}}$ est un $\Hh$-module plat. On applique alors le lemme \ref{Fur} et 
  l'on obtient la platitude du $\H$-module à gauche $\partial_1^0(\Bb^0_1(0, m))$.\end{proof}
Le corollaire \ref{Fiery} est démontré.
\end{proof}

\subsection{Etude de la platitude des modules universels lorsque $q=p$.}

On suppose désormais que $q=p$.

\subsubsection{Le $\Hh$-module $\Cc$}

\begin{prop}\label{Novecento} 
Supposons que $q=p\neq 2$.
Le $\Hh$-module $\Cc$ n'est pas plat.
\end{prop}

\begin{proof}
 On a la suite exacte de $\Hh$-modules~:
\begin{equation}\label{Exacte1}
0 \to\Im\,\partial_{1}\to\Ff_{0}(\X,\Cc)\to\Cc\to 0. 
\end{equation}
Supposons 
 que $q=p\neq 2$. Alors $\Cc^{\U_{\sigma_1}}$ est un $\Hh$-module plat (voir l'exemple \ref{Sweetgrass}). 
Le corollaire \ref{cock} donne donc la platitude du $\Hh$-module 
$\Im\partial_{1}$.  Ainsi, d'après le lemme \ref{Bki1}, 
si $\Cc$ était un $\Hh$-module plat, alors  $\Ff_{0}(\X,\Cc)$ serait également un $\Hh$-module plat. D'après les résultats rassemblés dans l'exemple \ref{Yellowstone}, ce n'est pas le cas. Donc 
le  $\Hh$-module $\Cc$ n'est pas plat.
\end{proof}

\subsubsection{Le $\Hh$-module $\CcI$}

Rappelons que l'on a un complexe exact de représentations de $\HhI$
\begin{equation}
\label{ComplexeGL3I}
0\to\Ff_{2}(\X,\CcI)\xrightarrow{\partial_{2}}\Ff_{1}(\X,\CcI) 
\xrightarrow{\partial_{1}}\Ff_{0}(\X,\CcI)\to\CcI\to0
\end{equation}
obtenu en appliquant au système de c\oe fficient (\ref{ComplexeGL3})
la projection $\e_1: \Cc\rightarrow \CcI$  qui est un idempotent central de $\HhI$.  Par ce même procédé de projection, les résultats du paragraphe \ref{anya} sont préservés   si l'on remplace
$\Cc$ par $\CcI$, $\C$ par $\CP$, $\Hh$ par $\HhI$ et
$\H$ par $\HP$. 
D'après la remarque \ref{parme}, le $\Hh$-module $\CcI^{\U_{\sigma_1}}$ est 
plat que $q$ soit égal à $p$ ou non. On en déduit le résultat suivant~:

\begin{prop}
 Le $\HhI$-module à gauche $\partial_1(\Ff_{1}(\X,\CcI))$ est un facteur direct de $\Ff_{0}(\X,\CcI)$.

\end{prop}

Le  $\HhI$-module  $\Ff_{0}(\X,\CcI)$  est  une  somme directe  de  copies  de
$\CcI^{\U_{\sigma_0}}$  qui est  un $\HhI$-module  projectif  sous l'hypothèse
$q=p$ (Exemple  \ref{Yellowstone} (4)).  D'après la proposition  précédente et
par exactitude  du complexe (\ref{ComplexeGL3I}), le  $\HhI$-module $\CcI$ est
un facteur direct de $\Ff_{0}(\X,\CcI)$, et l'on a prouvé~:

\begin{coro}
Si $q=p$, le $\HhI$-module à gauche $\CcI$ est projectif.\\
En particulier, si $q=p=2$,  le $\Hh$-module à gauche $\Cc$ est projectif.
\end{coro}

\providecommand{\bysame}{\leavevmode ---\ }
\providecommand{\og}{``}
\providecommand{\fg}{''}
\providecommand{\smfandname}{\&}
\providecommand{\smfedsname}{\'eds.}
\providecommand{\smfedname}{\'ed.}
\providecommand{\smfmastersthesisname}{M\'emoire}
\providecommand{\smfphdthesisname}{Th\`ese}


\begin{thebibliography}{10}

\bibitem{BL} {\scshape L.~Barthel {\normalfont \smfandname} R.~Livné} 
  -- {\og Irreducible modular 
representations of ${\rm GL}\sb 2$ of a local field\fg}, 
\emph{Duke Math. Journal} {\bf 75}, no. 2, (1994) p.~261--292.

\bibitem{BO}{\scshape J.~Bella{\"i}che {\normalfont \smfandname} 
A.~Otwinowska} -- {\og Platitude du module universel pour 
${\rm GL}\sb 3$ en caractéristique non banale\fg}, 
\emph{Bull. Soc. math. France} {\bf 131} (4), 2003, p.507--525.

\bibitem{Benson}{\scshape D.~J.~Benson} -- \emph{Representations and 
    Cohomology I, Basic Representation Theory of Finite Groups and Associative 
    Algebras}, Cambridge Studies in Advanced Mathematics (1991). 

\bibitem{Borel} {\scshape A.~Borel} -- {\og Admissible representations of a 
semi-simple group over a local field with vectors fixed under an Iwahori 
subgroup\fg}, \emph{Inventiones Math.} {\bf 35}, (1976) p.233--259.

\bibitem{Bki} {\scshape N.~Bourbaki} -- \emph{Éléments de mathématiques. 
Fascicule XXVII. Algèbre commutative. Chapitre~1~: Modules plats. 
Chapitre~2~: Localisation}.  Hermann, Paris (1961). 

\bibitem{B}{\scshape C.~Breuil} -- {\og Sur quelques repr{\'e}sentations
modulaires et $p$-adiques de ${\rm GL}\sb 2(\bold Q\sb p)$ I\fg}, 
\emph{Compositio Mathematica.} {\bf 138}, no. 2 (2003) p.165--188. 

\bibitem{Columbia}{\scshape C.~Breuil} -- {\og Representations of Galois 
  and of GL2 in characteristic $p$\fg}, Cours à l'université de Columbia, 
  \url{http://www.ihes.fr/~breuil/}, (2004)

\bibitem{BP} {\scshape C.~Breuil {\normalfont \smfandname} V.~Paskunas} 
  -- {\og Towards a modulo $p$ Langlands correspondence for GL(2)\fg}.  
A paraître dans \emph{Memoirs of Amer. Math. Soc.}

\bibitem{Broussous} {\scshape P.~Broussous} -- 
{\og Acyclicity of Schneider and Stuhler's coefficient systems: 
another approah in the level $0$ case\fg}, 
\emph{J. Algebra} {\bf 279} (2004), p.~737--748. 

\bibitem{Brown}
 {\scshape K.~Brown} -- \emph{ Buildings}, Springer-Verlag, New York, (1989). 

\bibitem{CE}
{\scshape M.~Cabanes {\normalfont \smfandname} M.~Enguehard} -- 
\emph{Representation theory of finite reductive groups}, Cambridge
University Press, (2004).

\bibitem{Carter}
{\scshape R.~W.~Carter} -- \emph{Finite Groups of Lie Type.} 
Wiley Interscience, (1985).

\bibitem{CL} 
{\scshape R.~W.~Carter {\normalfont \smfandname} G.~Lusztig} -- 
{\og Modular representations of finite groups of Lie type\fg}, 
{\it Proc. London Math. Soc.} {\bf 32} (1976), p.~347--384.

\bibitem{Diamond} {\scshape F.~Diamond} --
\og A correspondence between representations of local Galois 
groups and Lie-type groups.\fg {\it L-functions and Galois representations / 
  edited by David Burns, Kevin Buzzard and Jan Nekovár.} Cambridge University 
Press, (2007). 

\bibitem{C} {\scshape P.~Colmez} --
\og Représentations de ${\rm GL}_2({\mathbb Q}_p)$ et $(\phi,\Gamma)$-modules 
\fg, \emph{Astérisque} {\bf 330} (2010), p.~281--509.

\bibitem{GK} {\scshape E.~Grosse-Klönne} --
\og On the universal module of $p$-adic spherical Hecke algebras \fg. 
\emph{Prépublication}, 
\url{www.math.hu-berlin.de/~zyska/Grosse-Kloenne/Preprints.html}.  

\bibitem{GK2} \bysame
\og $p$-torsion coefficient systems for ${\rm SL}_2({\mathbb Q}_p)$ and 
${\rm GL}_2({\mathbb Q}_p)$\fg, \emph{Prépublication}, 
\url{www.math.hu-berlin.de/~zyska/Grosse-Kloenne/Preprints.html}.  

\bibitem{Herzig} 
{\scshape F.~Herzig} -- {\og The weight in a Serre type conjecture 
for tame $n$-dimensional Galois representations\fg}, 
\emph{Duke Math. J.} {\bf 149} (2009), p.~37--116.

\bibitem{Hu} {\scshape Y.~Hu} -- \og
Diagrammes canoniques et représentations modulo $p$ de ${\rm GL}_2(F)$\fg 
\emph{Prépublication}, à paraître au Journal de l'Institut de Mathématiques de 
Jussieu.

\bibitem{Jantzen}
{\scshape J.~C.~Jantzen} -- \emph{Representations of algebraic groups}, 
Second edition.  Mathematical Surveys and Monographs, 107.  American 
Mathematical Society, Providence, RI, 2003.

\bibitem{Jeya}
{\scshape A.~V.~Jeyakumar} -- {\og Principal indecomposable 
representations for the group {${\rm SL}(2,q)$}\fg}, 
{\it J. Algebra} {\bf 30} (1974), p.~444--458.

\bibitem{Lam} 
{\scshape T.~Y.~Lam} -- \emph{A First course in noncommutative rings}, 
Springer {(1991)}.

\bibitem{LMR} 
\bysame, {\og Lectures on modules and rings. }
\emph{Graduate Texts in Mathematics}, {\bf 189}. Springer-Verlag,  (1999).

\bibitem{Lu} {\scshape G.~Lusztig} -- \fg Affine Hecke algebras and their 
  graded 
version,\og \emph{Journal of A.M.S}. {\bf{Vol. 2, No.3}}, (1989). 

\bibitem{OGL3}
{\scshape R.~Ollivier} -- {\og Modules simples en caractéristique p sur
   l'algèbre de Hecke du pro-p-Iwahori de ${\rm GL(3,F)}$\fg}.  \emph{Journal 
  of Algebra} {\bf 304} p.~1-38 (2006). 

\bibitem{Platitude}
\bysame, {\og Platitude du pro-$p$-module universel 
de ${\rm GL}_2(F)$ en caractéristique $p$\fg}, \emph{Compositio
Math.} {\bf 143} p.~703--720  (2007).

\bibitem{Inv}
\bysame , {\og Le foncteur des invariants sous l'action du 
pro-$p$-Iwahori de ${\rm GL}_2(\mathbb{Q}_p)$\fg},
\emph{J. für die reine und angewandte Mathematik} {\bf 635} p.~149--185 
(2009). 

\bibitem{OParab}
\bysame , {\og Parabolic induction and Hecke modules in 
characteristic $p$ for $p$-adic ${\GL_n}$\fg}, \emph{Prépublication} 
\url{http://www.math.uvsq.fr/~ollivier/} (2009). 

\bibitem{Paskunas}
{\scshape V.~Paskunas} -- {\og Coefficient systems and supersingular 
representations of ${\rm GL}_2(F)$\fg}, \emph{Mém. Soc. Math. Fr. (NS)} 
{\bf 99} (2004).

\bibitem{SS} 
{\scshape P.~Schneider {\normalfont \smfandname} U.~Stuhler} \blue{--}
{\og Representation theory and sheaves on the Bruhat-Tits building\fg}, 
{\it Publ. Math. IHES} {\bf 85} p.~97--191 (1997).

\bibitem{SSReso} 
{\scshape P.~Schneider {\normalfont \smfandname} U.~Stuhler} -- 
{\og Resolutions for smooth representations of the general linear 
group over a local field\fg}, \emph{J. Reine Angew. Math.} \textbf{436} 
p.~19--32 (1993).

\bibitem{SV} {\scshape P.~Schneider {\normalfont \smfandname} M.-F.~Vignéras} 
  -- \og A functor from smooth o-torsion representations to (phi, 
  Gamma)-modules\fg, \emph{Prépublication}, 
  \url{www.math.uni-muenster.de/u/schneider/} (2009).

\bibitem{Serre} 
{\scshape J.-\P.~Serre} -- {\og Linear representations of finite groups\fg}, 
Graduate Texts in Math. 42, Springer.


\bibitem{VigBook} 
{\scshape M.-F.~Vignéras} --
{\og Repr{\'e}sentations $l$-modulaires d'un groupe r{\'e}ductif $p$-adique
avec $l\neq p$.\fg} 
Progress in Mathematics, 137. Birkh{\"a}user Boston, Inc., Boston, MA, (1996).

\bibitem{Selecta}\bysame --
{\og Induced {$R$}-representations of {$p$}-adic reductive groups\fg}, 
{\it Selecta Math. (N.S.)} {\bf 4} no.~4, p.~549--623 (1998). 

\bibitem{VigGL2} \bysame
\og Representations modulo $p$ of the $p$-adic group ${\rm GL}(2,F)$\fg, 
\emph{Compositio Math.} {\bf 140} p.~333--358, (2004). 

\bibitem{VigProp} \bysame
\og Pro-$p$-Iwahori Hecke ring and supersingular 
$\overline{\mathbb F}_{p}$-representations\fg, \emph{Mathematische Annalen} 
{\bf 331} p.~523--556. Erratum volume 
{\bf 333}, no. 3, p.~699--701, (2005). 

\bibitem{VigCrit} \bysame \og A criterion for integral structures and coefficient systems on the
tree of $\rm{PGL}(2, F )$ \fg, \emph{Pure and Applied Mathematics Quarterly} 
{\bf 4} p.~1291--1316 (2008). 

\end{thebibliography}
\end{document}